%% file: bandit_diffusion.tex
\def\useplain{1} 

\ifx \useplain\undefined
\documentclass[mnsc]{style/informs3aa}
\usepackage{style_set/informs_set}
\input{title_info/informs_title}

\else
\documentclass{article}
\usepackage[margin=1.5in]{geometry}
\usepackage{style_set/plain_set}
\input{title_info/plain_title}

\fi
\usepackage{style/eie_pckg}
\usepackage{style/macro}

\graphicspath{{figures/}}

\usepackage{refcount}

\def\rev#1{{\color{black} #1}}
\usepackage{ulem}

\begin{document}
{\tiny }
\ifx \useplain\undefined
\else
\maketitle
\fi

\def\atxt{
We use the lens of weak signal asymptotics to study a class of sequentially randomized experiments, including those that arise in solving multi-armed bandit problems. In an experiment with $n$ time steps, we let  the mean reward gaps between actions scale to the order $1/\sqrt{n}$ so as to preserve the difficulty of the learning task as $n$ grows.  In this regime, we show that the sample paths of a class of sequentially randomized experiments---adapted to this scaling regime and with arm selection probabilities that vary continuously with state---converge weakly to a diffusion limit, given as the solution to a  stochastic differential equation. The diffusion limit enables us to derive refined, instance-specific characterization of stochastic dynamics, and to obtain several insights on the regret and belief evolution of a number of sequential experiments including Thompson sampling (but not UCB, which does not satisfy our continuity assumption). We show that all sequential experiments whose randomization probabilities have a Lipschitz-continuous dependence on the observed data suffer from sub-optimal regret performance when the reward gaps are relatively large. Conversely, we find that a version of Thompson sampling with an asymptotically uninformative prior variance achieves near-optimal instance-specific regret scaling, including with large reward gaps, but these good regret properties come at the cost of highly unstable posterior beliefs.
	
	\noindent
	\textit{Keywords}: Diffusion approximation, Multi-armed bandit, Thompson sampling.}

\ifx \useplain\undefined
\ABSTRACT{\atxt}
\else
\begin{abstract}
	\atxt
\end{abstract}
\fi


\ifx \useplain\undefined
\maketitle
\fi


\section{Introduction}

Sequential experiments,\blfootnote{\hspace{-6.3mm}
Version: June 2023. An earlier draft of this paper was  circulated under the title {\it Diffusion Asymptotics for Sequential Experiments}. Xu Kuang published under a different full name in earlier versions of this manuscript. Please use {\it X.~Kuang and S.~Wager} when citing this paper.  \\
We are grateful to the referees and editors at Management Science for their detailed comments. 
We are also grateful for valuable feedback and suggestions from Susan Athey, Lin Fan, Peter Glynn, David Goldberg,
Keisuke Hirano, Michael Harrison, David Hirshberg, Max Kasy, Emilie Kaufmann, Tor Lattimore, Neil Walton,
Jack Porter, Daniel Russo, David Simchi-Levi, as well as seminar participants at a number of venues.}
pioneered by \citet{wald1947sequential} and \citet{robbins1952some}, involve collecting
data over time using a design that adapts to past experience.
Relative to classical randomized trials, they
can effectively concentrate power on studying the most promising alternatives and save on
costs by helping us avoid repeatedly taking sub-optimal actions.
Such experiments have now become widely adopted for automated decision making;
for example, \citet{hill2017efficient} show how sequential experiments can be used to optimize the layout and
content of a website, while \citet{ferreira2018online} discuss applications to pricing and online revenue management. 

Most existing research in sequential experimentation has focused on proving robust,  worst-case guarantees. This has been motivated in a large part by automated decision making applications,  such as website optimization and dynamic pricing, where it is common to run a large number of experiments 
in parallel while searching over a decision space.  In these settings, it is important to use robust methods that perform reasonably well across all problem instances, so that they do not  require
frequent human oversight; see  \citet{bubeck2012regret} for a review and discussion.

More recently, however, there has been growing interest in settings where an agent performs a small number of high-stake sequential experiments: \citet{athey2021increasing} discuss the use of sequential experiments for learning better
public health interventions, \citet{caria2020adaptive} deploy them for targeting job search assistance for
refugees, and \citet{kasy2020adaptive} consider whom to test for an infectious disease in a setting where testing
capacity is limited. In these applications where humans are likely to be closely monitoring and fine-tuning a small number of adaptive
experiments, the relatively crude performance guarantees provided by worst-case analysis fall short in guiding the decision maker towards choosing the best  policies. This gives rise to the need to develop a more refined, instance-specific understanding
of the stochastic behavior of adaptive experiments. In particular, it would be of interest to understand:
\begin{itemize}
\item 
How to quantify the performance of an experiment in an instance-specific manner, that is, as a function of features of the environment in which the experiment is being conducted, such as gaps in mean rewards among different actions.
\item Beyond just the mean, how to characterize the distributional properties of key performance measures of an adaptive experiment.
\item How the behavior of adaptive experiments evolves over time, and what sample paths of actions taken
by adaptive experiments look like.
\end{itemize}
Available worst-case-focused formal results, however, do not provide sharp answers to these questions that would
apply broadly to popular algorithms for adaptive experimentation, including the ones used in the studies cited above.
A central difficulty here is simply that sequential experiments induce intricate dependence patterns in the data that
make sharp finite-sample analysis of their behavior exceedingly delicate.

In this paper, we introduce a new approach to studying sequential experiments based on \rev{weak signal asymptotics}.
Specifically, we consider the behavior of adaptive experiments in a sequence of problems where, as the number of time steps
$n$ grows to infinity, the  gap between the mean rewards of different actions decays as $1/\sqrt{n}$.
\rev{We use this asymptotic regime to study the sequentially randomized Markov experiments,}
a general class of sequential experiments that includes several popular algorithms such as Thompson
sampling \citep{thompson1933likelihood,russo2018tutorial}. 
\rev{Our first result is that, under regularity conditions including continuity of how the arm selection
probabilities vary with certain state variables, the sample paths of appropriately scaled
sequentially randomized Markov experiments}
converge weakly to a diffusion limit characterized as the solution to a stochastic differential equation (SDE).
We then show that this diffusion limit enables us to derive practical, instance-specific and
distributional insights about the behavior of sequential experiments. 
The $ 1/\sqrt{n} $ arm  gap scaling considered here is an important one. It can be thought of as a ``moderate data'' regime of sequential experimentation, where the problem's inherent ``signal strength,'' expressed in terms of the mean reward gaps between the various actions, is not sufficiently stronge relative to sample size so as to make identifying the optimal action asymptotically trivial, but strong enough so that using a well designed policy is crucial as sub-optimal policies can lead to large regret. 

In using diffusion approximation to obtain insights about the behavior of large, dynamic processes, our approach
builds on a long tradition in operations research, especially in the context of queuing networks \citep{iglehart1970multiple,
harrison1981reflected,reiman1984open,harrison1988brownian,glynn1990diffusion,kelly1993dynamic,gamarnik2006validity}.
A key insight from this line of work is that, by focusing on a heavy-traffic limit where server utilization approaches full capacity and wait times diverge,
the behavior of a queuing network can be well approximated by a Brownian limit where, as argued by \citet{kelly1993dynamic},
``important features of good control policies are displayed in sharpest relief.'' Likewise, in our setting, we find that focusing on
a scaling regime with small gaps in the mean rewards between the arms and long time horizons enables us to capture key aspects of sequential experiments in terms
of a tractable diffusion approximation.

\subsection{Overview of Main Results}

Throughout this paper, we focus on the following $K$-armed setting,
also known as a stochastic $K$-armed bandit. There is a sequence of decision points $i = 1, \, 2, \, \ldots, \, n$
at which an agent chooses which action $A_i \in \cb{1, \, \ldots, \, K}$ to take and then observes a
reward $Y_i \in \RR$. Here, $Y_i$ is assumed to be drawn from a distribution $P_{A_i}$, where the
action $A_i$ may depend on past observations, and is conditionally independent from all other aspects of
the system given the realization of $A_i$. Our goal is to use diffusion approximation to characterize the
behavior of $K$-armed bandits in terms of a number of metrics, including regret
\begin{equation}
	\label{eq:regret}
	R^n = n \sup_{1 \leq k \leq K}\cb{\mu_k} - \EE{\sum_{i = 1}^n \mu_{A_i}}, \ \ \ \ \mu_k = \EE[P_k]{Y},
\end{equation}
i.e., the shortfall in rewards incurred by the bandit algorithm relative to always taking the best action. Unless otherwise specified, we will use ``arm gap'' as a short-hand for the gap between the mean rewards associated with the two best arms. 

The first part of our paper, Section \ref{sec:main}, \rev{introduces a weak signal scaling regime for the reward distributions, where arm gaps become small as the sample size of the experiment increases. We then establish that for a broad class
of sequential experiments, which we refer to as sequentially randomized Markov
experiments, the system dynamics converge to a diffusion limit under this weak signal asymptotic regime.}  This result applies to a wide variety of adaptive experimentation rules arising in statistics, machine learning and
behavioral economics \rev{including Thompson sampling
(see Examples \ref{ex:tempered_greedy}--\ref{ex:exploration} given below).}
We also find that this diffusion limit can be characterized using random time-change
applied to a Brownian motion, which enables much of our subsequent theoretical analysis.
\rev{Throughout our analysis, we assume that the arm selection probabilities in the 
sequentially randomized Markov experiments we consider vary continuously with relevant state variables.}
We provide convergence results for both when the limiting sampling function is Lipschitz continuous,
and when it is continuous but non-Lipschitz; \rev{however, our results do not apply to bandit
algorithms with a discontinuous deterministic component, such as UCB or $\varepsilon$-greedy}.

In the second part of the paper, we use our diffusion limit to obtain insights on the regret performance of a host of practical adaptive experiments. First, we turn our attention to the class of sequentially randomized experiments that admit an asymptotically Lipschitz limiting sampling function, which includes among other examples  the popular Gaussian Thompson sampling algorithm with an asymptotically proper prior. These  experiments tend to have favorable stability properties, making them less likely to overreact to the observed data and their dynamics more amenable to analysis. Surprisingly, we prove in Section \ref{sec:super_diffusive} that {any} experiment with an asymptotically Lipschitz sampling function necessarily has poor regret performance: in the diffusion limit, the scaled regret of any such experiment remains bounded from below as the arm gap gets large, whereas the
corresponding minimax risk goes to zero.
Our result thus shows that the benefits of stability of Lipschitz sampling come at the expense of regret performance, as the algorithm cannot not abandon a bad action quickly enough in the early stages of the experiment. 

As a counterpoint to the first finding, our next result demonstrates that, in a two-armed setting, we can achieve nearly instance-optimal regret by adopting a non-Lipschitz sampling function. We examine an ``undersmoothed'' version of the Thompson sampling algorithm where the agent uses an asymptotically uninformative prior distribution. 
We prove in Section \ref{sec:undersmoothedTS} that the regret of undersmoothed Thompson sampling not only vanishes as the arm gap grows, but does so at a rate that is nearly optimal in comparison to known instance-specific lower bounds \citep{mannor2004sample}, a sharp contrast to the negative result on Lipschitz sampling described earlier.
This is, to the best of our knowledge, the first result to highlight good instance-dependent
regret properties for Thompson sampling in the moderate signal-to-noise regime (i.e., with
$1/\sqrt{n}$-scale arm gaps).

Finally, we leverage the diffusion limit to obtain insights on distributional properties of sample paths in an adaptive experiment. We prove in Section \ref{sec:instability} that, for a family of one-armed problems, experiments with a ``sensitive'' sampling function---ones that respond quickly to observed data---will be highly unstable, where the sampling probability of choosing the optimal arm can swing wildly between very large to the very small, and that this can happen even with very large arm gaps. This result covers undersmoothed Thompson sampling, in which case our analysis implies large swings in its posterior beliefs. This finding complements our regret analysis to render a more complete and nuanced picture, suggesting that an agent may have to face a fundamental trade-off between achieving optimal regret and maintaining a relatively stable sample path.

On the methodological front, our work introduces new tools to the study of   adaptive experiments. We show weak convergence to the diffusion limit for sequentially randomized Markov experiments  using the martingale framework developed by  \citet{stroock2007multidimensional}, which hinges on showing that an appropriately scaled generator of the discrete-time Markov process converges to the infinitesimal generator of the diffusion process.  Our analysis of the regret profiles of Thompson sampling under the weak signal scaling relies on novel proof arguments that heavily exploit properties of Brownian motion, such as the Law of Iterated Logarithm.

\subsection{Related Work}

The multi-armed bandit problem is a popular framework for studying sequential experimentation; see
\citet{bubeck2012regret} for a broad discussion focused on bounds on the regret \eqref{eq:regret}.
An early landmark result in this setting is due to \citet{lai1985asymptotically}, who showed that given any
fixed set of arms \smash{$\cb{P_k}_{k = 1}^K$}, a well-designed sequential algorithm can achieve regret that scales logarithmically
with the number $n$ of time steps considered, i.e., \smash{$R^n = \oo_P(\log(n))$}. Meanwhile, given any fixed time horizon $n$, it
is possible to choose probability distributions \smash{$\cb{P_k}_{k = 1}^K$} such that the expected regret $\EE{R^n}$
of any sequential algorithm is lower-bounded to order \smash{$\sqrt{Kn}$} \citep{auer2002nonstochastic}.
It is worth-noting that the problem instance that achieves the $\sqrt{Kn}$ regret lower bound in \citet{auer2002nonstochastic}
involves the same mean reward scaling as we use for our diffusion limit, suggesting that the \rev{weak signal} scaling
proposed here captures the behavior of the most challenging (and thus potentially most interesting) sub-family of learning tasks.

Thompson sampling \citep{thompson1933likelihood, russo2018tutorial} has gained considerable popularity in recent years thanks to its simplicity and impressive empirical performance \citep{chapelle2011empirical}.  Regret bounds for Thompson sampling have been established in the frequentist \citep{agrawal2017near, kaufmann2012thompson} and Bayesian \citep{bubeck2014prior, lattimore2019information, russo2016information} settings; the setup here belongs to the first category. None of the existing {instance-dependent} regret bounds, however, appear to have sufficient precision to yield meaningful characterization in our regime. For example, the  instance-dependent  upper bound in \citet{agrawal2017near} contains a constant to the order of $1/\Delta^4$, where $\Delta$ is the gap in mean reward between optimal and sub-optimal arms. This would thus lead to a trivial bound of $\mathcal{O}(n^2)$ in our regime, where mean rewards scale as $1/\sqrt{n}$. Furthermore, many of the existing bounds also  require delicate assumptions on the reward distributions (e.g., bounded support). In contrast, the \rev{weak signal} asymptotics adopted in this paper are universal in the sense that they automatically allow us to obtain  approximations for a much wider range of reward distributions, requiring only a bounded fourth moment.

Our choice of the {weak signal} scaling and the ensuing {diffusion} limit are motivated by insights from both
queueing theory and statistics. The scaling   plays a prominent role in heavy-traffic diffusion
approximation in queueing networks \citep{gamarnik2006validity,harrison1981reflected,reiman1984open};
in particular, \citet{harrison1988brownian} uses diffusion approximation to study a dynamic control problem.
Here, one considers a sequence of queueing systems in which the excessive service capacity, defined as
the difference between arrival rate and service capacity, decays as $1/\sqrt{T}$, where $T$ is the time horizon.
Under this asymptotic regime, it is shown that suitably scaled queue-length and workload processes converge
to reflected Brownian motion. Like in our problem, the diffusion regime here is helpful because it captures the most
challenging problem instances, where the system is at once stable and exhibiting non-trivial performance variability.
See also \citet{glynn1990diffusion} for an excellent survey for the use of diffusion approximation  in operations research.

The \rev{weak signal asymptotics}  are further inspired by  a recurring insight from statistics that,
in order for asymptotic analysis to yield a normal limit that can be used for finite-sample insight,
we need to appropriately down-scale the arm gap as the sample size gets large. One concrete
example of this phenomenon arises when we seek to learn optimal decision rules from (quasi-)experimental
data. Here, in general, optimal behavior involves regret that decays as $1/\sqrt{n}$ with the sample size
\citep{athey2020policy,hirano2009asymptotics,kitagawa2018should}; however, this worst-case regret is only achieved if we let
effect sizes decay as $1/\sqrt{n}$. For any fixed sampling design, it's possible to achieve faster than
$1/\sqrt{n}$ rates asymptotically \citep{luedtke2017faster}.

Our work is broadly related to work on adaptive learning and experimentation in the statistics and operations research literatures. The seminal work of  \citet{chernoff1959sequential} studies how to adaptively choose from a menu of available experiments for hypothesis testing, and the main result consists of a dynamic experiment that achieves the asymptotic optimal rate of convergence to the true hypothesis as the number of samples grows. A number of papers expand upon this framework, including \cite{naghshvar2013active} who drive a range of refinements and extensions on the minimax bounds, \cite{russo2020simple} who adapts the algorithm for best-arm identification, and \cite{massoulie2018capacity} who consider simultaneous experiments run in a resource-limited setting. The asymptotic regime considered in this line of work is different from ours: They focus on  the large-deviation regime where the problem instance stays fixed while the number of samples grows to infinity, whereas we consider the \rev{weak signal} regime where the arm gap decreases along with the sample size. Another distinction is that this literature tends to treat what has become known as the post-experiment simple regret, where the agent's cost is determined by any error incurred in declaring the incorrect hypothesis after the completion of the experiment, rather than the in-experiment regret accrued due to using sub-optimal arms during the experiment. 
	
There are some existing results where diffusion-based analysis has been used for action selection and optimal stopping in sequential experiments. Some of the earlier studies assume randomization is fixed over the horizon \citep{siegmund1985sequential}. In contrast, in our multi-armed bandit setting the probabilities  in the randomization depend on the history which creates a qualitatively different limit object.   \citet{chick2009economic} consider optimal stopping in ranking and selection where the agent runs a sequence of experiments to determine which action carries the maximum expected payoff. They then apply diffusion approximation to the resulting Bellman equation and derive tractable structural results for simple cases. \cite{araman2019diffusion} study an adaptive optimal stopping problem with a binary hypothesis and use the diffusion limit to derive insights on how and when the agent should act, as well as the optimal choice of experiments.  \cite{wang2020adaptive} consider a  model similar to \cite{araman2019diffusion}  in the context of adaptive clinical trials, with the additional feature that the experiments can have different costs. \cite{harrison2015investment} use diffusion approximation to study an optimal stopping problem with a binary hypothesis when a firm decides between making an investment versus taking further learning actions to gather more information. The approaches taken in these papers are however different from ours. First, they tend to focus on discrete (and often binary) hypothesis testing whereas we study multi-armed experiments with continuously-valued mean rewards. 
Second, the optimal stopping formulation is centered around the expected cost incurred by a one-shot, post-experiment action choice, in contrast to the in-experiment regret analyzed here. Third, most of these papers focus on solving exactly, or approximately, the Bayesian Bellman-optimal decision rule based on dynamic programming, whereas we  are more interested in characterizing the dynamics of a range of heuristic policies, such as Thompson sampling, which are used in practice when the Bayesian priors of the model parameters are unknown or the agent operates under a frequentist, instance-specific performance mandate. Lastly, in this literature, diffusion approximation is sometimes evoked directly  or introduced as a  modeling assumption without a proof of convergence \citep[e.g.,][]{chick2009economic, harrison2015investment}.

Finally, after we posted a first version of this paper, a number of authors have disseminated work that
intersects with our main result. In the special case of Thompson sampling, \citet{fan2021diffusion} derive
the diffusion limit given in Theorems \ref{theo:SDE} and \ref{theo:timechange} using a different proof
technique from us: As discussed further in Appendix \ref{sec:proofs}, we derive the diffusion limit by studying
the asymptotic behavior of the Markov chain associated with our sequential experiments following
\citet{stroock2007multidimensional}, whereas \citet{fan2021diffusion} develop a direct argument based on
the continuous mapping theorem. Meanwhile, \citet{hirano2021asymptotic} extend the ``limits of experiments'' analysis
pioneered by \citet{lecam1972limits} to discrete-time, batched sequential experiments. The limits of experiments paradigm
also involves a $1/\sqrt{n}$-scale local parametrization, and their results can be applied
to sequentially randomized Markov experiments whose the randomization probabilities only
change at a finite set of pre-specified times. Both \citet{fan2021diffusion} and \citet{hirano2021asymptotic}
are based on research developed independently from ours.

\section{Asymptotics for K-Armed Sequential Experiments}
\label{sec:main}

As discussed above, the first goal of this paper is to establish a diffusion limit for a class
of sequential experiments.
To this end, we first introduce a broad class of sequential experimentation
schemes in Section \ref{sec:SRME}, which we refer to as sequentially randomized Markov
experiments. We describe a \rev{weak signal} scaling for sequential experiments in Section \ref{sec:scaling};
then, in Section \ref{sec:limit}, we establish conditions under which---in this limit---sample paths of sequentially
randomized Markov experiments converge weakly to the solution of a stochastic differential equation.

Throughout our analysis, we work within the following $K$-armed model. This model captures a number of interesting
problems and is widely used in the literature \citep[e.g.,][]{bubeck2012regret,lai1985asymptotically}. However, we note
that it does rule out some phenomena that may arise in applications; for example, we do not allow for distribution shift
in the reward distribution of a given arm over time, and we do not allow for long-term consequences of actions, i.e.,
an action taken in period $i$ cannot affect arm-specific period-$i'$ reward distributions for any $i' > i$. Extending
our asymptotic analysis to allow for distribution shift or long-term effects would be of considerable interest, but
we do not pursue this line of investigation in the present paper.

\begin{defi}[Stochastic $K$-Armed Bandit]
\label{def:bandit}
A stochastic $K$-armed bandit is characterized by time horizon $n$ and a set of $K$ reward distributions $P_k$ for
$k = 1, \, \ldots, \, K$. At each decision points $i = 1, \, 2, \, \ldots, \, n$, an agent chooses which action
$A_i \in \cb{1, \, \ldots, \, K}$ to take and then observes a reward $Y_i \in \RR$. The action $A_i$ is a random
variable that is measurable with respect to the observed history \smash{$\cb{A_{i'}, \, Y_{i'}}_{i'=1}^{i-1}$}.
Then, conditionally on the chosen action $A_i$, the reward $Y_i$ is drawn from the distribution $P_{A_i}$,
independently from the observed history.
\end{defi}

\subsection{Sequentially Randomized Markov Experiments}
\label{sec:SRME}

In the interest of generality, we state our main results in the context of sequentially randomized experiments
whose sampling probabilities depend on past observations only through the state variables
\begin{equation}
\label{eq:QS_raw}
Q_{k,i} = \sum_{j = 1}^i \mathbf{1}\p{\cb{A_j = k}}, \quad  S_{k,i} = \sum_{j = 1}^i \mathbf{1} \p{\cb{A_j = k}} Y_j,
\end{equation}
where $Q_{k,i}$ counts the cumulative number of times arm $k$ has been chosen by the time
we collect the $i$-th sample, and $S_{k,i}$ measures its cumulative reward. When useful,
we use the convention $Q_{k,0} = S_{k,0} = 0$.
Working with this class of algorithms,
which we refer to as sequentially randomized Markov experiments,
enables us to state results that cover many popular ways of running sequential experiments without
needing to derive specialized analyses for each of them.

\begin{defi}[Sequentially Randomized Markov Experiment]
\label{def:SRME}
A $K$-armed sequentially randomized Markov experiment chooses the $i$-th action $A_i$
by taking a draw from a distribution
\begin{equation}
\label{eq:randomized}
A_i \cond \cb{A_1, \, Y_1, \, \ldots, \, A_{i-1}, \, Y_{i-1}} \sim \operatorname{Multinomial}(\pi_i),
\end{equation}
where the sampling probabilities are computed using a measurable sampling function $\psi$,
\begin{equation}
\label{eq:psi}
\psi :  \R^K \times \R^K \rightarrow \mathbf{\Delta}^K, \ \ \ \pi_i = \psi\p{Q_{i-1}, S_{i-1}  },
\end{equation} 
where $(Q_{i-1}) =  (Q_{k,i-1})_{k=1, \ldots, K}, S_{i-1} = (S_{k,i-1})_{k=1, \ldots, K} $, and $\mathbf{\Delta}^K$ is the $K$-dimensional unit simplex.

\end{defi}

\begin{rema}[Capturing Time Dependence] It is useful to note that the family of experiments described in Definition \ref{def:SRME}  contains those algorithms whose sampling probabilities may depend on the time period, $ i $. This is captured through the fact that $ \sum_{k} Q_{i,k} =i $, that is, the time period $ i $ can be simply read off by calculating the $ L_1 $ norm of the vector $ Q_i $.
\end{rema}

We now examine several examples of popular algorithms that fit under the sequentially randomized Markov experiment framework. 

\begin{exam}
\label{ex:tempered_greedy}
A greedy agent may be tempted to always pull the arm with the highest apparent mean, $S_{i,k}/Q_{i,k}$;
however, this strategy may fail to experiment enough and prematurely discard good arms due to early
unlucky draws. A tempered greedy algorithm instead chooses arm $ k $ with probability: 
\begin{equation}
\label{eq:tempered_greedy}
\pi_{i,k} = \exp\sqb{\alpha \frac{S_{k,i}}{Q_{k,i} + c}} \, \bigg/ \, \sum_{l =1}^K \exp\sqb{\alpha \frac{S_{l,i}}{Q_{l,i} + c}},
\end{equation}
where $\alpha, \, c > 0$ are tuning parameters that serve to govern the strength of the extent to
which the agent focuses on the greedy choice and to protect against division by zero, respectively.
The selection choices \eqref{eq:tempered_greedy} satisfy \eqref{eq:psi} and thus Definition \ref{def:SRME}
by construction.
\end{exam}

\begin{exam}
\label{ex:erev_roth}
Similar learning dynamics arise in human psychology and behavioral economics where an agent chooses future actions with a bias towards those that have accrued higher (un-normalized) cumulative reward \citep{erev1998predicting, xu2020reinforcement}. A popular example, known as Luce's rule \citep{luce1959individual}, uses sampling probabilities 
\begin{equation}
\pi_{i,k} = (S_{k,i} \vee \alpha) \bigg/ \sum_{l=1}^K  (S_{l,i} \vee \alpha), 
\label{eq:luce_rule}
\end{equation}
where $\alpha \geq 0$ is a tuning parameter governing the amount of baseline exploration. More generally, the agent may ascribe to arm $k$ a weight $f(S_{k,i})$, where $f$ is a non-negative potential function, 
and sample actions with probabilities proportional to the weights. The decision rule in \eqref{eq:luce_rule} only depends on $S$ and thus satisfies \eqref{eq:psi}. 
\end{exam}

\begin{exam}
\label{ex:thompson}
Thompson sampling is a popular Bayesian heuristic for running sequential experiments \citep{thompson1933likelihood}.
In Thompson sampling an agent starts with a prior belief distribution $G_0$ on the reward distributions \smash{$\cb{P_k}_{k = 1}^K$}.
Then, at each step $i$, the agent samples a reward distribution from their posterior belief $G_{i-1}$, pulls the optimal arm that has the highest mean according to this sampled distribution, and uses the resulting realized reward to update the posterior $G_i$ using Bayes' rule.  The motivation behind Thompson sampling is that it quickly converges to pulling the best arm, and thus achieves low regret
\citep{agrawal2017near,chapelle2011empirical}.
Thompson sampling does not always satisfy Definition \ref{def:SRME}. However, widely used modeling choices involving exponential
families for the \smash{$\cb{P_k}_{k = 1}^K$} and conjugate priors for $G_0$ result in these posterior probabilities $\rho_{k,i}$
satisfying the sufficiency condition \eqref{eq:psi} \citep{russo2018tutorial}, in which case Thompson sampling yields a sequentially
randomized Markov experiment in the sense of Definition \ref{def:SRME}.

In this paper, we will focus on one of the most widely used variants of this design, Gaussian Thompson sampling. In this
setting the agent models $P_k^n$ as a Gaussian distribution with (unknown) mean $\mu_k^n$ and (known) variance $\sigma_k^2$,
and sets the prior distribution $G_0$ on the mean rewards $\mu_k^n$ of each arm to be Gaussian with mean 0 and variance $(\nu^n)^2 I_{K \times K}$. (Note that the Gaussian assumption only affects the behavior of the agent's algorithm; the actual reward distributions are not necessarily Gaussian.) 
The posterior sampling step described above then takes on a simple form in this model: For each step $i = 1, \, \ldots, \, n$ we generate draws
\begin{equation}
	\label{eq:TS1_finite}
	\begin{split}
		&\tilde{\mu}_k^i \cond G_{k, i}^n \sim \nn\p{\frac{\sigma_k^{-2}S_{k,i}}{\sigma_k^{-2} Q_{k,i}+ \p{\nu^n}^{-2}}, \ \ \frac{1}{\sigma_k^{-2} Q_{k,i}+ \p{\nu^n}^{-2}}}, 
	\end{split}
\end{equation}
and then pull the arm corresponding to the largest \smash{$\tilde{\mu}_k^i$}. This sampling step clearly corresponds to a sequentially randomized
Markov experiment. Formally, here,
\begin{equation}
\label{eq:psup}
\psi\p{ Q_{i}, \, S_{i}} = \Pi^{\sup}\p{Q_{i}, \, S_{i},; \, \sigma, \, \nu^n},
\end{equation}
where \smash{$\Pi^{\sup}_k(q, \, s; \, \sigma, \, \nu^n)$} denotes the probability that the $k$-th draw \smash{$\tilde{\mu}_k^i$}
from \eqref{eq:TS1_finite} is largest, with $S_{l,i} = s_l$ and $Q_{l,i} = s_l$ for all $l = 1, \, \ldots, \, K$. 
\end{exam}

\begin{exam}
\label{ex:exploration}
Exploration sampling is a variant of Thompson sampling where, using notation from the above example, 
the agent pulls each arm with probability $\pi_{k,i} = \rho_{k,i} (1 - \rho_{k,i}) / \sum_{l = 1}^K \rho_{l,i} (1 - \rho_{l,i})$
instead of $\pi_{k,i} = \rho_{k,i}$ \citep{kasy2019adaptive}. Exploration sampling is preferred to Thompson
sampling when the analyst is more interested in identifying the best arm than simply achieving low regret
\citep{kasy2019adaptive,russo2020simple}. Exploration sampling satisfies Definition \ref{def:SRME}
under the same conditions as Thompson sampling.
\end{exam}

\begin{exam}
\label{ex:EXP3}
The Exp3 algorithm, proposed by \citet{auer2002nonstochastic}, uses sampling probabilities
\begin{equation}
\label{eq:exp3}
\pi_{i,k} = \exp\sqb{\alpha \, \sum_{j = 1}^{i-1} \frac{\mathbf{1}\p{\cb{A_j = k}} Y_j}{\pi_{j,k}}} \, \bigg/ \, \sum_{l =1}^K \exp\sqb{\alpha \, \sum_{j = 1}^{i-1} \frac{\mathbf{1}\p{\cb{A_j = l}} Y_j}{\pi_{j,l}}},
\end{equation}
where again $\alpha > 0$ is a tuning parameter. The advantage of Exp3 is that it can be shown to achieve
low regret even when the underlying distributions \smash{$\cb{P_k}_{k = 1}^K$} may be non-stationary and
change arbitrarily across samples. The sampling probabilities \eqref{eq:exp3} do not satisfy \eqref{eq:psi},
and so the Exp3 algorithm is not covered by the results given in this paper; however, it is plausible that
a natural extension of our approach to non-stationary problems could be made accommodate it. We leave
a discussion of non-stationary problems to future work.
\end{exam}

\subsection{A \rev{Weak Signal} Scaling}
\label{sec:scaling}

Next, we describe  a sequence of experiments, indexed by $n$, and how the dynamics associated with these experiments converge to a diffusion limit under appropriate scaling. 
In order for behaviors of sequentially randomized Markov experiments to admit a limit distribution, we will require both the reward distributions $ P^n $ and sampling  functions $\psi^n$ used in the $n$-th
experiment to converge in an appropriate sense. 
First, we will assume that reward distributions satisfy the following  \rev{weak signal} scaling.
Unless otherwise stated, all reward distributions are assumed to be in  this scaling regime for the remainder of the paper. 

\begin{defi}[\rev{Weak Signal} Regime of Reward Distributions]
	\label{def:reward_diff_scale}
	Consider a sequence of $K$-armed stochastic bandit problems in the sense of Definition \ref{def:bandit},
	with reward distributions $ \{P^n_k\}_{k,n \in \N} $. We say that this sequence resides in the \rev{weak signal} regime if there exist
	$ \mu, \sigma \in \R_+^K $ such that 
	\begin{equation}
		\label{eq:Pn_def}
		\limn \sqrt{n} \mu_k^n = \mu_k,  \ \ \ \ \ \ \ \ \
		\limn \p{\sigma_k^n}^2 = \sigma_k^2, 
	\end{equation}
where $\mu_k^n = \EE[P^n_k]{Y}$ and $\p{\sigma_k^n}^2 = \Var[P^n_k]{Y}$.
\end{defi}

A simple sequence of models in the \rev{weak signal} regime is the Gaussian model where rewards from arm-$k$ in the $n$-th system
are normally distributed, $Y_i \cond A_i = k \sim \calN(\mu^n_k, \sigma^2_k) $, where $ \mu^n_k \sqrt{n} \to \mu_k$ as $n$ grows
whereas the variance $ \sigma^2_k>0 $ remains fixed.
We can also construct a Bernoulli model in the \rev{weak signal} regime by considering centered rewards. If $Y_i = Z_i - \zeta_0$ where
$Z_i \cond A_i = k \sim \text{Bernoulli}(\zeta_k^n)$ and $\sqrt{n}(\zeta_k^n - \zeta_0) \to \delta_k$, then we are in the diffusion
regime with $\mu_k = \delta_k$ and $\sigma_k^2 = \zeta_0 (1 - \zeta_0)$. We discuss the further role of centering in our setting
in Remark \ref{rema:translation} below.

Next, we require the sequence of sampling functions to converge in an appropriate sense. As discussed further below, the natural scaling
of the the $Q_{k,i}$ and $S_{k,i}$ state variables defined in \eqref{eq:QS_raw} is
\begin{equation}
\label{eq:QS}
Q_{k,i}^n = \frac{1}{n} \sum_{j = 1}^i \mathbf{1}\p{\cb{A_j = k}}, \quad  S_{k,i}^n = \frac{1}{\sqrt{n}} \sum_{j = 1}^i \mathbf{1}\p{\cb{A_j = k}} Y_j.
\end{equation}
We then say that a sequence of sampling functions $\psi^n$ is convergent if it respects this scaling.

\begin{defi}[Convergent Sampling Function]
\label{defi:psi_conv}
Writing sampling functions in a scale-adapted way as follows,
\begin{equation} 
\bar\psi^n(q,s)= \psi^n\p{nq, \sqrt{n}s}, \quad q  \in [0, \, 1]^K, s \in \R^K,
\end{equation}
we say that a sequence of sampling functions $\psi^n$ satisfying \eqref{eq:psi} is convergent if
\begin{equation}
\bar\psi^n \stackrel{n \to \infty}{\longrightarrow} \psi
\end{equation}
uniformly over all compact subsets of $   [0, \, 1]^K \times \R^K $, for a limiting sampling function $\psi$.
\end{defi}

When using our results in applications, one key practical consideration is in understanding conditions under which
it is natural to consider a sequence of sampling functions $\psi^n$ that is
convergent in the sense of Definition \ref{defi:psi_conv}. We here revisit some of our earlier examples in the
context of these scaling considerations.

\setcounter{contex}{\getrefnumber{ex:tempered_greedy}}
\addtocounter{contex}{-1}

\begin{contex}[continued]
The tempered greedy method can immediately be seen to be convergent, provided we use a sequence
of tuning parameters $\alpha_n$ and $c_n$ satisfying $\limn \sqrt{n} \alpha_n = \alpha$ and $\limn n c_n = c$
for some $\alpha, \, c \in \RR_+$, resulting in a limiting sampling function
\begin{equation}
\label{eq:tempered_greedy_limit}
\psi_k\p{q,s} = \exp\sqb{\alpha \frac{s_k}{q_k + c}} \, \bigg/ \, \sum_{l =1}^K \exp\sqb{\alpha \frac{s_l}{q_l + c}}.
\end{equation}
For tempered greedy sampling to be interesting, we in general want
the limit $\alpha$ to be strictly positive, else the claimed diffusion limit \eqref{eq:diffLim} will be trivial.
Conversely, for the second parameter, both the limits $c > 0$ and $c=0$ may be interesting, but working
in the $c=0$ limit may lead to additional technical challenges due to us getting very close to taking an ill-conditioned ratio.
\end{contex}

\setcounter{contex}{\getrefnumber{ex:erev_roth}}
\addtocounter{contex}{-1}

\begin{contex}[continued]
Luce's rule is convergent whenever $\sqrt{n}\alpha_n \to \alpha > 0$.
\end{contex}

\setcounter{contex}{\getrefnumber{ex:thompson}}
\addtocounter{contex}{-1}

\begin{contex}[continued]
There are two ways to make Gaussian Thompson sampling \eqref{eq:TS1_finite} convergent.
The first is to scale the prior variance $(\nu^n)^2$ in a manner that matches the order-of-magnitude of 
the mean rewards in the \rev{weak signal} regime. Since the mean reward scales on the order of $1/\sqrt{n}$, this would amount to having 
\begin{equation}
\limn (\nu^n)^{-2}/n = c > 0.
\label{eq:nu_smooth}
\end{equation}
With this scaling and using notation from \eqref{eq:psup}, we can check the probability of the $k$-th posterior sample
being largest satisfies\footnote{The key observation here is that, in \eqref{eq:TS1_finite}, if we divide $Q_{i,k}$ and $\nu^{-2}$
by $n$ and $S_{i,k}$ by $\sqrt{n}$, then the sampling distribution of the \smash{$\tilde{\mu^i_k}$} is the same up to scaling up by
a factor $\sqrt{n}$. However,  scaling all the \smash{$\tilde{\mu^i_k}$} by $\sqrt{n}$ doesn't change which one is largest, and so
the middle equality in \eqref{eq:TS_scale} holds.}
\begin{equation}
\label{eq:TS_scale}
\begin{split}
\psi^n(nq, \, \sqrt{n}s)
&= \Pi^{\sup}(nq, \, \sqrt{n}s; \, \sigma^2, \, \nu^{-2}) \\
&= \Pi^{\sup}(q, \, s; \, \sigma^2, \, \nu^{-2}/n)
\to \Pi^{\sup}(q, \, s; \, \sigma^2, \, c),
\end{split}
\end{equation}
and thus the induced sampling functions are in fact convergent. The second is to scale the prior variance such that
\begin{equation}
\limn (\nu^n)^{-2}/n = 0, \ \ \ \ \ \psi^n(nq, \, \sqrt{n}s) \to \Pi^{\sup}(q, \, s; \, \sigma^2, \, 0),
\label{eq:nu_undersmooth}
\end{equation}
in which case the prior is asymptotically uninformative about diffusion scale mean rewards. We will again see that the 
\eqref{eq:nu_smooth} and \eqref{eq:nu_undersmooth} lead to markedly different statistical properties; we refer to
the former as \textit{smoothed} Thompson sampling and to the latter as \textit{undersmoothed} Thompson sampling.
\end{contex}

\setcounter{contex}{\getrefnumber{ex:exploration}}
\addtocounter{contex}{-1}

\begin{contex}[continued]
As before, exploration sampling is convergent whenever the associated Thompson sampling rule is.
\end{contex}

\begin{rema}[Non-Zero Mean Rewards]
	\label{rema:translation}
	The scaling condition in Definition \ref{def:reward_diff_scale} implies that, in large samples, all arms have roughly
	zero rewards on average, i.e., $\limn \mu_k^n = 0$. In some applications, however, it may be more natural to
	consider a local expansion around non-zero mean rewards, where
	\begin{equation}
		\label{eq:local_expansion}
		\limn \sqrt{n} \p{\mu_k^n - \mu_0} = \delta_k
	\end{equation}
	for some potentially non-zero $\mu_0 \neq 0$ and $\delta \in \RR^K$.
	In our general results, we focus on the setting from Definition \ref{def:reward_diff_scale};
	however, we note that, when applied to any translation-invariant algorithm (i.e., a sequential experiment whose sampling
	function is invariant to adding a fixed offset to all rewards), any results proven under the setting of
	Definition \ref{def:reward_diff_scale} will also apply under \eqref{eq:local_expansion}. In Section \ref{sec:undersmoothedTS}, we will use this fact when studying a translation-invariant Thompson sampling algorithm.
\end{rema}

\section{Convergence to a Diffusion Limit under \rev{Weak Signal Scaling}}
\label{sec:limit}

We are now ready to state our first main result: Given a sequence of reward distributions satisfying \eqref{eq:Pn_def} and
under a number of regularity conditions discussed further below, the sample paths of the scaled
statistics $Q_{k,i}^n$ and $S_{k,i}^n$ of a sequentially randomized Markov experiments with convergent
sampling functions converge in distribution to the solution to a stochastic differential equation
\begin{equation}
\label{eq:diffLim}
\begin{split}
&dQ_{k,t} = \psi_k\!\p{Q_t, S_t} dt, \\
&dS_{k,t} = \mu_k \, \psi_k\!\p{Q_t, S_t} dt + \sigma_k \sqrt{\psi_k\!\p{Q_t, S_t}} dB_{k,t},
\end{split}
\end{equation}
where $B_{\cdot, \, t}$ is a standard $K$-dimensional Brownian motion, $\mu_k$ and $\sigma_k$
and the mean and variance parameters given in \eqref{eq:Pn_def}, and the time variable $t \in [0, \, 1]$
approximates the ratio $i/n$. A formal statement is given in Theorem \ref{theo:SDE}.

We will consider separately two cases depending on whether the limiting sampling function, $\psi$,
is Lipschitz-continuous, or continuous but non-Lipschitz.
We begin with the case where $ \psi $ is Lipschitz, which allows us to obtain stronger convergence guarantees,
before extending the results to the non-Lipschitz case. As will be discussed further below, tempered greedy and Thompson
sampling with $c > 0$ as well as Luce's rule are all examples with Lipschitz limiting sampling functions, while tempered
greedy and Thompson sampling with $c = 0$ are examples of the non-Lipschitz case.

\begin{rema}
In our current analysis, we do not consider cases where the limiting sampling function $\psi$ is not continuous
with respect to the underlying state $(Q,S)$. Considering the discontinuous case could also be of interest,
as there are a number of bandit algorithms---including upper-confidence bound (UCB) and  $\epsilon$-greedy algorithms---that are
consistent with Definition \ref{def:SRME} but with sampling functions $\psi$ that have sharp cutoffs.

\end{rema}

\begin{rema}
We generally assume that the reward variances are known and do not consider the problem of estimating them. This is because, \rev{under the weak signal scaling}, the arm variances stay constant while the number of samples scales, and therefore the variances can be accurately estimated. 
\end{rema}

\subsection{Lipschitz Limiting Sampling Functions}
Given any Lipschitz sampling functions, we will show that a suitably
scaled version of the process $(Q^n_i, S^n_i)$ converges to an It\^o diffusion process. 
Define $\bQ^n_t$ to be the linear interpolation of $Q^n_{\lfloor t n \rfloor}$,
\begin{align}
\bQ^n_{k,t} = (1-tn + \lfloor tn \rfloor)Q^n_{k, {\lfloor tn \rfloor}} +  (tn  - \lfloor tn \rfloor )Q^n_{k, \lfloor tn \rfloor +1},  \quad t \in [0,1], k = 1, \ldots, K,
\label{eq:linear_interp}
\end{align}
and define the process $\bS^n_t$ analogously.\footnote{Note that the above process is not adapted to its natural filtration as a result of the linear interpolation; this feature, however, will not impact our results. }
Here, $n$ is the total number of samples in the experiment, and the continuous  index $t$ can be thought of as the fraction of total samples collected up to a certain point. For instance,  $(\bQ^n_{0.5}, \bS^n_{0.5})$ can be viewed as the state of the system at the mid-point of an experiment. 

Let $\calC$ be the space of continuous functions $[0,1] \mapsto \R^{2K}$ equipped with the uniform metric: 
$d(x,y) = \sup_{t \in [0,1]} |x(t)-y(t)|$, $x,y \in \calC$.  
We have the following result; the proof is given in Section \ref{sec:proof:theo:SDE}. 

\begin{theo}
\label{theo:SDE}
Fix $K\in \N$, $\mu \in \R^K$ and $\sigma \in \rp^K$. 
Suppose that we have a sequence of $K$-armed bandit problems as in Definition \ref{def:bandit}, whose
reward distributions $ \{P^n_k\}_{k,n \in \N} $ reside in the \rev{weak signal} regime as per Definition \ref{def:reward_diff_scale}
and admit  fourth moments that are uniformly bounded across $ k $ and $ n $.
 Suppose that we have a convergent sequence of sequentially randomized Markov experiments following
Definitions \ref{def:SRME} and \ref{defi:psi_conv}.
Suppose furthermore that the limit sampling function $ \psi $ is Lipschitz-continuous, and $(\bQ^n_0, \bS^n_0) = 0$. Then, as $n \rightarrow \infty$,
$(\bQ^n_t, \bS^n_t)_{t\in [0,1]}$ converges weakly to $(Q_t, S_t)_{t\in [0,1]} \in \calC$,
which is the unique solution to \rev{the stochastic differential equation \eqref{eq:diffLim}} over $t \in [0,1]$,
where $(Q_0, S_0) = 0$  and $B_t$ is a standard Brownian motion in $\R^K$ with independent components.\footnote{Here and henceforth, all standard Brownian motions in this paper are assumed to have independent components. } 
Furthermore,  
\begin{equation}
\label{eq:bdd_func}
\lim_{n \rightarrow \infty} \EE{f {  (\bQ^n_t, \bS^n_t) }} =  \EE{f(Q_t, S_t)}, \quad \forall t \in [0,1]
\end{equation}
for any bounded continuous function $f: \R^{2K} \mapsto \R$.
\end{theo}

As an immediate corollary of Theorem \ref{theo:SDE}, we obtain the following characterization of the finite-sample expected  regret; the proof follows simply by setting $ f(Q,S) :=  (\max_k \mu_k)-  \left \langle Q, \mu \right\rangle $. 

\begin{coro}[Convergence of Expected Regret]
	\label{cor:convg_Rn}
	Fix $K\in \N$, $\mu \in \R^K$ and $\sigma \in \rp^K$.   Suppose that the limit sampling function $ \psi $ is Lipschitz-continuous.  Define the diffusion regret
	\begin{equation}\label{eq:scaled_regret}
		R = (\max_k \mu_k)-  \left \langle Q_1, \mu \right\rangle, 
	\end{equation}
	where $ \{Q_t\}_{t\in [0,1]} $ is given by the solution to \eqref{eq:diffLim}.  Then,
	\begin{equation}\label{key}
     {\EE{R^n}} =   \EE{R } \sqrt{n}+ o(\sqrt{n}),
	\end{equation}
i.e., $ \lim_{n\to \infty} \EE{R^n}/\sqrt{n}= \EE{R} $. 
\end{coro}

Finally, the following theorem gives a more compact representation of the stochastic differential equations in Theorem \ref{theo:SDE}, showing that they can be written as a set of ordinary differential equations driven by a Brownian motion with a random time change, $t \Rightarrow Q_t$.  The result will be useful, for instance, in our subsequent analysis of Thompson sampling in the super diffusive regime. The proof is given in Section \ref{sec:proof:theo:timechange}. 

\begin{theo}
\label{theo:timechange}
The limit stochastic differential equation in \eqref{eq:diffLim} can be equivalently written as 
\begin{align}
dQ_{k,t} =& \psi_k(  Q_t, Q_{t}\mu + \sigma W_{Q_{t}}) \, dt, \quad k = 1, \ldots, K, 
\label{eq:Qt_timechange}
\end{align}
with $Q_0=0$, where $W$ is a $K$-dimensional standard Brownian motion. Here, $Q_t\mu$ and $\sigma W_{Q_t}$ are understood to be vectors of element-wise products, with $Q_t\mu = (Q_{k,t}\mu_k)_{k=1, \ldots, K}$, and $\sigma W_{Q_t} = (\sigma_k W_{k, Q_{k,t}})_{k =1, \ldots, K}$.  In particular, we may also represent $S_t$ explicitly as a function of $Q$ and $W$: 
\begin{equation}
S_{k,t} =  Q_{k,t} \mu_k + \sigma_k W_{k,Q_{k,t}}, \quad k = 1, \ldots, K, \, t \in [0,1]. 
\label{eq:St_timechange}
\end{equation}
\end{theo}

All proofs are deferred to Section \ref{sec:proofs}.
Our proof of Theorem \ref{theo:SDE} uses the Stroock-Varadhan program which is in turn based on the martingale characterization of diffusion \citep{durrett1996stochastic, stroock2007multidimensional}. The main technique is based on showing that the the generator of the Markov chain associated with the sequentially randomized Markov experiment converges, in an appropriate sense, to the generator of the desired limit diffusion process. 
Meanwhile, Theorem \ref{theo:timechange} builds upon the convergence result in Theorem \ref{theo:SDE}. The key additional step is
to use the Skorohod's representation theorem  so as to allow us to couple all relevant sample paths, including, $ S $, $ Q $
and a noise process $ U_t $ under a random time-change $ U_t \to U_{Q_t} $, to a single Brownian motion, and
show that they converge to the appropriate limits jointly. 


\subsection{Non-Lipschitz Limiting Sampling Functions}
\label{sec:nonLip}
In some cases,  such as the undersmoothed Thompson sampling studied later in this paper, the sequence of sampling functions may converge to a limit that is continuous but non-Lipschitz. 
In this subsection, we show how to characterize the limiting sample paths under a
 non-Lipschitz  sampling function. The main takeaway here is that now the resulting diffusion limits may not be unique, but they can still be characterized by the same set of SDEs. 

Denote by $\calD = [0, 1]^K \times \R^K \to \mathbf{\Delta}^K$ the domain of the sampling functions. Fix a continuous limiting sampling function, $ \psi $. By a version of the Stone–Weierstrass theorem, we can find a sequence of function $\{ \psi^j\}$, such that $ \psi^j$ is Lipschitz-continuous for all $ j $, and $ \psi^j$ converges to $ \psi $ uniformly over compact subsets of $ \calD $, as $ j\to \infty $. We have the following lemma; the proof is given in Appendix \ref{app:lemm:tightness}. 
\begin{lemm}
	\label{lemm:tightness}
	Fix a continuous function $ \psi : \calD \to \mathbf{\Delta}^K $ and a sequence of Lipschitz-continuous functions  $\{\tilde \psi^j\}_{j\in \N} $ that converge  to $ \psi $  uniformly over  compact sets. Denote by $(Q^j, S^j)$ the unique solution to the ODEs associated with a limiting sampling function $\tilde \psi^j$ in Theorem \ref{theo:timechange}. The following holds almost surely: 
	\begin{enumerate}
	\item  $ \{(Q^j, S^j)\}_{j\in \N} $ is tight, in the sense that any of its subsequences  admits a further subsequence that converges uniformly over $ [0,1] $ to a limit that is continuously differentiable over $[0,1]$. We say that $ (Q,S)  $ is a limit function if it is a limit point for one of these convergent subsequences.
	
	\item 	Every such limit function $ (Q,S) $ is a solution to the ODEs in Theorem \ref{theo:timechange} with sampling function $ \psi $. 
\end{enumerate}
\end{lemm}

We now use Lemma \ref{lemm:tightness}  to show for every continuous limiting sampling function $ \psi $ we can construct  an appropriately convergent sequence of (pre-limit) sampling functions; the proof follows directly from Lemma \ref{lemm:tightness}   via a triangular array argument \rev{and is given in Appendix \ref{app:theo:conv_triangular_psi}. }
	\begin{theo}
		\label{theo:conv_triangular_psi}		
Fix $K\in \N$, $\mu \in \R^K$ and $\sigma \in \rp^K$. Suppose that we have a sequence of $K$-armed bandit problems as in Definition \ref{def:bandit} whose
reward distributions $ \{P^n_k\}_{k,n \in \N} $ reside in the \rev{weak signal} regime as per Definition \ref{def:reward_diff_scale}
and admit  fourth moments that are uniformly bounded across $ k $ and $ n $.  Fix a continuous limiting sampling function $ \psi $. Let $\{\psi^{n,j}\}_{n,j \in \N}$ be an array of sampling functions such that 
\begin{enumerate}
	\item For all $ j $, the scale-adjusted sampling function $ \bar \psi^{n,j} $  (Definition \ref{defi:psi_conv}) converges to a Lipschitz-continuous function $ \psi^j $ uniformly over compact sets as $ n\to \infty $ 
	\item $\psi^j$ converges to $ \psi $ uniformly over compact sets as $ j \to \infty $. 
\end{enumerate}
Let  $(\bar Q^{n,j}_t,\bar S^{n,j}_t)_{t \in [0,1]} $  be the system dynamics in the $ n $-th experiment under the sampling function $ \psi^{n,j} $, where $(\bQ^n_0, \bS^n_0) = 0$.  Then, there exists a sequence  $ j_n  \in \N$,  with $ j_n \to \infty $ as $ n\to \infty $, such that almost surely   $ \{ (\bar Q^{n, j_n},\bar S^{n, j_n}) \}_{n \in \N}$  converges uniformly to a solution of the ODEs in Theorem \ref{theo:timechange} with sampling function $ \psi $. 
	\end{theo}

Practically, Theorem \ref{theo:conv_triangular_psi} shows that for a continuous limiting sampling function $ \psi $, \rev{there exists} a sequence of pre-limit sampling functions $ \psi^n $, such that the resulting system dynamics converge \rev{under the weak signal scaling} to a solution of the ODEs in Theorem \ref{theo:timechange}. While these solutions may not be unique when $ \psi $ is non-Lipschitz, the theorem demonstrates that characterizing their behavior would nevertheless offer meaningful insights for the behavior of the finite-sample, pre-limit system. Indeed, we will use these characterizations in subsequent sections for analyzing  the regret of Thompson sampling 
algorithms  associated with a  non-Lipschitz limiting sampling function. \rev{Finally, we should note that the approximation result in Theorem \ref{theo:conv_triangular_psi} does not automatically apply to every arbitrary sequence of pre-limit sampling functions that converges to a continuous, though non-Lipschitz, limit, as such a statement would likely require  further restrictions on how quickly the smoothness of the sampling functions is allowed to decrease as the sample size $ n $ grows large. Establishing the conditions  on the joint scaling of the Lipschitz-continuity of the sampling function and the sample size that would ensure such convergence is an interesting question for future research.}

\section{Applications to Regret Analysis}

Having established a general diffusion limit for sequentially experiments \rev{in the weak signal asymptotic regime}, \rev{in this section, we apply
this limit} to gain insights about the regret properties of a number of practical sequential designs.

We first turn our attention to those designs with Lipschitz sampling functions, i.e., ones satisfying the
setting of Theorem \ref{theo:SDE} directly. \rev{These designs  have the desirable stability property where the sampling probabilities  do not vary drastically as a function of the state variables. Furthermore, they arise naturally in some of the algorithms  considered from Examples~\ref{ex:tempered_greedy}--\ref{ex:thompson}, including Thompson sampling. Interestingly, we will demonstrate that while
these designs are perhaps the most amenable to formal study, 
they do not have  robust regret guarantees: In Section \ref{sec:super_diffusive}, we show that any sequential
experiment with a Lipschitz sampling function will perform poorly for problems with \rev{``large''
arm gaps} (the meaning of ``large'' will be made precise below).}

In Section \ref{sec:undersmoothedTS}, however, we provide a counterpoint \rev{to this negative result}. We study \rev{
a variant of undersmoothed Thompson sampling (i.e., with an asymptotically uninformative
prior on treatment effects) for the two-armed bandit problem. Importantly, this algorithm does not have a Lipscthiz-continuous sampling function. Using
a mix of numerical and analytic arguments starting from our diffusion limit, we show that this undersmoothed Thompson sampling algorithm has excellent regret performance {under the weak signal scaling}---both for ``moderate'' and ``large''
effects sizes.}

In addition to being of independent interest, we hope that these results help highlight the versatility
of the diffusion limit: It is helpful both as a tool for unified analysis of large classes of sequential designs
(as in the case of Lipschitz experiments), and as a tool for sharp characterization of specific
designs of interest (as in the case of undersmoothed Thompson sampling).

\subsection{Lower Bounds for Asymptotically Lipschitz Experiments}
\label{sec:super_diffusive}

Sequentially randomized experiments whose sampling functions converge to a Lipschitz limit seem like
an a-priori desirable class of experiments to consider: Their sampling functions are not too sensitive
to observed data, and in many settings this type of stability can lead to good behavior. In the
\rev{weak signal asymptotic regime}, however, we find that this is generally not the case---at least if our goal is to achieve a small 
regret, \rev{as defined in} \eqref{eq:scaled_regret}. We find that, by considering arm reward gaps  that are large \rev{in the diffusion limit},
\rev{ \textit{any} Lipschitz experiment admits a regret that exceeds the minimax regret by an arbitrarily large factor}.

\rev{We start by specifying a class of problems where the arm reward gaps are large in the diffusion limit, as follows.} Fix $ \tilde \mu \in  \R^{K-1} $, where, without the loss of generality, we assume that the entries of $ \tilde \mu$ are listed in a decreasing order. We start with a mean reward vector parameterized by $ \delta >0 $, where 
\begin{equation}
	\mu^\delta_1 = \tilde\mu_1 + \delta,  \mbox{ and }   \p{\mu^\delta_2, \ldots, \mu^\delta_K} = \tilde \mu, 
\label{eq:mu_delta}
\end{equation}
that is, the reward gap between the top two arms in $ \mu^\delta $ is equal to $ \delta $. Finally, we will
consider problems where $\delta$ grows, i.e., the gap between the top two arms gets large, while we
keep $\tilde \mu$ and the variances $\sigma^2$ fixed.  
\rev{
We refer to the above problem setting, i.e., where
the experiment with $n$ samples has $\mu^n = \mu^\delta/\sqrt{n}$ for large values of $\delta$, as the \textit{super diffusive} regime. This regime can be thought of taking the pre-limit system through two scalings in the following order: First, for a fixed $\delta$, and therefore fixed $\mu^\delta$, we take the weak signal scaling as per Definition \ref{def:reward_diff_scale}; then, we look at how the system behaves in the diffusion limit as the arm gap $\delta$ grows large.  }

\rev{Next, we introduce a regret performance baseline in this super diffusive regime so that we can use it to subsequently evaluate the regret of asymptotically Lipscthiz experiments.  We recall a well known}
fact that sequential experiments are generally able to achieve lower regret in the presence of
large arm gaps (which makes identifying the best arm easier). The following \rev{result} characterizes the minimax
\rev{regret} in the super diffusive regime. The key observation here is that problems with large gap $\delta$
on the diffusion scale allow for smaller regret: The scaled regret \rev{vanishes at the rate $\log(\delta)/\delta$ as $\delta \to \infty$}.
The lower bound \eqref{eq:superdiff_minimax_lb} follows directly from \citet{mannor2004sample} (Eq.~(17)),
while the upper bound follows from using a modified UCB algorithm proposed in \citet{auer2010ucb} (Theorem 3.1).
We note that neither \citet{mannor2004sample} nor \citet{auer2010ucb} considered the \rev{weak signal asymptotics}
 in their analysis; however, we see that their bounds take on a particularly simple
form when considered in this limit.

\begin{prop}
\label{prop:superdiff_minimax}
Fix $\delta>\sqrt{e}$ and $ K \geq 2 $. There exist constants $ c_1, c_2>0 $, such that for any $n$-indexed family of sequentially randomized experiments, there exists a sequence of $ n $-indexed, $ K $-armed bandit instances with arm gaps $ \{ \delta / \sqrt{n}\}_{n \in \N}$, such that: 
\begin{equation}
\label{eq:superdiff_minimax_lb}
\liminf_{n \rightarrow \infty} \frac{\EE{R^n}}{\sqrt{n}} \geq  \frac{c_1K\log(c_2\delta^2/K)}{\delta}. 
\end{equation}
Moreover, there exists a sequentially randomized experiment and constants $C >0 $, such that for any sequence of $ n $-indexed, $ K $-armed bandit instances with arm gaps $ \{ \delta / \sqrt{n}\}_{n \in \N}$,
\begin{equation}
\label{eq:superdiff_minimax_ub}
\limsup_{n \rightarrow \infty} \frac{\EE{R^n}}{\sqrt{n}} \leq \frac{C K\log (\delta) }{\delta}. 
\end{equation}
\end{prop}

The guarantees given in Proposition \ref{prop:superdiff_minimax} characterize the \rev{minimax regret} for 
arbitrary sequential experiments; not necessarily ones covered by our theory or admitting a diffusion limit.
However, in order to understand the behavior of convergent sequentially randomized experiments with
a Lipschitz limiting sampling function, we will use Proposition \ref{prop:superdiff_minimax} as our
benchmark: Can asymptotically Lipschitz experiments attain the behavior \eqref{eq:superdiff_minimax_ub}?

\rev{Our main result in this subsection provides a (strongly)}  negative answer to this question.
\rev{We show that} Lipschitz experiments are not just unable to achieve the optimal $\log(\delta)/\delta$ \rev{regret} scaling
in the arm gap; they \rev{in fact cannot even ensure a vanishing regret} as $\delta$ gets large.
We say that  a limiting sampling function $ \psi $ is non-trivial if $ \max_{k=1, \ldots, K} \psi_k(0,0) <1$,
i.e., the  sampling function does not assign overwhelming probability  to any single action in the absence of any
observations. We show below that under any non-trivial, Lipschitz-continuous limiting
sampling function, the regret is always bounded away from zero in the super diffusive regime; the proof is
given in Appendix \ref{app:lip_bad}.

\begin{theo} 
	\label{theo:lip_bad}
For any $ L>0 $, $ \tilde \mu \in \R^K $, and sampling function $ \psi $ that is $ L $-Lipschitz with $\max_{k=1, \ldots, K} \psi_k(0,0)  = \tilde p <1  $, we have that, almost surely, 
\begin{equation}\label{eq:lip_bad}
	\liminf_{\delta \to \infty } R \geq \frac{(1-\tilde p)^2}{8KL}. 
\end{equation}
where $ R $ is the diffusion regret defined in \eqref{eq:scaled_regret}. 
\end{theo}

We end this section by considering implications of Theorem \ref{theo:lip_bad} in the context of a number
of convergent sequential experiments discussed above. We note that the lower bound in Theorem \ref{theo:lip_bad}
also provides a quantitative description of the regret plateau, which depends inverse-proportionally on the
sampling function's Lipschitz-constant, $L$. Thus, the smaller Lipschitz constant a limiting sampling
function has, the more \eqref{eq:lip_bad} limits any good behavior we could hope for in the presence of large arm gaps.

\setcounter{contex}{\getrefnumber{ex:tempered_greedy}}
\addtocounter{contex}{-1}

\begin{contex}[continued]
Consider convergent tempered greedy sampling as discussed above with $\alpha, \, c > 0$, and
limiting sampling function
\begin{equation}
	\psi_k\p{q,s} = \exp\sqb{\alpha \frac{s_k}{q_k + c}} \, \bigg/ \, \sum_{l =1}^K \exp\sqb{\alpha \frac{s_l}{q_l + c}}. 
\end{equation}
This function is Lipschitz with constant $L = 2K\alpha/c$, and so Theorem \ref{theo:lip_bad}
implies that $ \liminf_{\delta \to \infty } R \geq \frac{c(1-\tilde p)^2}{16  K^2 \alpha }  $ under tempered greedy sampling. 
\end{contex}

\setcounter{contex}{\getrefnumber{ex:erev_roth}}
\addtocounter{contex}{-1}

\begin{contex}[continued]
The convergent version of Luce's rule has limiting sampling function
\begin{equation}
\psi_k(q,s) = (s_k \lor \alpha) \, \bigg/\, \sum_{l = 1}^K (s_l \lor \alpha)
\end{equation}
for some $\alpha > 0$.
This function is Lipschitz with constant $L =  1/\alpha + K$, and so Theorem \ref{theo:lip_bad}
implies that $ \liminf_{\delta \to \infty } R \geq \frac{(1-\tilde p)^2}{8(K^2 + K/\alpha)}  $ under Luce's rule. 
\end{contex}

\setcounter{contex}{\getrefnumber{ex:thompson}}
\addtocounter{contex}{-1}

\begin{contex}[continued]
Smoothed Gaussian Thompson sampling, i.e., with prior variance satisfying \eqref{eq:nu_smooth}
can also be verified to be Lipschitz.
Now, recall that using the smooth scaling \eqref{eq:nu_smooth} involves choosing priors that are
informative on the diffusion scale. Thus, the following result implies that, on the diffusion scale, Gaussian
Thompson sampling with an informative prior is not able to achieve vanishing regret with large arm gaps; the proof is given in Appendix \ref{app:coro:smoothGaussBad}. 

\begin{coro}
	\label{coro:smoothGaussBad}
Fix $K>0$. For any Gaussian Thompson sampling agent with independent priors such that \eqref{eq:nu_smooth} holds, its sampling function converges to a Lipschitz-continuous limiting sampling function in the \rev{weak signal} regime. 
Thus, by Theorem \ref{theo:lip_bad}, its diffusion limit regret in the  super diffusive regime is non-vanishing: 
	\begin{equation}\label{key}
		\liminf_{\delta \to \infty } R  >0. 
	\end{equation}
	\end{coro}
\end{contex}

\subsection{Undersmoothed Thompson Sampling}
\label{sec:undersmoothedTS}

Above, we showed that sequentially randomized experiments with asymptotically
Lipschitz sampling functions can perform extremely poorly in the presence of
large arm gaps. However, the regret lower bounds in Theorem \ref{theo:lip_bad}
decay as the Lipschitz constant of $\psi$ gets larger, leading to a natural
question: \rev{for some of the algorithms considered from Examples~\ref{ex:tempered_greedy}--\ref{ex:thompson}, can we achieve better regret performance in the {weak signal regime} by simply choosing tuning
parameters  so as to make their limiting sampling functions less smooth? }
Here, we provide a positive
answer in the case of Thompson sampling (Example~\ref{ex:thompson}): We use our results to show
that an undersmoothed (and thus non-Lipschitz) variant of two-armed Thompson sampling has 
excellent \rev{regret performance under both small and large arm gaps in the weak signal regime.}

We consider the following algorithm, which we refer to as translation-invariant Thompson sampling.
In periods $i = 1, \, \ldots, \, n$, an agent chooses which of two distributions $P^n_{1}$ or $P^n_{2}$
to draw from, each with (unknown) mean $\mu^n_k$ and (known) variance $(\sigma^n_k)^2$. To streamline notation, we will present the analysis for the case where $ \sigma^n_1=\sigma^n_2=\sigma^n$, with the understanding that all results stated in this section will generalize to the case with heterogeneous reward variances in a straightforward manner. 

In a finite-horizon, pre-limit system, the agent uses the following version of translation-invariant Thompson sampling based on reasoning about the posterior distribution of the arm difference $\delta^n = \mu^n_1 - \mu^n_{2}$. The agent starts with one draw from each arm, and subsequently pulls arm $1$ in period $i$ with probability: 
\begin{equation}
	\label{eq:TS2_finite}
	\pi^n_i = \Phi\p{\frac{\alpha_i^{-2} \Delta_i}{\sqrt{\alpha_i^{-2} + (\nu^n)^{-2}}}}, \ \ \  \mbox{with } \ \  \Delta_i = {\frac{S_{1,i}}{Q_{1,i}} - \frac{S_{2,i}}{Q_{2,i}}}, \ \ \ \ \alpha_i^2 = \frac{\sigma^2i}{Q_{1,i}Q_{2,i}},
\end{equation}
where $ \nu^n $ is interpreted as the prior standard deviation for $ \delta^n$, $ \Delta_i $ the empirical mean of $ \delta^n $, and $ \alpha_i^2 $ the variance associated with the noisy realizations of rewards.
Here, we note that $Q_{1,i} + Q_{2,i} = i$. 
Now, to obtain a diffusion limit, we again need to choose a scaling for the prior standard
deviation, $\nu^n$, that yields a convergent sequence of sampling functions. As before,
we will assume that, $\limn (\nu^n)^{-2} / n = c$, for some $ c \geq 0 $, resulting
in limiting sampling probabilities
$$ \limn \pi^n_{\lfloor tn \rfloor} = \Phi\p{\frac{\sigma^{-2} Q_{1,t} (t - Q_{1,t}) \p{S_{1,t} / Q_{1,t} - S_{2,t} / (t - Q_{1,t})}}{\sqrt{\sigma^{-2} t Q_{1,t} (t - Q_{1,t}) + t^2 c}}}. $$
Theorem \ref{theo:conv_triangular_psi} implies that there exists a sequence of sequentially
randomized experiments with the same limiting sampling functions whose sample paths converge to:
\begin{align}
\label{eq:TS2}
\begin{split}
	&dS_{1,t} = \delta \pi_t dt + \sqrt{\pi_t} \sigma dB_{1,t},  \ \ \ \ dS_{2,t} =  \sqrt{1 - \pi_t}  \sigma  dB_{2,t}, \\
	&dQ_{1,t} = \pi_t dt, \ \ \ \ \pi_t = \Phi\p{\frac{\sigma^{-2} Q_{1,t} (t - Q_{1,t}) \p{S_{1,t} / Q_{1,t} - S_{2,t} / (t - Q_{1,t})}}{\sqrt{\sigma^{-2} t Q_{1,t} (t - Q_{1,t}) + t^2 c}}} \\
	& Q_{2,t} = t- Q_{1,t},
\end{split}
\end{align}
with $S_{k,0} = Q_{k,0} = 0$. Below, we will focus on the behavior of translation-invariant two-armed
Thompson sampling in the undersmoothed ($c = 0$) \rev{setting}, but will also empirically evaluate
the smoothed ($c > 0$) setting for different effect sizes.

\begin{rema}
As usual, we focus on the behavior of \eqref{eq:TS2_finite} in the \rev{weak signal} regime as in Definition \ref{def:reward_diff_scale}.
However, because this algorithm is translation invariant, we can always without loss of generality assume that $\mu^n_2 = 0$
when studying its behavior. In that case---as discussed in Remark \ref{rema:translation}---our results apply
as long as $\limn \sqrt{n}(\mu^n_1 -  \mu^n_2) = \delta$ and $(\sigma^n_k)^2 = \sigma^2$, even if the mean arm rewards
$\mu_k^n$ themselves may not converge to 0. When stating results below, we assume that $\mu^n_2 = 0$, in which
case $\mu^n_1 = \delta^n$.
\end{rema}


\subsubsection{Super Diffusive Analysis}

Given that \rev{in Theorem \ref{theo:lip_bad}} we found that \rev{sequential experiments with Lipschitz sampling functions} to perform particularly
poorly in the super diffusive regime (i.e., with a large arm gap on the diffusion scale), we start by
looking at how undersmoothed translation-invariant Thompson sampling does in this regime. \rev{Our main result in this subsection shows that}, as we had hoped, exiting the Lipschitz regime enables us to get considerably better
behavior in terms of regret; the proof is given in Appendix \ref{app:theo:undersmooth_regret}. 

\begin{theo}
	\label{theo:undersmooth_regret}
	Consider the diffusion limit associated with two-armed undersmoothed Thompson sampling, with $ c=0 $,
	and assume that $\sigma = 1$ in \eqref{eq:TS2}.
	Then, the following holds almost surely for the \rev{limiting} regret \eqref{eq:scaled_regret} in the super diffusive regime:
	\begin{align}
		R \prec 1/\delta , \quad \mbox{as  $\delta \to \infty$}, 
	\end{align}
where $ \delta = \mu_1 -\mu_2 $ is the arm gap
and we write  $f(x) \prec g(x)$ to indicated that
for any $\beta \in (0,1)$, we have $f(x)/g(x)^\beta \to 0$.
\end{theo}

We note that the above theorem not only establishes that
the regret of undersmoothed Thompson sampling vanishes in the super diffusive regime, a stark contrast to the regret scaling of smoothed Thompson sampling  shown in Corollary \ref{coro:smoothGaussBad}; it also provides a quantitative characterization of the regret profile
when the arm gap is large. In particular,  we find that regret decays faster than
$1/\delta^{1-\epsilon}$ for any $\epsilon>0$ as $ \delta \to \infty$ in the super diffusive regime,
thus nearly matching the minimax rate given in Proposition \ref{prop:superdiff_minimax}.\footnote{
Theorems \ref{theo:lip_bad} and \ref{theo:undersmooth_regret}  are given in an almost-sure sense (with respect to the measure associated with the Brownian motion), and it is natural to ask whether analogous statements can be established for expected regret $ \EE{R} $ as well. Since regret is always non-negative, almost-sure regret lower bounds immediately extend to expected regret. For instance, it follows from  Theorem \ref{theo:lip_bad} that  the expected limit regret under a Lipschitz-continuous  sampling function satisfies
	\begin{align}
  \liminf_{\delta \to + \infty} \EE{R} >  0. 
	\end{align}
	Unfortunately, almost-sure regret upper bounds, on the other hand, do not  extend immediately. This is because $ R $ can be as large as  $ \delta $ in the worst case, which diverges as $ \delta \to \infty $, and as such we do not have an easy tightness property to rely on in order to extend the almost-sure guarantees to expected regret. Showing the same super diffusive upper bound holds for expected regret is an important question for further work. 
}

To the best of our knowledge, this is the first formal result suggesting that Thompson sampling achieves
anything close to instance-optimal behavior as the effect size grows large in the regime with moderate sample size. 
Algorithms currently known to attain regret upper bounds on the order of \eqref{eq:superdiff_minimax_ub}
tend to rely on substantially more
sophisticated mechanisms, such as adaptive arm elimination and time-dependent confidence intervals \citep{auer2010ucb}.
It is thus both surprising and encouraging that such a simple and easily implementable heuristic as Thompson sampling should
achieve near-optimal instance-dependent regret.

\subsubsection{Regret Profiles}

\begin{figure}
	\begin{tabular}{cc}
		\includegraphics[width=0.45\textwidth]{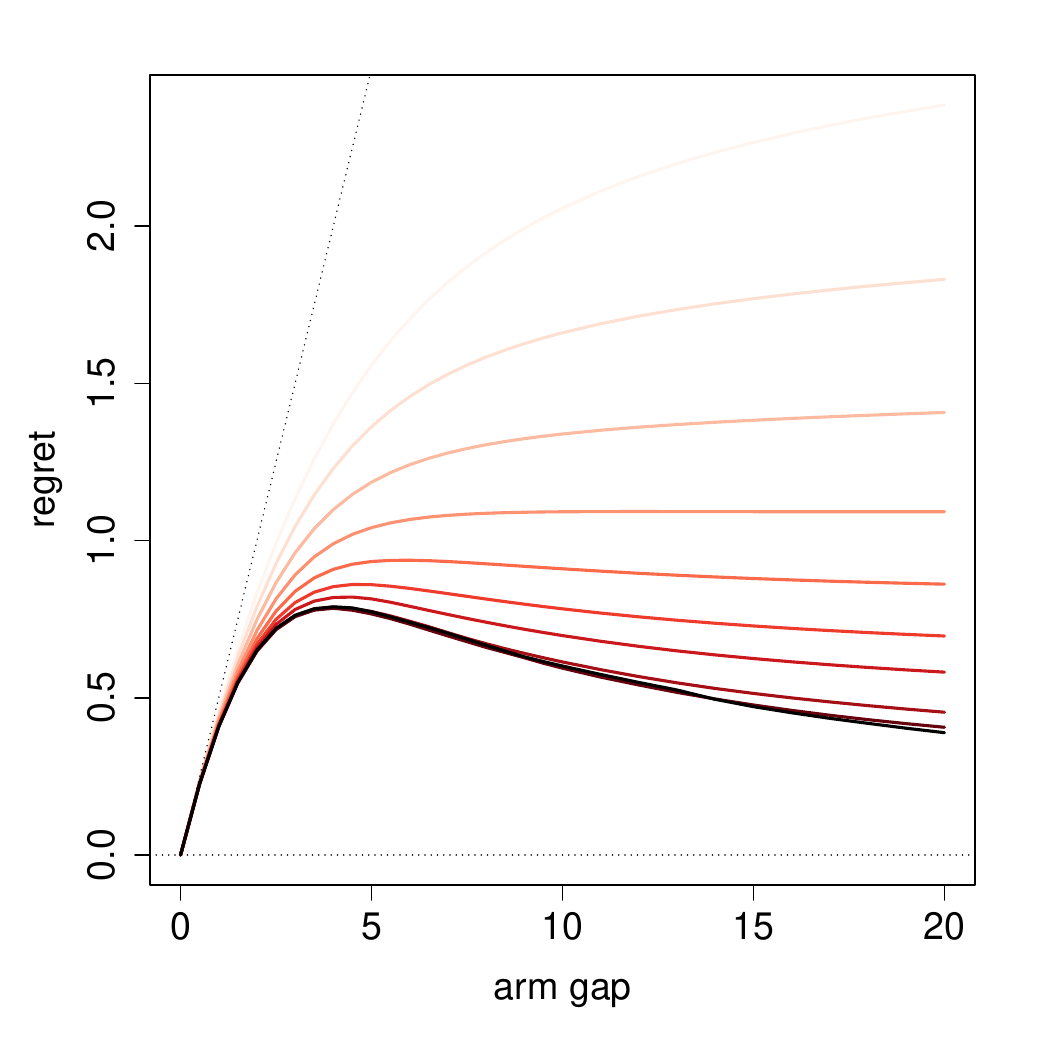} &
		\includegraphics[width=0.45\textwidth]{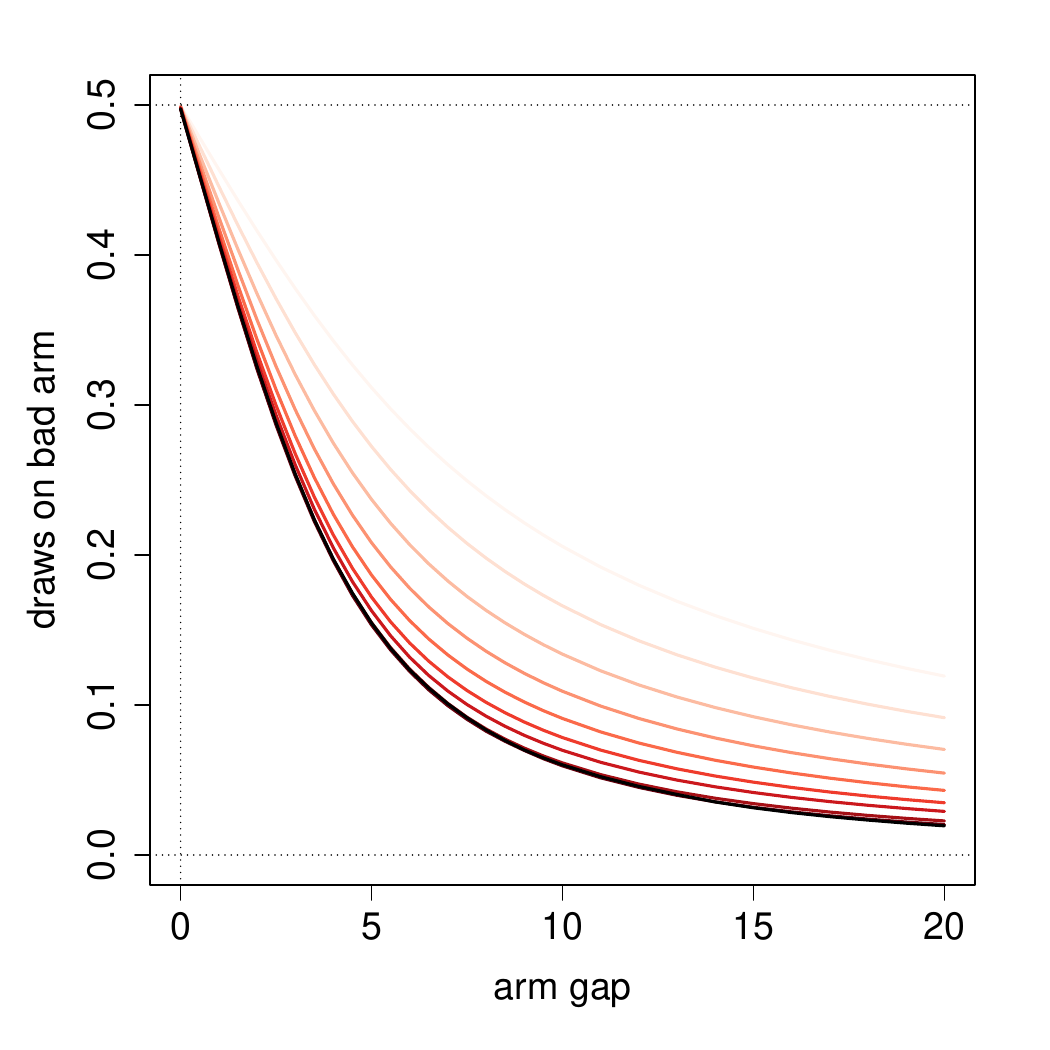}
	\end{tabular}
	\caption{Regret profile for two-armed Thompson sampling, for $c = 1$, 1/2, 1/4, 1/8, 1/16, 1/32, 1/64, 1/256, 1/1024,
		and finally $c = 0$. We use $\sigma^2 = 1$ throughout.
		The left panel shows expected regret, while the right panel shows $\EE{Q_1}$.
		The curves with positive values are shown in hues of red with
		darker-colored hues corresponding to smaller values of $c$, while $c = 0$ is shown in black. The
		algorithm \eqref{eq:TS2_finite} is translation-invariant and symmetric in its treatment of the arms, so regret only
		depends on the scaled arm gap $\delta = \abs{\mu_1 - \mu_2}$.}
	\label{fig:two_arm_profile}
\end{figure}

Having verified that undersmoothed Thompson sampling performs well in the
super diffusive limit, we next use the diffusion limit \rev{that arise in the weak signal regime} to numerically derive an exact, instance-specific
characterization of the mean scaled regret of Thompson sampling, i.e., the limit of
$\sqrt{n} \EE{R^n}$. The result, shown in Figure \ref{fig:two_arm_profile}, gives
us a sharp picture of how the  of Thompson sampling varies
with the arm gap  in a sequential experiment, and also helps us understand the
effect of the prior choice $\nu^n$ on performance.

This analysis immediately reveals an intriguing relation between the regularization parameter $c$ and regret.
The Bayesian heuristic behind Thompson sampling appears to be mostly uninformative about which choices
of $c$ will perform well in terms of regret; and instead, in the region and to the resolution displayed, regret
appears instance-wise monotone increasing in $c$ (i.e., for each arm gap $\delta$ and each pair of regularization
parameters $c_1 < c_2$ displayed, we obtain lower regret with $c_1$).
We also see that regret in fact converges as $c \rightarrow 0$,
and that undersmoothed Thompson sampling appears instance-wise optimal. Thus, \rev{our analysis based on the diffusion limit suggests} 
that if a practitioner wants to use two-armed Thompson sampling and has prior information about the size of the
gap between the two arms, then they may actually be best off ignoring this information in choosing the prior for
Thompson sampling and instead using an undersmoothed specification.

\renewcommand{\arraystretch}{1.6}
\begin{table}[t]
	\resizebox{\textwidth}{!}{\begin{tabular}{|c|ccc|}
			\hline
			& Regret (original) & Scaled regret & Reference \\ \hline
			UCB & $3\Delta_n + 16 \log(n) / \Delta_n$ & $+\infty$ &   LS20, \S 7.1\\
			Thompson  &  ${\log(n) \Delta_n / }{D(\mu^n_{1}, \, \mu^n_{2})}$  & $+\infty$ &  AG17  \\
			MOSS &  $39 \sqrt{2n} + \Delta_n$ & $55.2$ &  LS20, \S9.1 \\ 
			Impr.~UCB & $ \min \cb{\Delta_n + 32 \frac{\log((\Delta^n)^2 n)}{\Delta^n} + \frac{96}{\Delta^n},  n \Delta^n }$ &  $\min \cb{64 \frac{\log(\Delta)}{\Delta}+ \frac{96}{\Delta},   \Delta }$ & AO10\\
			Oracle ETC &$\min \cb{\frac{4}{\Delta^n}\p{1+  \log\p{ \frac{n(\Delta^n)^2 }{4}} ^+} , n \Delta^n }$ & $\min\cb{\frac{4}{\Delta}\p{1+ \log (\Delta^2/4)^+}, \Delta}$& LS20, \S6.1\\ 
			\hline
	\end{tabular}}
	\caption{Comparison with existing finite-time instance-dependent regret bounds. The column of scaled
		regret corresponds to what the bound would become under \rev{weak signal} scaling, where
		$ \mu_k^n = \mu_k / \sqrt{n} $ and $ \Delta^n = \Delta /\sqrt{n} $, with $ n\to \infty $. The algorithms
		under consideration are the upper confidence bound (UCB) algorithm of \citet{lai1985asymptotically},
		Thompson sampling, the Minimax Optimal Strategy in the Stochastic case (MOSS) from \citet{audibert2009minimax},
		improved UCB by \citet{auer2010ucb}, and an oracle explore-then-commit (ETC) baseline that takes uniformly random actions
		up to a deterministic time chosen using a-priori knowledge of the effect size $\abs{\delta}$, and then commits to the
		most promising arm for the rest of time. The specific bounds are as reported in \citet[AG17]{agrawal2017near},
		\citet[AO10]{auer2010ucb} and \citet[LS20]{lattimore2020bandit}.}
	\label{tab:bounds}
\end{table}
\renewcommand{\arraystretch}{1}

\begin{figure}[t]
	\centering
	\includegraphics[scale=.25]{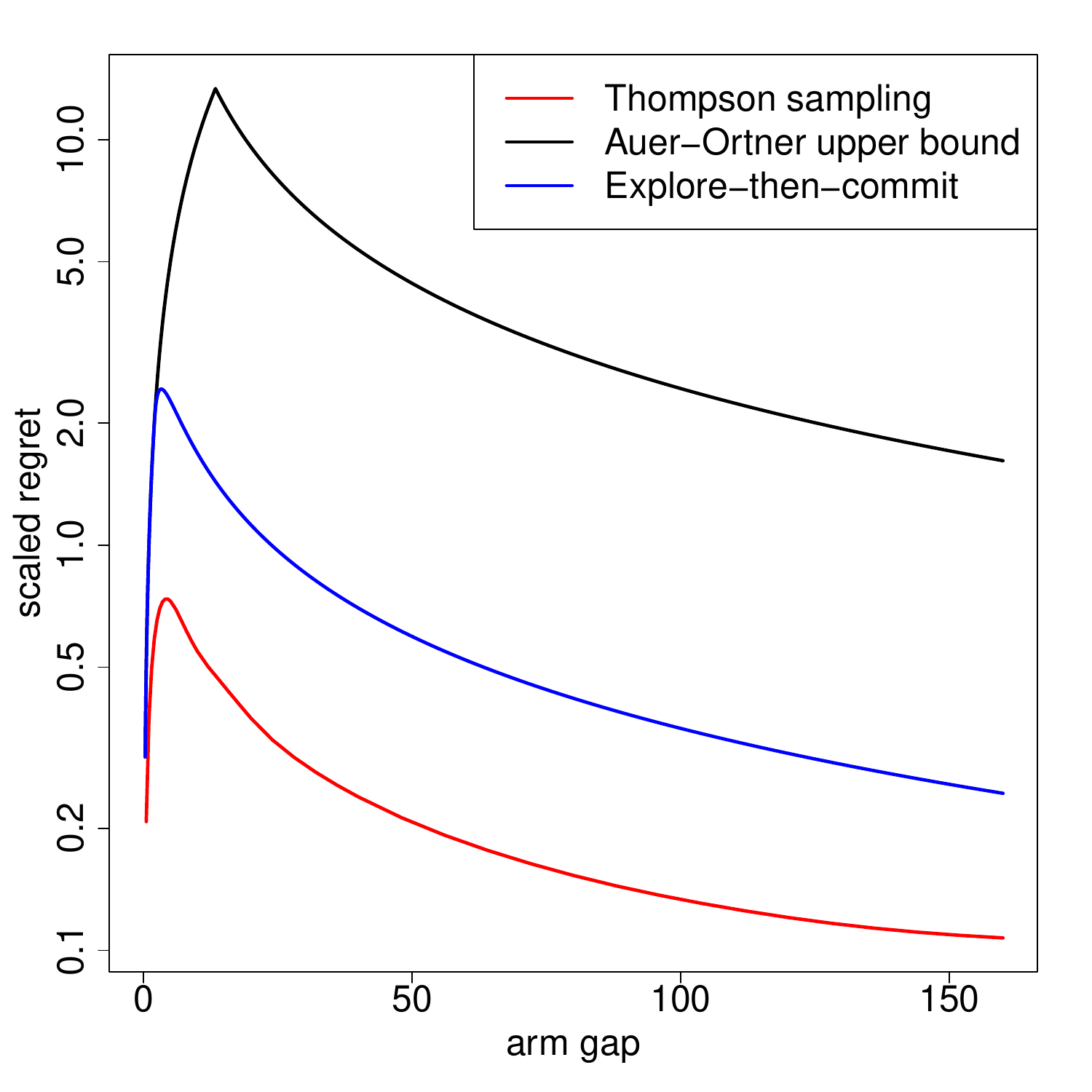}
	\caption{Comparison between the scaled  regret \rev{under the weak signal scaling} for undersmoothed Thompson sampling and existing bounds in Table \ref{tab:bounds}, as function of the arm gap. Here, we have a two-armed bandit with $ \sigma=1 $. We plot all bounds with a finite scaled regret, with the exception of MOSS, which has a constant value of 55.2 in this case.  }
	\label{fig:bounds_comp}
\end{figure}

Next, we can move beyond internal comparisons of different prior choices for Thompson sampling, and
use our diffusion limit to compare large sample behavior of (undersmoothed, translation-invariant)
Thompson sampling to other available baselines.
In Table \ref{tab:bounds}, we collect a number of state-of-the-art finite-sample bounds for two-armed bandits
(we do not have access to exact diffusion-scale regret limits for these methods). We report results
both in the original finite-sample form of the bound, and the (scaled) limit of the
bound under the \rev{weak signal} regime from Definition
\ref{def:reward_diff_scale}.

A first interesting finding is that many available regret bounds are in fact vacuous in the \rev{weak signal} regime, and do not provide any meaningful control on regret.
For example, while we know from our \rev{analysis of the diffusion limit} that Thompson
sampling gets bounded (and in fact quite good) regret under \rev{the weak signal scaling},
the strongest available finite-sample regret bound for
Thompson sampling, due to \citet{agrawal2017near}, diverges under \rev{the weak signal scaling}. 
The only upper bound given in Table \ref{tab:bounds} that both remains
finite under \rev{the weak signal scaling} and has meaningful instance-dependent
behavior (i.e., that improves as the scaled arm gap $\delta$ gets large) is the bound of \citet{auer2010ucb} for improved UCB.
We also note that the oracle explore-then-commit gets good instance-dependent behavior; however, this algorithm relies
on a-priori knowledge of the arm gap size $\abs{\delta}$, and so it is not a feasible baseline.
Thus, any non-trivial bounds obtained via a diffusion-based analysis are likely to be stronger then
many available bounds in regimes where they apply.
Overall, this highlights the insight that the moderate data regime captured by the \rev{weak signal}
asymptotics is in fact a challenging one, where the sample size and arm gap scale together
in such a way that only the sharpest analyses yield non-trivial guarantees on performance.

Finally, in cases where the bounds from Table \ref{tab:bounds} are not vacuous \rev{under the weak signal scaling},
Figure \ref{fig:bounds_comp} compares our limiting regret to the relevant upper bounds.
Interestingly, we see that the limiting regret we obtain for Thompson sampling
is much lower than any of the available upper bounds, and Thompson sampling even outperforms the oracle explore-then-commit algorithm where the duration of exploration is optimized with a-priori knowledge of the effect size and time horizon. These examples suggest that existing finite-sample instance-dependent regret upper bounds could still be improved substantially, possibly by leveraging the \rev{weak signal} asymptotics advanced in this work.

\section{The (In)stability of Sensitive Sampling Functions}
\label{sec:instability}

The diffusion limits \rev{that arise from the weak signal scaling} also allow us to conduct refined performance analysis that goes beyond mean rewards. In this section, we use the  distributional characterization of the diffusion limit to unearth some interesting instability properties of a certain family of ``sensitive'' sampling functions, where the sampling probability of choosing the optimal arm can swing wildly between very large to the very small. 
This family also includes Thompson sampling as a special case.  Both our numerical results and the theoretical analysis in the preceding sections point to the fact that undersmoothed
Thompson sampling ($c=0$) yields far superior total regret than its smoothed counterpart. 
However, this performance improvement from undersmoothing, as it turns out, does not come for free. Although undersmoothed Thompson sampling identifies and
focuses on the correct best action often enough to achieve low average regret, it is also liable to fail completely and double
down on a bad arm.

\begin{figure}[t]
	\centering
	\includegraphics[width=0.7\textwidth]{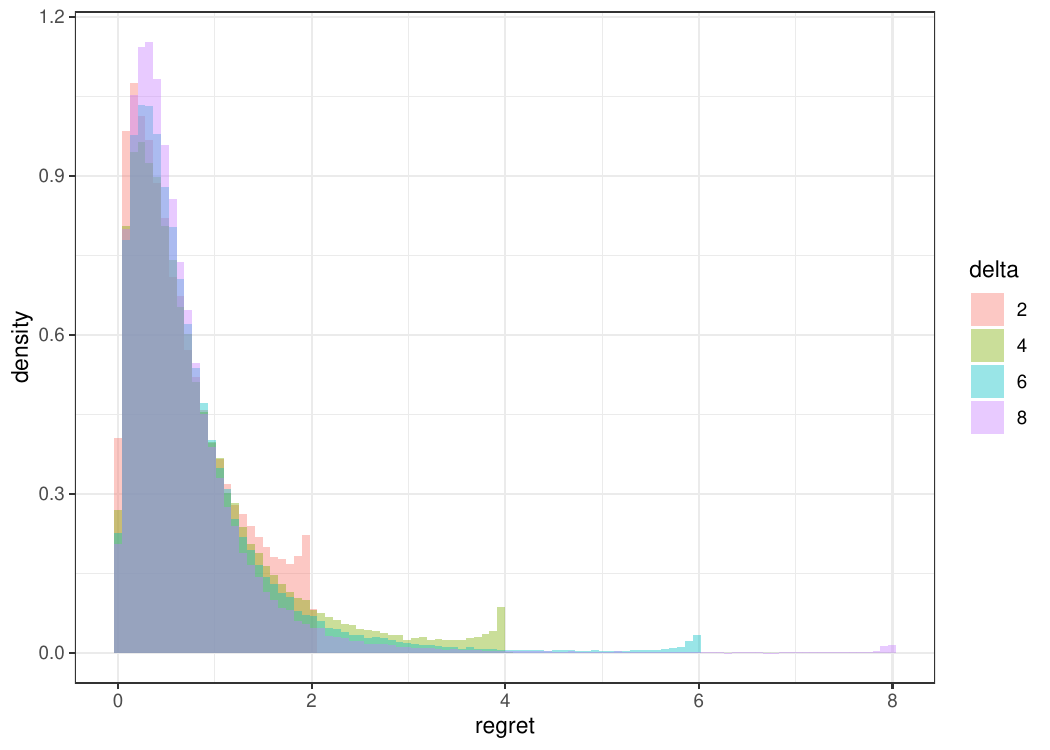}
	\caption{Distribution of the (scaled) regret for two-armed Thompson sampling in the undersmoothed
		regime (i.e., with $c = 0$), as a function of (scaled) arm gap $\delta$. The histograms are aggregated
		over 100,000 realization of the limiting stochastic differential equation.}
	\label{fig:regret_distr}
\end{figure}

As a first lens on the instability of Thompson sampling, Figure \ref{fig:regret_distr} displays the distribution
of regret for undersmoothed two-armed Thompson sampling for a variety of arm gaps $\delta$. Interestingly,
we see that for all considered values of $\delta$, the distribution of regret is noticeably not unimodal. Rather,
there is a primary mode corresponding to the bulk of realizations where Thompson sampling gets reasonably
low regret, but there is also a second mode near $R \approx \delta$. Recall that, if $\mu_1 > \mu_2$, then regret
measures the frequency of draws on the second arm $R = \abs{\delta} Q_{2,1}$, and in particular
$\abs{R} \leq \abs{\delta}$ almost surely. Thus realizations of Thompson sampling with $R \approx \delta$ correspond
to cases where the algorithm almost immediately settled on the bad arm, and never really even gave the good arm
a chance.

To probe deeper into the underlying cause of this instability phenomenon, we will analyze a family of simplified, \textit{one-armed} adaptive experiments. In this setting, the agent is tasked with choosing between an uncertain arm (arm 1), with unknown mean reward $\mu^n$ and known reward variance $ (\sigma^n)^2 $, and a default option (arm 2), with mean reward $\mu_0$. Without loss of generality, we will assume that the default mean reward is equal to $0$. 
We are interested in the \rev{weak signal} regime where  
$\limn \sqrt{n} \mu^n = \mu $ for some fixed $\mu \in \RR$ while $\limn \sigma^n = \sigma$ remains constant. We will refer to $ \mu $ as the {effect size} of the experiment, which can also be thought of as the diffusion-scaled arm gap between the unknown and default arms.

\begin{exam}[One-Armed Thompson sampling]
	\label{ex:onearmTS}
As an example of a one-armed experiment, we can apply the Thompson sampling algorithm described  \eqref{eq:TS1_finite}, where the default option would have a posterior mean reward distribution that is  concentrated at zero at all times. We assume that 
	the agent takes the reward distribution of the uncertain arm $P^n$ to be a Gaussian distribution with (unknown) mean $\mu^n$ and (known) variance $\sigma^2$,
	and sets $G_0^n$ to be a Gaussian prior on $\mu^n$ with mean 0. Thus, writing $Q_i$ for the total number of samples we have drawn the uncertain arm by period $ i $ and $ S_i $ for the sum of these samples, we have
	\begin{equation}
		\label{eq:TS_onearm}
		\begin{split}
			&\mu^n \cond G_i^n \sim \nn\p{\frac{\sigma^{-2}S_i}{\sigma^{-2} Q_i+ \p{\nu^n}^{-2}}, \ \ \frac{1}{\sigma^{-2} Q_i+ \p{\nu^n}^{-2}}}, \\
			&\pi_i =\Phi\p{\frac{\sigma^{-2} S_i}{\sqrt{\sigma^{-2}  Q_i + \p{\nu^n}^{-2}}}},
		\end{split}
	\end{equation}
	where $(\nu^n)^2$ is the prior variance and $\Phi$ is the standard Gaussian cumulative distribution function. In particular, we will assume that $ 		\limn (\nu^n)^{-2}/n = c, \  c \geq 0.  $ As before, this Thompson sampling algorithm is considered {smoothed} if $ c>0 $, and {undersmoothed}, if $ c=0 $. With this choice of prior variance, we obtain a limiting sampling probability for the uncertain arm in the \rev{weak signal} regime: 
	\begin{equation}\label{key}
		\psi_1(q,s) = \Phi\p{\frac{s}{\sigma \sqrt{q + \sigma^2c}}}. 
	\end{equation} 
%
	
	\end{exam}

We now formulate the following notion of sensitive sampling functions within the family of one-armed experiments. 

\begin{defi}
	\label{def:sensitive}
	Consider a one-armed experiment. We say that a limiting sampling function, $ \psi(q,s) $, is sensitive, if 
	\begin{enumerate}
		\item for any sequence $ \{ (q_n,s_n) \}_{n \in \N} $ such that $ \lim_{n\to \infty} s_n /\sqrt{q_n} = \infty  $, we have 
		\begin{equation}\label{key}
			\lim_{n \to \infty  }		\psi(q_n,s_n) =   1; 
		\end{equation}
		\item for any  sequence $ \{ (q_n,s_n) \}_{n \in \N} $ such that $ \lim_{n\to \infty} s_n /\sqrt{q_n} = - \infty  $, we have
		\begin{equation}\label{key}
			\lim_{ n \to \infty}	\psi(q_n,s_n) =   0. 
		\end{equation}
	\end{enumerate}
\end{defi}

In words, a sampling function is sensitive if the policy would assign probability either one or zero to drawing from the uncertain distribution, in the limit as $ s/\sqrt{q} $ tends to positive or negative infinity, respectively. The quantity $ s/\sqrt{q} $  arises, for instance, under a Bayesian interpretation of the uncertain mean under a diffusive prior. Consider the Gaussian one-armed experiment described in Example \ref{ex:onearmTS}. 
 In this case, when the prior variance $ (\nu^n)^2 $ is very large,  the agent's posterior belief on the effect size $ \mu^n $ admits a Gaussian distribution  in period $ i $  with mean $ S_i/Q_i $ and standard deviation $ 1/\sqrt{Q_i} $. The value $ S_i/\sqrt{Q_i} $ is therefore the z-score associated with $ \mu^n$.

The following theorem characterizes the the evolution of the sampling probabilities $\pi_t$ over time under any sensitive sampling function in a one-armed experiment. 
Qualitatively, these sampling probabilities $\pi_t$ can be thought of as the agent's subjective preference about which arm is the best. In this case, we find that  a sensitive sampling function will always lead the \rev{agent to committing} to the ``wrong'' superior arm with arbitrarily high intensity at some point during the sampling process, and that this is true no matter how large the magnitude of the actual effect size. 

\begin{theo}
	\label{theo:sample_path}
	Fix a one-armed experiment and a continuous limiting sampling function $\psi$ that is also sensitive as per Definition \ref{def:sensitive}.  Consider the sampling probability \smash{$\pi_t = \psi(Q_t, S_t)$} associated with the diffusion limit under $\psi$, where $ Q $ and $ S $ correspond to the state of the uncertain arm. Then,  fixing any effect size $\mu \in \RR$ and confidence level $ \eta \in (0,1) $, we have, for all $ \epsilon \in (0,1) $, 
	\begin{equation}
		\PP{ \sup_{ t \in [0,\epsilon)} \pi_t \geq 1 - \eta} = \PP{ \inf_{ t\in [0,\epsilon)} \pi_t \leq \eta} = 1 . 
	\end{equation}
\end{theo}

%
%

\begin{figure}[t]
	\begin{tabular}{cc}
		\includegraphics[width=0.45\textwidth]{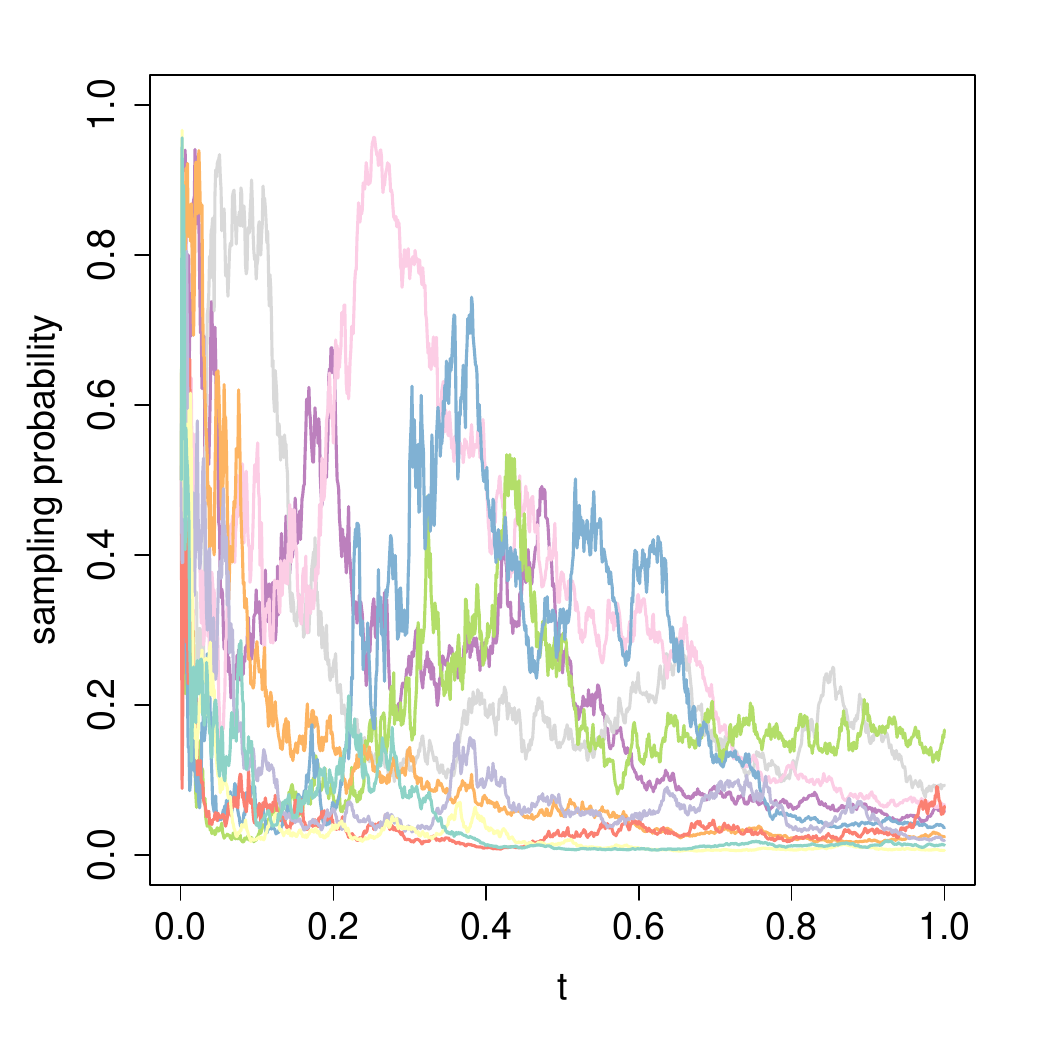} &
		\includegraphics[width=0.45\textwidth]{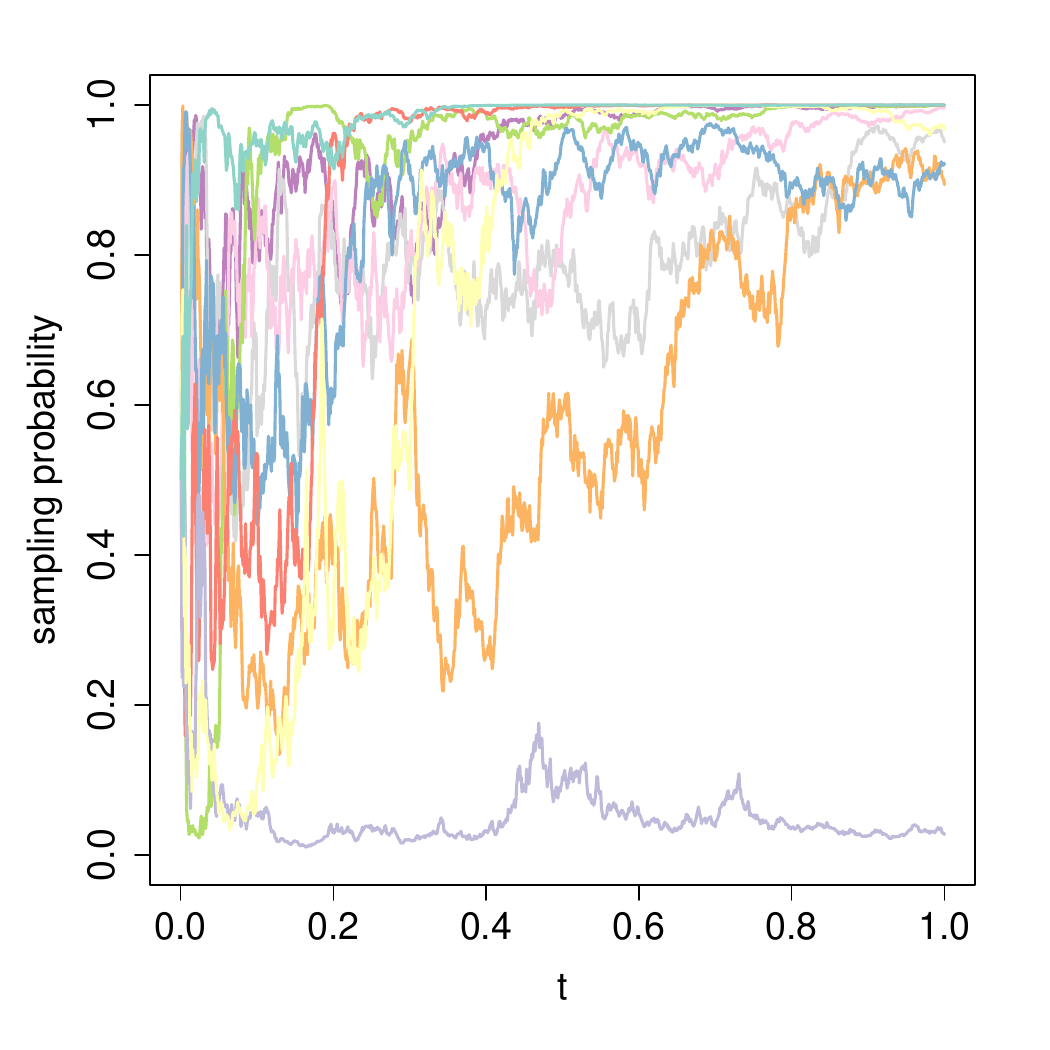} \\
		$\mu = -3$ & $\mu = 3$
	\end{tabular}
	\caption{Sample paths of the sampling probability $\pi_t$ in one-armed Thompson sampling as defined in \eqref{eq:TS_onearm},
		in the undersmoothed regime (i.e., with $c = 0$).}
	\label{fig:one_arm_TS_path}
\end{figure}

The key idea in the proof of Theorem \ref{theo:sample_path}, given in Appendix \ref{proof:theo:sample_path}, is to use the time-changed form
of the diffusion limit given in Theorems \ref{theo:timechange}. In particular, we show that under the weak signal scaling, the quantity $S_t/\sqrt{Q_t}$ can be written via the time-change formula
such that a leading determinant of its value is $ W_{Q_t}/\sqrt{Q_t} $, 
where $W_q$ is a standard Brownian motion. By the law of iterated logarithm, we know that
$W_q / \sqrt{q}$ will get arbitrarily close to $\pm \infty$ as $q \downarrow 0$, which, using the definition of a sensitive sampling function,  we can use
to show that $\pi_t$ will in turn spend time arbitrarily close to $0$ and $1$ as $t \downarrow 0$. 
A similar result also holds for the two-armed case, although its proof is not as immediate so we do not provide it here.

In  words, Theorem \ref{theo:sample_path} shows that an agent using a sensitive sampling function will almost always at some early point in time
be arbitrarily convinced about $\mu$ having the wrong sign; and this holds no matter how large the effect size
$\abs{\mu}$ really is. 
One consequence of this finding is that it challenges a perspective that would take Thompson sampling literally
as a principled Bayesian algorithm (since in this case we'd expect belief distributions to follow martingale
updates), and instead highlights that Thompson sampling has subtle and unexpected behaviors that
can only be elucidated via dedicated methods. 

\section{Discussion}
\label{sec:discuss}

In this paper, we introduced a \rev{weak signal asymptotic  regime} under which sequentially randomized experiments converge to a diffusion limit. In particular, the limit cumulative reward is obtained by applying a random time change to a constant-drift Brownian motion, where the time change is in turn given by cumulative sampling probabilities (Theorem \ref{theo:timechange}). We then applied this result to derive sharp insights about the behavior of Lipschitz-continuous sampling functions and undersmoothed two-armed Thompson sampling. Notably, we demonstrate that Lipschitz-continuous sampling always leads to non-vanishing, and thus sub-optimal, regret when the arm gaps are large. In contrast, we show that  using an asymptotically uninformative prior with two-armed Thompson sampling leads to  desirable regrets across all arm gaps in the \rev{weak signal} regime,  and especially when the arm gap is large,  where we show the regret is nearly instance-optimal. This latter result is one of the first recent attempts at  formally justifying the folklore wisdom of using diffusive priors in Thompson sampling \citep[e.g.,][]{liu2022gaussian}, where, for instance, it is recommended to set the prior variance to a fixed large constant \citep{agrawal2017near}. 

A first class of natural follow-up questions is in seeing whether diffusion limits hold for broader classes of sequential experiments.  The current diffusion analysis leaves out those experiments that have a discontinuous sampling function, which include, for instance, UCB-type experiments.  A key technical challenge in analyzing discontinuous sampling functions is that the number of arm pulls under them can be highly erratic,  which makes showing the convergence to a diffusion limit more challenging.  After we posted a version of this paper online,  the recent work by \cite{kalvit2021closer} made progress towards this question by showing a diffusion limit for a two-armed UCB algorithm.  Their analysis relies on a refined regularity property of UCB where they show that both arms will be pulled approximately an equal number of times as the number of samples grows. It remains unclear whether such analysis can be generalized beyond UCB with two arms, or to other  discontinuous sampling functions.  Furthermore, it would be interesting to see whether our results can be extended to the case of contextual bandits, or to bandit problems with continuous action spaces that
arise, e.g., with pricing. Finally, another practical question is whether the approach used here can be used to build confidence
intervals using data from sequential experiments, thus adding to the line of work pursued by \citet{hadad2019confidence},
\citet{howard2018uniform}, \citet{zhang2020inference}, and others.

Another set of questions could aim to characterize the rate of convergence to the diffusion limit, which we do not consider in the current paper. Such results can be useful in quantifying the quality of  \rev{weak signal} asymptotics for finite-sample experiments. It would be interesting, for instance, to understand how the quality of approximation depends on the smoothness of the sampling function. 

Further along, a potentially interesting avenue for investigation is whether the diffusion limit derived here is useful
in understanding human---as opposed to algorithmic---learning. Throughout this paper, we have considered Thompson
sampling and related algorithms as a class of sequential experiments designed by an investigator, and have discussed how
different design choices (e.g., around smoothing) affect the performance of the learning algorithms.
However, following \citet{erev1998predicting} and \cite{xu2020reinforcement}, we could alternatively use sequentially randomized Markov experiments as
models for how humans (or human communities) learn over time, and use our results to make qualitative predictions about
their behavior.
For example, it may be of interest to examine Thompson sampling as a model for how a scientific community
collects and assimilates knowledge about different medical treatments. Qualitatively, this would correspond to a hypothesis that the number of
scientists investigating any specific treatment should be proportional to the consensus beliefs that the treatment is best given the
available evidence at the time (i.e., that scientists at any given time prioritize investigating the treatments that appear most promising).
In this case, our Theorem \ref{theo:sample_path} would provide conditions under
which we predict consensus beliefs to first temporarily concentrate around sub-optimal treatments before eventually reaching the
truth. More broadly, the many unintuitive phenomena arising from our diffusion limits could yield a number of valuable insights
on how we collect and process information.

\bibliographystyle{plainnat}
\bibliography{references}

\newpage

\begin{center}
	\Large Supplemental Material for\\ \textit{Weak Signal Asymptotics for Sequentially Randomized Experiments}
	
	\vspace{20pt} 
\end{center}

\ifx \useplain\undefined
\begin{APPENDICES}
	\else
	\begin{appendix}
		\fi
		
		\normalsize

\section{Proofs of Main Results}
\label{sec:proofs}

\textit{Terminology}: We will use the term  ``almost all sample paths'' throughout the proofs to refer to sample paths of the Brownian motion that belong to a set of probability measure one. We will also use the following asymptotic notation: $ f(x) \ll g(x)$ indicates that $ f(x) / g(x) \to 0 $, similarly for $ f (x) \gg g (x) $. 

\subsection{Proof of Theorem \ref{theo:SDE}}
\label{sec:proof:theo:SDE}

\bpf The proof is based on the martingale framework of  \citet{stroock2007multidimensional}. Let us first introduce some notation to streamline the presentation of the proof.  Define $Z_t = (\bQ_t, \bS_t)$. Denote by $\calI_S$ and $\calI_Q$ the indices in $Z_t$ corresponding to the coordinates of $S$ and $Q$, respectively. Both sets are understood to be an ordered set of $K$ elements, where the subscript is for distinguishing whether the ordering is applied to $S$ versus $Q$.  For $z \in \rp^K\times \R^K$, define the functions $(b_k)_{k \in \calI_Q \cup \calI_S}$, 
\begin{align}
	b_k(z) =  &   \begin{cases} 
		\psi_k(z), & k \in \calI_Q, \\
		\psi_k(z)\mu_k, & k \in \calI_S. 
	\end{cases}
	\label{eq:bi}
\end{align}
For $1\leq k,l \leq K$, define 
\begin{align}
	\eta_{k,l}(z) = & 
	\begin{cases}
		\sqrt{\psi_k(z)}\sigma_k, & \text{if $k=l \in \calI_S$,}\\
		0, & \text{otherwise.}
	\end{cases} 
	\label{eq:etaij}  
\end{align}
Then, the It\^o diffusion SDE in \eqref{eq:diffLim} can be written more compactly as
\begin{equation}
	dZ_t = b(Z_t) dt + \eta(Z_t) dB_t, \quad t \in [0,1], 
	\label{eq:ItoDiff}
\end{equation}
with $Z_0 = 0$.

Next, we briefly review the relevant results of the Stroock and Varadhan program.  Fix $d\in \N$. Let $(Z^n_i)_{i\in \N}$ be a sequence of time-homogeneous Markov chains taking values in $\R^d$, indexed by $n\in \N$.  Denote by $\Pi^n$ the transition kernel of $Z^n$: 
\begin{equation}
\Pi^n(z, A) = \pb\p{Z^n_{i+1} \in A \cond Z^n_i = z}, \quad z \in \R^d, A \subseteq \R^d.
\end{equation}
Let $\bZ^n_t$ be the piece-wise linear interpolation of $Z^n_{nt}$: 
\begin{equation}
\bZ^n_t = (1-tn+ \lfloor tn \rfloor) Z^n_{\lfloor tn \rfloor} + (tn - \lfloor tn \rfloor ) Z^n_{\lfloor tn \rfloor +1},  \quad t \in [0,1].  
\end{equation}
Define $K^n(z,A)$ to be the scaled transition kernel: 
\begin{equation}
K^n(z,A) = n \Pi^n(z,A). 
\end{equation}
Finally, define the functions 
\begin{align}
a^n_{k,l}(z) = & \int_{x: |z-x|\leq 1} (x_k-z_k)(x_l-z_l)K^n(z,dx), \nln
b^n_{k}(z) = & \int_{x: |z-x|\leq 1} (x_k-z_k)K^n(z,dx), \nln
\Delta^n_{\epsilon}(z) = &  K^n(z, \{x: |x-z|> \epsilon\}). \nnb
\end{align}
We will use the following  result. A proof of the theorem can be found in \citet[Chapter 11]{stroock2007multidimensional} or \citet[Chapter 8]{durrett1996stochastic}. For conditions that ensure the uniqueness and existence of the It\^o diffusion  \eqref{eq:ito_sde}, see \citet[Chapter 5, Theorem 2.9]{karatzas2005brownian}.
\begin{theo} 
\label{theo:StroVar}
Fix $d$. Let $\{a_{k,l}\}_{1\leq k,l \leq d}$ and $\{b_k\}_{1\leq k \leq d}$ be bounded Lipschitz-continuous functions from $\R^d$ to $\R$.  Suppose that for all $k,l \in \{1, \ldots, d\}$ and $\epsilon, R >0$
\begin{align}
& \lim_{n \to \infty} \sup_{z: |z|<R} |a^n_{k,l}(z) - a_{k,l}(z)| = 0, \label{eq:an} \\
& \lim_{n \to \infty} \sup_{z: |z|<R} |b^n_{k}(z) - b_{k}(z)| = 0, \label{eq:bn} \\
& \lim_{n \to \infty} \sup_{z: |z|<R}  \Delta^n_\epsilon(z) = 0. \label{eq:deltan}
\end{align}
If $\bZ^n_0 \to z_0$ as $n\to \infty$, then $(\bZ^n_t)_{t\in [0,1]}$ converges weakly  in $\calC$ to the unique solution to the stochastic differential equation
\begin{equation}
dZ_t = b(Z_t) dt + \eta(Z_t) dB_t, \label{eq:ito_sde}
\end{equation}
where  $Z_0 = z_0$ and $\{\eta_{k,l}\}_{1\leq k,l\leq d}$ are dispersion functions such that $a
(z) = \eta(z) \eta^\intercal(z) $.\footnote{The decomposition from $a$ to $\eta$ is unique only up to rotation. However, the resulting stochastic differential equation \textit{is} uniquely defined by $a$. This is because the distribution of the standard Brownian motion is invariant under rotation, and hence any valid decomposition would lead to the same stochastic differential equation.} 

\end{theo}

We are now ready to prove Theorem \ref{theo:SDE}.  We will use the compact representation of the SDE given in \eqref{eq:ItoDiff}, with $\bZ^n = (\bQ^n, \bS^n)$, $Z_t = (Q_t, S_t)$ and $b$ and $\eta$ defined as in \eqref{eq:bi} and \eqref{eq:etaij}, respectively. To prove the convergence of $\bZ^n$ to the suitable diffusion limit, we will evoke Theorem \ref{theo:StroVar}  (here $d \leftrightarrow 2K$). It suffices to verify the convergence of the corresponding generators  in \eqref{eq:an} through \eqref{eq:bn}. To start, the following technical lemma \cite[Section 8.8]{durrett1996stochastic} will simplify the task of proving convergence by removing the need of truncation in the integral; the proof is given in Section \ref{app:lem:no_trunc_int}. Define 
\begin{align} 
m^n_{p}(z) =  \int  |x-z|^p K^n(z,dx), \ \ \ \ \ \tilde a^n_{k,l}(z) = & \int (x_k - z_k)(x_l-z_l)K^n(z,dx), \nln
\tilde b^n_{k}(z) =  \int (x_k-z_k)K^n(z,dx). \nnb
\end{align}
\begin{lemm}
\label{lem:no_trunc_int}
Fix $p \geq 2$ and suppose that for all $R < \infty$, 
\begin{align}
& \lim_{n \to \infty} \sup_{z: |z|<R}  m^n_p(z) = 0.  \label{eq:mnp} \\
& \lim_{n \to \infty} \sup_{z: |z|<R} |\tilde a^n_{k,l}(z) - a_{k,l}(z)| = 0, \label{eq:tild_an} \\
& \lim_{n \to \infty} \sup_{z: |z|<R} |\tilde b^n_{k}(z) - b_{k}(z)| = 0, \label{eq:tild_bn}
\end{align}
Then, the convergence in \eqref{eq:an} through \eqref{eq:deltan} holds. 
\end{lemm}

In what follows, we will use $z = (q,s)$ to denote a specific state of the Markov chain $\bZ^n$.  The transition kernel of the pre-limit chain $\bZ^n$ can be written as 
\begin{equation}
\Pi^n((q,s) , (q+{e_k}/{n},s+e_k ds/\sqrt{n}) ) =   \bar\psi^n_k(q,s) P^n	_k(ds), \quad k = 1, \ldots, K, 
\end{equation}
and zero elsewhere, where $e_k \in \{0,1\}^K$ is the unit vector where the $k$th entry is equal to $1$ and all other entries are $0$, and $\{P^n_k\}_{k = 1,\ldots, K}$ are the  reward  probability measures. Define $K^n(z, A) = n \Pi^n(z, A)$.

We next define the limiting functions $a$ and $b$. The function $b$ is defined as in \eqref{eq:bi}: 
\begin{align}
b_k(z) =  &   \begin{cases}
      \psi_k(z), & k \in \calI_Q, \\
\psi_k(z)\mu_k, & k \in \calI_S, 
    \end{cases}
\end{align}
and we let 
\begin{equation}
a_{ij}(z) = (\eta\eta^\intercal)_{k,l}(z)
\end{equation}
where $\eta$ is defined in \eqref{eq:etaij}. That is, 
\begin{align}
a_{k,l}(z) = & 
    \begin{cases}
      {\psi_k(z)}\sigma_k^2, & \text{if $k=l \in \calI_S$,}\\
      0, & \text{otherwise.}
    \end{cases} 
\end{align}

Fix $R>0$. We show that the corresponding $a^n$ and $b^n$ converge to the functions $a$ and $b$ defined above, uniformly over the compact set $\{z: |z| \leq R\}$. In light of Lemma \ref{lem:no_trunc_int}, it suffices to verify the convergence in \eqref{eq:mnp} through \eqref{eq:tild_bn} for $p= 4$. Starting with \eqref{eq:mnp}, we have that
\begin{align}
m^n_4(z) = & \int |z'-z|^4 n\Pi^n(z, dz') \nln
= & \sum_{k=1}^K n\bar\psi^n_k(z) \int_{w\in \R} \p{\frac{1}{n^2} + \frac{w^2}{n}}^{2} P^n_k(dw) \nln
\leq & \sum_{k=1}^K n\bar\psi^n_k(z) \p{\frac{2}{n^4}+ \frac{1}{n^2}\int_{w\in \R}w^4 P^n_k(dw)}\nln
= & \frac{2}{n} + \frac{1}{n}\EE[Z \sim P^n_k]{Z^4} \nln
\stackrel{n \to \infty} {\longrightarrow} & 0, 
\end{align}
as $n\to 0$, uniformly over all $z$, where the last step follows from the assumption that the reward distributions admit bounded fourth moments.  This shows \eqref{eq:mnp}. 

For the drift term $b$, we consider the following two cases; together, they prove \eqref{eq:tild_bn}. 

\textit{Case 1, $k\in \calI_Q$.} For all $k\in \calI_Q$, and $n\in \N$, 
\begin{align}
\tilde b^n_{k}(z)  =& \int (q'_k-q_k)K^n(z,dz') \nln
= & \frac{1}{n} (n \bar\psi^n_k(z) )  \nln
\stackrel{n \to \infty} {\longrightarrow}  & \psi_k(z). 
\end{align}

\textit{Case 2, $k\in \calI_S$.} For all  $k \in \calI_S$,  
\begin{align}
\tilde b^n_{k}(z)  :=& \int (s'_k-s_k)K^n(z,dz') \nln
= &  \bar\psi^n_k(z)n \int \frac{w}{\sqrt{n}} P^n_k(dw)  \nln
=&  \bar\psi^n_k(z) \mu_k  \nln
\stackrel{n \to \infty} {\longrightarrow}  & b_k(z). 
\end{align}

For the variance term $a$, we consider the following three cases: 

\textit{Case 1, $k,l \in \calI_Q$.} Note that under the multi-armed bandit model, only one arm can be chosen at each time step. This means that only one coordinate of $\bQ^n$ can be updated at time, immediately implying that for all $n$ and $k,l \in \calI_Q$, $k\neq l$,
\begin{equation}
\tilde a^n_{k,l} (z)= \int (q'_k- q_k)(q'_l-q_l)K^n(z,dz') = 0. 
\end{equation}
For the case $k=l$, we note that for all $k\in \calI_Q$, and all sufficiently large $n$
\begin{equation}
\tilde a^n_{k,k} (z) = \frac{1}{n^2} n \bar\psi^n_k(z)  \stackrel{n \to \infty} {\longrightarrow} 0.  
\end{equation}

\textit{Case 2: $k\in \calI_Q, l \in \calI_S$, or $k \in \calI_S$ and $l\in \calI_Q$}. 
\begin{align}
\tilde a^n_{k, l}(z) = & \int (q'_k- q_k)(q'_l-q_l)K^n(z,dz')  \nln
= & \frac{\bar\psi^n_k(z)}{n} \int  w (n P^n_k(\sqrt{n}dw)) \nln
= & \bar\psi^n_k(z) \EE[Z \sim P^n_k]{Z}  \nln
= & \bar\psi^n_k(z)\mu_k / \sqrt{n} \nln
\stackrel{n \to \infty} {\longrightarrow} & 0. 
\end{align}

\textit{Case 3: $k, l\in \calI_S$}.  This case divides into two further sub-cases. Suppose that $k\neq l$. Similar to the logic in Case 1, because only one coordinate of $\bQ^n$ can be updated at a given time step, we have 
\begin{align}
\tilde a^n_{k, l}(z) = 0, \quad k\neq l. 
\end{align}
Suppose now that $k=l$. We have
\begin{align}
\tilde a^n_{k, l}(z) = & \int (q'_k - q_k)^2 K^n(z,dz')  \nln
= & {\bar\psi^n_k(z)} \int  w^2 (n P^n_k(\sqrt{n}dw)) \nln
= & \bar\psi^n_k(z) \EE[Z \sim P^n_k]{Z^2 }  \nln
\stackrel{n \to \infty} {\longrightarrow} & \psi_k(z)\sigma_k^2 \nln
= &  a_{k,l}(z). 
\end{align}
We note that by assuming that the limiting sampling function $ \psi $ is Lipschitz-continuous, the convergence of $\tilde b^n$, $\tilde a^n$ and $m^n_p$ to their respective limits holds uniformly over compact sets.   We have thus verified the conditions in Lemma \ref{lem:no_trunc_int}, further implying \eqref{eq:an} through \eqref{eq:bn}.  
Note that because $\psi_k$ is bounded and Lipschitz-continuous, so are $a$ and $b$.  This proves the convergence of $\bZ^n$ to the  diffusion limit in $\calC$. 

Finally, to prove the convergence of $\EE{f(\bZ^n_1)} $ to $ \EE{f(Z_t)}$, note that the weak convergence of $\bZ^n$ in $\calC$ implies that the marginal distribution, $\bZ^n_t$ converges weakly to $Z_t$, as $n \to \infty$.  The result then follows immediately from the continuous mapping theorem and the bounded convergence theorem. This completes the proof of Theorem \ref{theo:SDE}. 
\qed

\subsection{Proof of Theorem \ref{theo:timechange}} 
\label{sec:proof:theo:timechange}

\bpf It suffices to show that \eqref{eq:St_timechange} holds.  We will begin with a slightly different, but equivalent, characterization of the pre-limit bandit  dynamics. Consider the $n$th problem instance. Denote by $\tilde Y_{k,j}$ the reward obtained from the $j$th pull of arm $k$. Then, we have that for a fixed $k$, $\tilde Y_{k, \cdot}$ is an i.i.d.~sequence, independent from all other aspects of the system, and 
\begin{equation}
S_{k,i} = \sum_{j=1}^{Q_{k,i}} \tilde Y_{k, j}.
\end{equation}
We can further write 
\begin{equation}
\tilde Y_{k,j} = \mu_k/\sqrt{n} + U_{k,j}, 
\end{equation}
where $U_{k,j}$ is a zero-mean random variable with variance $\sigma_k^2$. 
Define the scaled process: 
\begin{equation}
U^n_{k,i} = \frac{1}{\sqrt{n}}\sum_{j=1}^{i}  U_{k, j}.
\end{equation}
We thus arrive at the following expression for the diffusion-scaled cumulative reward:  
\begin{equation}
\bS^n_{k,t} =  \bQ^n_{k,t}  \mu_k + U^n_{k,n \bQ^n_{k,t}} , \quad i=1, \ldots, n.
\label{eq:Sn_alter}
\end{equation}

Denote by $\bar U^n_{k,t}$ to be the linear interpolation of $ U^n_{k, \lfloor n t \rfloor}  $ for $t \in [0,1]$.  By Donsker's theorem, there exists a $K$-dimensional standard Brownian motion $W$ such that $\bar U^n$ converges to $\{\sigma_k W_{k,\cdot}\}_{k=1, \ldots, K}$ weakly in $\calC$. Evoking Theorem \ref{theo:SDE} and the Skorohod's representation theorem \citep[Theorem 6.7]{billingsley1999convergence}, we may construct a probability space on which the following convergences in $\calC$ occur almost surely: 
\begin{equation}
\bS^n_{k,t} \to S_{k,t}, \, \bQ^n_{k,t} \to Q_{k,t}, \, \bar U^n_t{k,t} \to \sigma_k W_{k,t}, 
\label{eq:as_convg_SQU}
\end{equation}
as $n\to \infty$, where $S$ and $Q$ are diffusion processes satisfying the SDEs in Theorem \ref{theo:SDE}. 

We now combine \eqref{eq:Sn_alter} and \eqref{eq:as_convg_SQU}, along with the fact that $W$ is uniformly continuous in the compact interval $[0,1]$, to conclude that almost surely
\begin{equation}
U^n_{k,n \bQ^n_{k,i}} \to \sigma_k W_{k, Q_{k,t}}, 
\end{equation}
in $\calC$.  This further implies that $S$ also satisfies
\begin{equation}
S_{k,t} =  Q_{k,t} \mu_k + \sigma_k W_{k,Q_{k,t}}, 
\end{equation}
proving our claim.

\qed

\subsection{Proof of Lemma \ref{lemm:tightness}}
\label{app:lemm:tightness}

\bpf  For the first claim, we note that the $ Q^j $'s are $ 1 $-Lipschitz and therefore uniformly equicontinuous over $ [0,1] $. It follows from the Arzela-Ascoli theorem that $ \{Q^j\} $ admits a convergent subsequence $ Q^{j_m} $. The convergence of $ S^{j_m} $ follows directly from its characterization in \eqref{eq:St_timechange} and so does the fact that $ S $ satisfies \eqref{eq:St_timechange}. To show that $ Q $ is a solution to the ODEs, let $ \{Q^{j_m}\}_{m \in \N} $ be the subsequence that converges uniformly to $ Q $. Because the sample path of the Brownian motion $ W$ is uniformly bounded over $ t\in [0,1]$ almost surely, it follows from \eqref{eq:Qt_timechange} that the derivative $  (Q^{j_m})'$ converges uniformly over $ (0,1) $. Combining this fact with the uniform convergence of $ Q^{j_m} $, it follows from the fundamental theorem of calculus that $ Q $ is differentiable over $ (0,1) $, and $ (Q^{j_m})' $ converges uniformly to $Q'$. This proves the claim. 
\qed

\subsection{\rev{Proof of Theorem \ref{theo:conv_triangular_psi}}}
\label{app:theo:conv_triangular_psi}		

\rev{
\bpf By Theorem \ref{theo:timechange}, we have that for every $ j $, $ (\bar Q^{n,j}, \bar S^{n,j} ) $ converges to a unique solution to the ODEs in Theorem \ref{theo:timechange} with sampling function $ \psi^j$, $ (Q^j, S^j) $, as $ n\to \infty $. By Lemma \ref{lemm:tightness}, there exists a limit function $ (Q, S) $ for the sequence $ \{(Q_j, S_j)\}_{j \in \N}$   and furthermore that $ (Q,S) $ is a solution to the ODEs in Theorem \ref{theo:timechange} with sampling function $ \psi $. The above two factors combined imply that there exists a sequence $ \{j_n\}$, with $ \lim_{n\to \infty} j_n = \infty$, such that $ (\bar Q^{n, j_n} , \bar S^{n, j_n} ) $ converges to $ (Q, S) $ as $ n\to \infty $. This proves the claim. 
\qed
}
\subsection{Proof of Theorem \ref{theo:lip_bad}}
\label{app:lip_bad}

\bpf   Denote by $ R_t $ the cumulative  regret up to time $ t $ under the weak signal scaling. We have that for all $ t \in (0,1) $: 
\begin{equation}
	R = R_1 \geq R_t \geq (t-Q_{1,t}) \delta. 
	\label{eq:R_to_Rtdelta}
\end{equation}

Next, we will bound the term $  (t-Q_{1,t}) $ by analyzing its evolution over a small time interval near $ t=0 $. Define $ t_\delta= c/ \delta$, where $ c $  is a constant whose value will be chosen subsequently.  Recall that $ Z_t := (Q_t, S_t) $. By the Lipschitz continuity of $ \psi $, we have that for all $ t  < t_\delta $, and all sufficiently large $ \delta $, the following holds almost surely: 
\begin{align}\label{key}
	\psi_1(Z_t) \leq & \psi_1(0,0)+ L\|Z_t\|  \nln
	\leq &  \tilde p + L(\|Q_t\| +\|S_t\|) \nln
	\sk{a}{\leq} &  \tilde p + L( Kt_\delta +\|S_t\|) \nln
	\sk{b}{\leq} &  \tilde p + LK( t_\delta + (\tilde m_1 + \delta)t_\delta  + {t_\delta}^{2/3}) \nln
	\leq &  \tilde p +2  LK \delta t_\delta \nln
	= &  \tilde p + 2  LKc, 
\end{align}
where step $ (a) $ follows from each coordinate of $ Q_t $ being 1-Lipschitz by construction. Step $ (b) $ follows from the characterization in Theorem \ref{theo:timechange} that $ S_{k,t} =  Q_{k,t} \mu_k + \sigma_k W_{k,Q_{k,t}},$ where, evoking the law of iterated logarithm (Lemma \ref{lemm:LIL}), we conclude that, for all sufficiently large $ \delta $: 
\begin{equation}\label{key}
	\sigma_k W_{k, Q_{k,t}} \leq 	\sup_{0 \leq u \leq Q_{k,t_{\delta}} } \sigma_k W_{k, u} \leq  	\sup_{0 \leq u \leq  t_\delta } \sigma_k W_{k, u} \leq  {t_\delta}^{2/3} . 
\end{equation}
Set $ c = \frac{1-\tilde p}{4 LK} $. The above inequality implies that
\begin{equation}\label{key}
	1 - \psi_1(Z_t) \geq \frac{1-\tilde p}{2}, \quad \forall t \in (0, t_\delta). 
\end{equation}
Combining this with \eqref{eq:R_to_Rtdelta}, we have that for all large $ \delta $: 
\begin{equation}\label{key}
	R \geq  (t-Q_{1,t}) \delta \geq   \delta \int_{0}^{t_{\delta}} 1 - \psi_1(Z_t) ds = \delta \frac{(1-\tilde p )/4LK}{\delta }\cdot \frac{1-\tilde p}{2} = \frac{(1-\tilde p)^2}{8KL}. 
\end{equation}
This proves our claim. 
\qed

\subsection{Proof of Corollary \ref{coro:smoothGaussBad}}
\label{app:coro:smoothGaussBad}

\bpf  
  For simplicity of notation, we will let $\sigma_k=1$ for all $k$, while the proof can be extended easily to the case where $\sigma_k >0$ for all $k$. 
For $ (q,s)  \in [0,1]^K \times \R^K$, define 
\begin{align}
	\bar \mu_k (q,s) = \frac{ s_k}{ q_k+ c}, \quad \bar \sigma_k(q,s) = \sqrt{\frac{1}{q_k+c}}, \quad k =1, \ldots, K. 
\end{align}
%
Then,  the limiting sampling function associated with Gaussian Thompson sampling, given in  \eqref{eq:TS_scale}, can be written as: 
\begin{equation}\label{eq:GTSsample}
	\psi_k(q,s) = \Pi^{\sup}(q, \, s; \, \sigma^2, \, c) = \pb\p{N( \bar \mu_k(q, s) , \bar \sigma_k(q, s)^2 ) > \max_{l \neq k} N(\bar \mu_l(q, s) , \bar \sigma_l(q, s)^2 )  },
\end{equation}
where $N(\mu, \sigma^2)$ represents a Gaussian random variable with mean $\mu$ and variance $ \sigma^2 $, independent from the rest of the system.  

Next, we show that, for all $ c>0 $, the  sampling function in \eqref{eq:GTSsample} is Lipschitz-continuous.  
Denote by $ f^{(q,s)} $ the density of the Gaussian distribution $\calN(\bar\mu(q,s), I \bar\sigma(q,s)^2)$, and by $ f^{(q,s)}_k $ the density of its $k$th marginal. 
Then, 
\begin{align}
		\psi_k(q,s) = &  \int_{\R^K} f^{(q,s)}(x) \mathbb{I}(\max_{l \neq k} x_l < x_k)   d x  \nln
		= &  \int_{\R^K} \p{ \prod_{l =1, \ldots, K} f^{(q,s)}_l (x_l)} \mathbb{I}(\max_{l \neq k} x_l < x_k)   d x . 
		\label{eq:psik_gaussTS_integral}
\end{align}
Because $ c>0 $, we have that $  \bar\sigma_k(q,s)$ is uniformly bounded from below by $ \sqrt{1/c} $  for all $q \in [0,1] $. This implies $ \bar \mu $ is Lipschitz, and that there exists constant $ D>0 $, such that 
\begin{equation}\label{eq:partial_marginal}
\abs{\frac{\partial}{\partial q_l} f^{(q,s)}_l(x_l)}  \leq D, \ \ \abs{\frac{\partial}{\partial s_l} f^{(q,s)}_l(x_l)}  \leq D, \quad \mbox{$ \forall l =1, \ldots, K$,  $ x, s \in \R^K $, $ q \in [0,1]^K $.} 
\end{equation}
Note that by the product form of $f^{(q,s)}$ and the fact that $ f_l $ is solely determined by $ (q_l, s_l) $, we have that
\begin{equation}\label{key}
	\frac{\partial}{\partial q_l} f^{(q,s)}(x)  = 	\frac{\partial}{\partial q_l}   \prod_{l =1, \ldots, K} f^{(q,s)}_l (x_l) \leq D { \prod_{l' \neq l } f^{(q,s)}_{l'} (x_{l'})} \leq \frac{D}{(2\pi/c)^{(K-1)/2}}, 
\end{equation}
where step $ (a) $ follows from \eqref{eq:partial_marginal} and $ (b) $ from the fact that $  \bar\sigma_k(q,s) \geq 1/\sqrt{c}$ for all $ k $ and hence $ f^{(q,s)}_{l'} (x_{l'}) \leq \frac{1}{\sqrt{2\pi /c}} $. The same bound holds for $ 	\frac{\partial}{\partial s_l} f^{(q,s)}(x)  $ via an identical argument. Combining the above inequality with \eqref{eq:psik_gaussTS_integral} and applying the Leibniz integral rule, we conclude that the magnitude of every coordinate of the gradient $ \nabla 	\psi_k(q,s)  $ is bounded from above by $   \frac{D}{(2\pi/c)^{(K-1)/2}}$. This proves that $ \psi(\cdot, \cdot) $ is Lipschitz-continuous. The fact that $ 		\liminf_{\delta \to \infty } R  >0 $ thus follows immediately from Theorem \ref{theo:lip_bad}. 

\qed

\subsection{Proof of Theorem \ref{theo:undersmooth_regret}}
\label{app:theo:undersmooth_regret}

We will prove Theorem \ref{theo:undersmooth_regret} by showing the following stronger version, which  includes a similar super diffusive regret bound for undersmoothed one-armed Thompson sampling introduced in \eqref{eq:TS_onearm}, where  only one of the two arms is uncertain. We do so because the analysis for the one-armed version of the algorithm is simpler, and will provide the crucial ingredients for understanding  undersmoothed two-armed Thompson sampling.

\begin{theo}
	\label{theo:undersmooth_regret_v2}
	Consider the diffusion limit associated with one- and two-armed undersmoothed Thompson sampling, with $ c=0 $. Then, the following holds almost surely for the diffusion regret \eqref{eq:scaled_regret} in the super diffusive regime.\footnote{When taking the  super diffusive scaling,  for the two-armed setting, it suffices to let the arm gap $\delta = \mu_1 - \mu_2 $ tend to $+\infty$; this case where $ \delta \to -\infty $ is covered here, due to both arms' mean rewards being uncertain. For the  one-armed setting, however, the symmetry  breaks down because there there is only one arm with uncertain mean rewards. For this setting, our analysis will cover both limits as the unknown mean $ \mu $ tends to $ +\infty $ and $ - \infty $, respectively. }
	\begin{enumerate}
		\item  For one-armed Thompson sampling, 
		\begin{equation}
			R \prec 1/|\mu| , \quad \mbox{as $|\mu|  \to \infty$}, 
			\label{eq:c0Reg}
		\end{equation}
		where $ \mu \in \R$  is the mean of the unknown arm. 
		\item  For two-armed undersmoothed Thompson sampling, 
		\begin{align}
			R \prec 1/\delta , \quad \mbox{as  $\delta \to \infty$}, 
		\end{align}
		where $ \delta = \mu_1 -\mu_2 $ is the arm gap. 
	\end{enumerate}
\end{theo}

We divide the proof of the above theorem into two sections, for the one- and two-armed Thompson sampling, respectively.  We begin by presenting the equivalent ODE characterization of the diffusion limit for one- and two-armed Thompson sampling using the random-time change in Theorem \ref{theo:timechange}.  These ODEs will be used repeatedly in our proof.  Fix $ c \geq 0 $.  

For one-armed Thompson sampling \eqref{eq:TS_onearm}, it is not difficult to verify that  the sample paths of $S^n$ and $Q^n$ associated with the uncertain arm	converge weakly to the SDE: 
	\begin{equation}
			\label{eq:TS1}
			dQ_t = \pi_t dt, \ \ \ \ \ \ dS_t = \mu \pi_t dt + \sqrt{\pi_t} dB_t, \ \ \ \ \ \ \pi_t = \Phi\p{\frac{S_t}{\sigma \sqrt{Q_t + \sigma^2c}}}. 
		\end{equation}
which can be expressed in the ODE form as follows: 
\begin{align}
	dQ_{t} =&\Phi\p{\frac{Q_t\mu + W_{Q_t}}{\sigma \sqrt{Q_t + \sigma^2c}}}   \, dt, \ \ \ S_t = \mu Q_t + \sigma W_{Q_t}, 
	\label{eq:TS1_timechange}
\end{align}
with $ Q_0=S_0=0 $ and $W$ is a standard one-dimensional Brownian motion. Likewise, for two-armed Thompson sampling, the following  is the ODE corresponding to the SDE in \eqref{eq:TS2}: 	
\begin{align}
	\begin{split}
		d {Q_{1,t}} = & \Phi\p{\frac{     \sigma^{-2} \p{ Q_{1,t}Q_{2,t}\delta + Q_{2,t}W_{1, Q_{1,t}} - Q_{1,t} W_{2, Q_{2,t}} } }{  \sqrt{ \sigma^{-2}   t Q_{1,t} Q_{2,t}  + t^2  c}} } dt,  \\
		d {Q_{2,t}} = &\Phi\p{\frac{     - \sigma^{-2}  \p{ Q_{1,t}Q_{2,t}\delta - Q_{2,t}W_{2, Q_{1,t}} + Q_{1,t} W_{1, Q_{2,t}} } }{  \sqrt{  \sigma^{-2}  t Q_{1,t} Q_{2,t}  + t^2  c}}  } dt, \\
		S_{k,t} = & \mu_k Q_{k,t} + \sigma_k W_{k, Q_{k,t} } ,  \quad k = 1, 2, 
	\end{split}
	\label{eq:TS2_timechange}
\end{align}
with $ Q_{\cdot, 0}= S_{\cdot, 0}=0 $, where $\delta = \mu_1 - \mu_2$ and $W$ is a standard two-dimensional Brownian motion. For our subsequent analysis, we will take $\sigma=1 $.

\subsubsection{One-Armed Setting}
\label{sec:proof:undersmooth_one}

Define the limit cumulative regret $R_t$ as:
\begin{equation}
 R_t = (\mu)_+t - \mu Q_t, \quad t\in [0,1], 
\end{equation} 
where $Q_t$ is the diffusion limit associated with $Q^n_i$. Note that $R_1$ corresponds to the scaled cumulative regret $R$ in \eqref{eq:scaled_regret}.

We will be working with the random-time-change version of the diffusion, as per \eqref{eq:TS1_timechange}. Define the  following stopping time: 
\begin{equation}\label{eq:tau_stopping}
	\tau(q) = \inf\{t: Q_t \geq  q \},
\end{equation}
with $\tau(q):=1$ if $Q_1 \leq  q$. As mentioned earlier, the sample paths exhibit different dynamics depending on whether $ \mu $ tends to the positive or negative infinity; as such, our proof will also treat these two cases separately.

\textit{Case 1: $\mu \to -\infty$.}  Fix $\alpha \in (1,2)$. Recall the stopping time from \eqref{eq:tau_stopping} with threshold $ q $. We will be particularly interested in the case where
\[q = |\mu|^{-\alpha}. \]
Decompose $R_1$ into two components: 
\begin{equation}
R_1 = R_{\tau(|\mu|^{-\alpha})} + (R_1 - R_{\tau(|\mu|^{-\alpha})}). 
\end{equation}
We next bound the two terms on the right-hand side of the equation above separately. Since $Q_t$ is non-decreasing, we have
\begin{equation}
R_{\tau(|\mu|^{-\alpha})} = |\mu| \cdot Q_{\tau(|\mu|^{-\alpha})} = |\mu|^{-(\alpha-1)}. 
\label{eq:Reg_par1}
\end{equation}
For the second term, the intuition is that by the time $Q_t$ reaches $|\mu|^{-\alpha}$, the drift in $Q_t$ will have already become overwhelmingly small for the rest of the time horizon. To make this rigorous, note the following facts: 
\begin{enumerate}
\item By the law of iterated logarithm of Brownian motion, along almost all sample paths of $ W $, there exists constant $C$ such that 
\begin{equation}
\limsup_{\mu \to - \infty} \sup_{x \in [\tau(|\mu|^{-\alpha}), 1]} \left|\frac{W_x}{\sqrt{x}} \right| \leq C\sqrt{\log \log (|\mu|^\alpha)}. 
\end{equation}
\item $\mu \sqrt{Q_{\tau(|\mu|^{-\alpha})}} = - |\mu|^{1-\alpha/2}$, and therefore $\frac{\mu \sqrt{Q_{\tau(|\mu|^{-\alpha})}}} {2}  \ll -  \sqrt{\log \log (|\mu|^{\alpha})}$ as $\mu \to -\infty$.
\end{enumerate}
Combining these facts along with the normal cdf tail bounds from Lemma \ref{lemm:gaussian_cdf} in Section \ref{app:tech}, we have that along almost all sample paths of $ W $, there exists constant $b>0$, such that for all sufficiently small $\mu$, 
\begin{align}
R_1 - R_{\tau(|\mu|^{-\alpha})} = & |\mu|\int_{\tau(|\mu|^{-\alpha})}^1 \Pi(0, Q_t) dt \nln
\leq & |\mu| \p{ \sup_{t \in [\tau(|\mu|^{-\alpha}), 1]} \Phi\p{\frac{\mu \sqrt{ Q_t}}{\sigma} +  \frac{W_{ Q_t}}{\sigma \sqrt{ Q_t}}} } \nln
\sk{a}{\leq} & |\mu| \Phi\p{ - |\mu|^{1-\alpha/2} +  C\sqrt{\log \log (|\mu|^\alpha)} } \nln
\leq &  |\mu| \exp(-|\mu|^b),
\label{eq:Reg_par2}
 \end{align}
where $ (a) $ follows from the aforementioned facts. Putting together \eqref{eq:Reg_par1} and \eqref{eq:Reg_par2} shows that
\begin{equation}
R_1 \leq \mu^{-(\alpha-1)} +  |\mu| \exp(- |\mu| ^b)  \stackrel{\mu\to \infty}{\prec} \mu^{-(\alpha-1)} , \quad a.s.
\end{equation}
This proves the claim by noting that the above holds for all $ \alpha \in (1,2) $.

\textit{Case 2: $  \mu \to \infty$.}  In this case, we would like to argue that $Q_t$ will increase rapidly as $\mu$ grows.  Let $\eta$ be a function of the form: 
\begin{equation}
\eta(x)= 1- x^{-\alpha}, 
\label{eq:nu_def}
\end{equation}
where $\alpha \in (1,2)$ is a constant; the value of $\alpha$ will be specified in a later part of the proof.   

The remainder of the proof will be centered around the dynamics of $ Q $ before and after the following stopping time: 
\begin{equation}
\tau(\eta(\mu)) = \inf \{t: Q_t \geq 1 - \mu^{-\alpha}   \}.
\end{equation}
It follows from the definition that if 
\begin{equation}
\tau(\eta(\mu))  <  1,
\label{eq:tau_leq_1}
\end{equation}
then
\begin{equation}\label{key}
	Q_1 \geq \eta (\mu)   = 1- \mu^{-\alpha}. 
\end{equation}
Therefore, if we can show that almost surely \eqref{eq:tau_leq_1} holds for all sufficiently large $\mu$, then it follows that  for all large $\mu$, the desired inequality holds: 
\begin{equation}
R_1 = \mu(1-Q_1) \leq  \mu^{- (\alpha-1)}, \quad a.s.
\label{eq:R1_case4}
\end{equation}

The remainder of the proof is devoted to showing \eqref{eq:tau_leq_1}. A main challenge in this part of the proof is that the dynamics of $ W_{Q_t}/\sqrt{Q_t} $, and by consequence that of $ Q_t $, is highly volatile near $ t=0 $. To obtain a handle on the behavior of these quantities, we will use the following trick by performing a change of variable in time so that the integrand below changes from $ s $ to $ 1/s $. We have that 
\begin{align}
\tau(\eta(\mu))=&  \int_{0}^{\eta(\mu)} 1/ \Phi\p{\frac{s\mu + W_{s}}{\sigma \sqrt{s }}} ds  \nln
= & \int_{\eta(\mu)^{-1}}^\infty u^{-2} \p{{\Phi\p{ \frac{\mu}{\sigma\sqrt{u}} +  \frac{\tilde W_u}{\sigma \sqrt{u}}}} }^{-1}du  \nln
=& \int_{\eta(\mu)^{-1}}^\infty u^{-2} \xi(\mu, u) du 
\label{eq:chage_var_int}
\end{align}
where $ u=1/s $, 
\begin{equation}
 \xi(\mu, u) :=  \p{{\Phi\p{ \frac{\mu}{\sigma\sqrt{u}} +  \frac{\tilde W_u}{\sigma \sqrt{u}}}} }^{-1}.
 \end{equation} 
and 
\begin{equation}\label{key}
\tilde W_t = tW_{1/t}. 
\end{equation}
Importantly, it is well known that if $W_t$ is a standard Brownian motion, then so is $\tilde W_t$. 
 
 We now bound the above integral using a truncation argument. For $K>\eta(\mu)^{-1}$, we write
\begin{align}
\tau(\eta(\mu)) = & \int_{\eta(\mu)^{-1}}^\infty u^{-2} \xi(\mu, u) du   \nln
= & \int_{\eta(\mu)^{-1}}^K u^{-2} \xi(\mu, u) du + \int_{K}^\infty u^{-2} \xi(\mu, u) du  \nln
\leq & \p{\sup_{u \in [\eta(\mu)^{-1}, K] }\xi(\mu, u) } \int_{\eta(\mu)^{-1}}^\infty u^{-2} du + \int_{K}^\infty u^{-2} \xi(\mu, u) du \nln
= & \p{\sup_{u \in [\eta(\mu)^{-1}, K]} \xi(\mu, u)} \eta(\mu)   + \int_{K}^\infty u^{-2} \xi(\mu, u) du  .
\label{eq:tau_int}
\end{align}
The following lemma bounds the second term in the above equation; the proof is given in Section \ref{app:lem:trunc_int}. 
\begin{lemm} 
\label{lemm:trunc_int}  For any $\delta \in (0,1)$, there exists $C>0$  such that,  along almost all sample paths of $ W $,  for all  large $\mu$ and $K$: 
\begin{equation}
\int_{K}^\infty u^{-2} \xi(\mu, u) du \leq CK^{-(1-\delta)}.
\end{equation}
\end{lemm}

Bounding the first term in \eqref{eq:tau_int} is more delicate, and will involve taking $\mu$ to infinity in a manner that depends on $K$. Fix any $\gamma \in (0,1)$ and consider a parameterization of $\mu$ where
\begin{equation}
\mu_K = K^{\frac{1}{2}+\gamma}, \quad K \in \N.
\end{equation}
By law of iterated logarithm (Lemma \ref{lemm:LIL} in Section \ref{app:tech}), and noting that $ \eta(\mu) <1 $, we have that there exists $C>0$ such that for all sufficiently large $K$
\begin{equation}
\inf_{u \in [\eta(\mu_K)^{-1}, K]}\frac{\tilde W_u}{\sqrt{u}} \geq -C \sqrt{\log \log K}, \quad a.s.
\end{equation}
Combining this with the lower bound on the normal cdf (Lemma \ref{lemm:gaussian_cdf}), we have, for all large $K$, 
\begin{align}
\sup_{u \in [\eta(\mu_K)^{-1}, K]} \xi(\mu_K, u) \leq  &  1+ \frac{\exp(-(\mu_K / \sqrt{K} - C\sqrt{\log\log K})^2)}{\mu_K / \sqrt{K} -C\sqrt{\log\log K}} \nln
= & 1+ \frac{\exp(-( K^{1/2+\gamma} / \sqrt{K} - C\sqrt{\log\log K})^2)}{K^{1/2+\gamma} / \sqrt{K} -C\sqrt{\log\log K}} \nln
\leq & 1+ \exp \p{-( K^\gamma - C\sqrt{\log\log K})^2 - \gamma \log K } \nln
 \leq&  1+ \exp( - K^{\gamma}).
 \label{eq:sup_xi} 
\end{align}

Fix $\nu \in (0,1)$, and $\delta, \gamma \in (0,1/4)$ such that 
\begin{equation}
2 > \frac{1-\delta}{1/2+\gamma}> 2-\nu.
\label{eq:gammadeltaratio}
\end{equation}
Note that such $\delta$ and $\gamma$ exist for any $\nu$, so long as we ensure that both $\delta$ and $\gamma$  are sufficiently close to $0$.  Combining \eqref{eq:tau_int}, \eqref{eq:sup_xi} and Lemma \ref{lemm:trunc_int}, we have that there exist $c_1, c_2>0$ such that,  along almost all sample paths of $ W $, for all large $K$: 
\begin{align}
\tau_{\mu_{K}} \leq & \p{\sup_{u \in [\eta(\mu_K)^{-1}, K]} \xi(\mu_K, u)}  \eta(\mu_K)   + \int_{K}^\infty u^{-2} \xi(\mu_K, u) du \nln
\leq & \p{1+ \exp( - K^{\gamma})}\eta(\mu_K) + c_1K^{-(1-\delta)} \nln
\leq &\eta(K^{  1/2+\gamma  }) + c_2K^{-(1-\delta)}, \quad a.s.
\end{align}
Recall that $ \eta(\mu) = 1- \mu^{-\alpha}$ .  We now choose $\alpha$ to be such that
\begin{equation}
2-\nu < \alpha < \frac{1-\delta}{(1/2 +\gamma) } <2.
\label{eq:alpha_small}
\end{equation}
Under this choice of $ \alpha $ (which exists because of \eqref{eq:gammadeltaratio}), we have that for all sufficiently large $K$
\begin{equation}
\tau_{\mu_{K}} \leq 1 - K^{-\alpha(1/2+\gamma)} + c_2K^{-(1-\delta)} < 1, \quad a.s.,
\end{equation}
where the last inequality follows from \eqref{eq:alpha_small}.   Combining the above equation,  \eqref{eq:tau_leq_1}, \eqref{eq:R1_case4} and the fact that $\nu$ can be arbitrarily close to 0, we have thus shown that for all $\alpha \in (1,2)$, 
\begin{equation}
R_1 \leq \mu (1- \eta (\mu) )= \mu^{-(\alpha-1)}, \quad a.s.,
\end{equation}
for all large $\mu$. This proves our main claim in this case, that is, almost surely 
\begin{equation}
R_1 \prec 1/\mu, \quad \mbox{as $\mu \to \infty$}.
\end{equation}

This completes the proof for the one-armed setting.

\subsubsection{Two-Armed Setting}
\label{sec:proof:undersmooth_two}

We now prove that the $ 1/\delta  $ regret scaling holds under undersmoothed two-armed Thompson sampling, as the arm gap $ \delta $ tends to infinity. Before delving into the details, let us first point out an intriguing and useful connection between the two-armed bandit analyzed here and the one-armed version presented in the previous subsection.  In the one-armed setting, a crucial simplifying feature   is that there is no uncertainty associated with the second, default arm, whereas in the two-armed case, both arms' mean rewards are uncertain. This manifests in there being only one  Brownian motion in the description of the diffusion limit for one-armed Thompson sampling, versus  two independent Brownian motions in the two-armed case. 
 
Fortunately, this distinction also suggests a plan of attack for analyzing the two-armed bandit. Consider the diffusion process at a small positive time $ t = v > 0 $. At this point, we have little knowledge of the behavior of $ Q_t $ because its dynamics near $ t=0 $ is highly volatile due to a lack of data during this period. However, we do know that the total arm pulls up to this point have to add up to $ v $, and so at least one of the two arms has been pulled by the amount $ v/2 $. The key insight here is that, depending on which arm has been pulled more by this point, we obtain, from time $ t = v $ onward, a version of diffusion that mirrors one of the two super diffusive regimes under the \textit{one-armed} Thompson sampling, i.e., $ \mu \to -\infty $ or $ \mu \to \infty $. Specifically:
\begin{enumerate}
	\item If arm 1 (superior arm) has been pulled by at least $ v/2 $, then we can approximately treat arm 1 as the ``certain'' arm, and our problem can be approximately mapped to a one-armed bandit with $ \mu \to -\infty $. 
	\item If, on the other hand, arm 2 (inferior arm) has been pulled by $ v/2 $, then arm 2 can be viewed as the certain arm, and the problem can be roughly reduced to a one-armed bandit with $ \mu \to +\infty$. 
\end{enumerate}
The above heuristic argument sets the stage for how the proof will proceed: we will consider two separate cases depending on the realization of $ Q_{1,t} $ at a carefully chosen, early point in time, and subsequently manipulate the drift equation to exploit the above-mentioned symmetry. That being said, the reduction from the two-armed bandit into two separate one-armed bandits is not exact, and the remainder of our proof is centered around using delicate estimates to make the above connection precise. 
 
We now present the formal proof for the regret scaling of undersmoothed two-armed Thompson sampling.  We will be working with the random-time-change version of the diffusion as per \eqref{eq:TS2_timechange}.  For clarity of notation, we will fix $ \sigma = 1 $ throughout the proof.  All results extend easily to the case of an arbitrary, fixed $ \sigma $.  Fix $\alpha \in (1, 2) $. Define
 \begin{equation}
 	v = 2\delta^{-\alpha}, 
 \end{equation}
 and the events
 \begin{equation}
 	\calE_k = \{Q_{k,v} \geq \delta^{-\alpha}\}, \quad k = 1, 2. 
 \end{equation}
 Because $t=  Q_{1,t}+Q_{2,t}$ for all $t$, we have that either $\calE_1$ or $\calE_2$ occurs almost surely. The proof will be completed by showing the claimed regret bound by conditioning upon each of these two events separately. Without loss of generality, the analysis of each case assumes the corresponding event occurs with strictly positive probability. Should one of the events occurs with probability zero, that portion of the proof can simply be ignored and should not impact the overall claim.

 \textit{Case 1.} First, suppose that $\calE_2$ has occurred, and without stating otherwise, this conditioning will be assumed throughout this portion of the proof. From \eqref{eq:TS2_timechange}, we have that 
 \begin{equation}
 	d {Q_{1,t}} = \Phi\p{     \frac{Q_{2,t}}{\sqrt{t Q_{2,t}}} \p{\frac{\delta Q_{1,t} + W_{1, Q_{1,t}}}{\sqrt{Q_{1,t}}}} - \sqrt{\frac{Q_{1,t}}{t}} \frac{W_{2,Q_{2,t}}}{\sqrt{Q_{2,t}}} } dt, 
 \end{equation}
 Define functions 
 \begin{align}
 	g_1(t) =& \frac{Q_{2,t}}{\sqrt{t Q_{2,t}}}, \quad g_2(t) = -\sqrt{\frac{Q_{1,t}}{t}} \frac{W_{2,Q_{2,t}}}{\sqrt{Q_{2,t}}} . 
 \end{align}
 Then, 
 \begin{equation}
 	d {Q_{1,t}} = \Phi\p{     g_1(t) \p{\frac{\delta Q_{1,t} + W_{1, Q_{1,t}}}{\sqrt{Q_{1,t}}}} + g_2(t) } dt.
 	\label{eq:g1g2}
 \end{equation}
 
 The following lemma shows that the terms $g_1$ and $g_2$ are appropriately bounded, and will be crucial to our subsequent analysis; the proof is presented Appendix  \ref{proof:lemm:g1g2prop}. 
 \begin{lemm}
 	\label{lemm:g1g2prop}
 	The following is true almost surely: 
 	\begin{enumerate}
 		\item $1 \geq g_1(t) \geq \delta^{-\alpha/2} $ for all $t \geq v$.
 		\item There exists constant $C$,  such that for all large $\delta$ and $t\geq v$,
 		\begin{equation}
 			|g_2(t) | \leq C \sqrt{\log \log \delta}.
 		\end{equation}
 	\end{enumerate}
 \end{lemm}

 Recall the stopping time defined in \eqref{eq:tau_stopping}: 
 \begin{equation}
 	\tau(q) = \inf \{t: Q_{1,t} \geq q\}. 
 \end{equation}
 Using the fact that $Q_{1,\tau(t)} = t$, we have that
 \begin{align}
 	\tau(q) = & \int_0^q \tau'(s) ds  \nln
 	=& \int_0^q 1/\Phi \p{ g_1(\tau(s)) \p{\frac{\delta s + W_{1, s}}{\sqrt{s}}} + g_2(\tau(s))} ds
 \end{align}
 For $q > Q_{1,v}$, we have, from the change of variables $s=  u^{-1}$ (see e.g., the discussion preceding \eqref{eq:chage_var_int} where this transformation was first used), 
 \begin{align}
 	\tau(q) =& v + \int_{Q_{1,v}}^q 1/\Phi \p{ g_1(\tau(s)) \p{\frac{\delta s + W_{1, s}}{\sqrt{s}}} + g_2(\tau(s))} ds \nln
 	= & v + \int_{1/q}^{1/Q_{1,v}} u^{-2}/\Phi \p{ g_1(\tau(u^{-1})) \p{ \frac{\delta}{\sqrt{u}} + \frac{\tilde W_u}{\sqrt{u}}} + g_2(\tau(u^{-1}))} du. 
 	\label{eq:tauq_twoarm_0}
 \end{align}
 Note that  $  	\tau(s) \geq v$ for all $ s \geq Q_{1,v}. $ 
 
  We now employ a truncation argument similar to that in the proof of the one-armed case, Case 4 \eqref{eq:tau_int}. Define 
 \begin{equation}
 	\xi(\delta, u) = 1/\Phi \p{ g_1(\tau(u^{-1})) \p{ \frac{\delta}{\sqrt{u}} + \frac{\tilde W_u}{\sqrt{u}}} + g_2(\tau(u^{-1}))}.
 \end{equation}
 For $M\in \R$, $1/q < M < 1/Q_{1,v}$, we have, from \eqref{eq:tauq_twoarm_0},  that
 \begin{align}
 	\tau(q) \leq & v + \int_{1/q}^M u^{-2}\xi(\delta, u) du + \int_M^{1/Q_{1,v}} u^{-2} \xi(\delta, u) du \nln
 	\leq &   v +  \p{\sup_{u \in [1/q, M]}\xi(\delta, u)}q + \int_M^{1/Q_{1,v}} u^{-2} \xi(\delta, u) du.
 	\label{eq:tauq_twoarm}
 \end{align}
 The last term in the above equation can bounded by the following lemma. The result is analogous to Lemma \ref{lemm:trunc_int}, with the key difference being that now the instance-specific parameter $\delta$ also features in the bound; the proof is given in Section \ref{proof:lemm:trunc_int_twoarm}. 
 \begin{lemm} 
 	\label{lemm:trunc_int_twoarm}  For any $\beta \in (0,1)$, there exist constants $B, C>0$  such that for all  large $\delta$ and $M$: 
 	\begin{equation}
 		\int_{M}^{1/Q_{1,v}} u^{-2} \xi(\delta, u) du \leq C(\log \delta)^B M^{-(1-\beta)}, \quad a.s. 
 	\end{equation}
 \end{lemm}

 Fix $\epsilon	 \in (1,2)$ and define 
 \begin{equation}
 	q^*(\delta) = 1- \delta^{-\epsilon}. 
 \end{equation} 
 Applying Lemma \ref{lemm:trunc_int_twoarm} to \eqref{eq:tauq_twoarm}, we have that for all $\beta \in (0,1)$ and sufficiently large $M$ and $\delta$: 
 \begin{align}
 	\tau(q^*(\delta)) \leq &   v +  \p{\sup_{u \in [q^*(\delta)^{-1}, M]}\xi(\delta, u)}q^*(\delta) + C(\log \delta)^B M^{1-\beta}.
 	\label{eq:tauq_twoarm2}
 \end{align}
 To bound the  term in the middle, we again resort to a double limit, in which $M$ and $\delta$ will tend to infinity simultaneously.  Fix $\gamma > \frac{\alpha/2}{ 2- \alpha}$,  and consider a parameterization of $\delta$: 
 \begin{equation}
 	\delta_M = M^{1/2+\gamma}.
 \end{equation}
 By the law of iterated logarithm (Lemma \ref{lemm:LIL} in Section \ref{app:tech}), we have that there exists $C>0$ such that for all large $\delta$ and $M$
 \begin{equation}
 	\inf_{u \in [ q^*(\delta_M)^{-1}, M]} \frac{\tilde W_u }{u} \geq - C\sqrt{\log \log M}. 
 	\label{eq:tildeWLILtwoarm}
 \end{equation}
 We have that for all large $M$
 \begin{align}
 	& \sup_{u \in [q^*(\delta_M)^{-1}, M]} \xi(\delta_M, u)  \nln
 	= & \sup_{u \in [q^*(\delta_M)^{-1}, M]}  1/\Phi \p{ g_1(\tau(u^{-1})) \p{ \frac{\delta}{\sqrt{u}} + \frac{\tilde W_u}{\sqrt{u}}} + g_2(\tau(u^{-1}))} \nln
 	\sk{a}{\leq}  & \sup_{u \in [q^*(\delta_M)^{-1}, M]}  1/\Phi \p{ g_1(\tau(u^{-1})) \p{ \frac{\delta}{\sqrt{u}} + \frac{\tilde W_u}{\sqrt{u}}} - c_1\sqrt{\log \log  \delta_M}} \nln
 	\sk{b}{\leq} & \sup_{u \in [q^*(\delta_M)^{-1}, M]}  1/\Phi \p{(\delta_M)^{1- \alpha/2} / \sqrt{u} -  c_2\p{\sqrt{\log \log u} + \sqrt{\log \log  \delta_M}  }} \nln
 	\leq &   1/\Phi \p{(\delta_M)^{1- \alpha/2} / M^{1/2}-  c_2\p{\sqrt{\log \log M} + \sqrt{\log \log \delta_M}  }}  \nln
 	= &   1/\Phi \p{M^{\gamma(1-\alpha/2) - \alpha/4}-  c_2\p{\sqrt{\log \log M} + \sqrt{\log \log \delta_M}  } }  \nln
 	\sk{c}{\leq} & 1+ \frac{\exp \p{-\p{M^{\gamma(1-\alpha/2) - \alpha/4}-  c_2\p{\sqrt{\log \log M} + \sqrt{\log \log \delta_M}  } }^2 } }{M^{\gamma(1-\alpha/2) - \alpha/4}-  c_2\p{\sqrt{\log \log M} + \sqrt{\log \log \delta_M}  } } \nln 
 	\leq & 1+ \exp(- M^{\gamma(1-\alpha/2) - \alpha/4}), 
 	\label{eq:sup_xi_two_arm}
 \end{align}
 where we note that the exponent $(\gamma(1-\alpha/2) - \alpha/4)$ is strictly positive based on the definition of $\gamma$. The steps are based on: 
 \begin{enumerate}[(a): ]
 	\item $g_2(t) \leq C\sqrt{\log \log \delta}$ when $t\geq v$.
 	\item $1 \geq g_1(t) \geq \delta^{-\alpha/2}$ for $t\in [v,1)$, and \eqref{eq:tildeWLILtwoarm}. 
 	\item The lower bound on the normal cdf for $x>0$ in Lemma \ref{lemm:gaussian_cdf}, Section \ref{app:tech}.
 \end{enumerate}
 
 We now substitute \eqref{eq:sup_xi_two_arm} into \eqref{eq:tauq_twoarm2}, and recall that $v = 2\delta^{-\alpha}$ and $q^*(\delta) = 1-\delta^{-\epsilon}$. We have 
 \begin{align}
 	\tau(q^*(\delta_M)) \leq &   v +  \p{\sup_{u \in [q^*(\delta_M)^{-1}, M]}\xi(\delta_M, u)}q^*(\delta_M) + C(\log \delta_M)^B M^{1-\beta} \nln
 	\leq & 2\delta_M^{-\alpha} +  \p{1+ \exp(- M^{\gamma(1-\alpha/2) - \alpha/4})}(1- \delta_M^{-\epsilon})+ C(\log \delta)^B M^{1-\beta}\nln
 	\leq & 1- \delta_M^{-\epsilon} + 2\delta_M^{-\alpha} + C(\log \delta_M)^B M^{1-\beta} +\exp(- M^{\gamma(1-\alpha/2) - \alpha/4}).
 	\label{eq:tauq_twoarm3}
 \end{align}
Let us choose the parameters $\epsilon $  and $\beta$   as follows:
 \begin{enumerate}
 	\item Let $\epsilon \in (1,\alpha)$, so that $\delta^{-\epsilon} \gg \delta^{-\alpha} $ as $\delta \to \infty$. 
 	\item Choose $\beta$ to be sufficiently close to $1$ such that $C(\log \delta_M)^B M^{1-\beta} < \delta_M ^{-\alpha}$ for all large $M$. 
 \end{enumerate}
 Under these choices of parameters, we see that $\delta_M^{-\epsilon} $ is orders-of-magnitude larger than the sum of the last three terms in \eqref{eq:tauq_twoarm3}, and we have that for all large $M$: 
 \begin{align}
 	\tau(q^*(\delta_M))
 	\leq & 1- \delta_M^{-\epsilon} /2 <1.
 	\label{eq:tauq_twoarm4}
 \end{align}
 
 We can now turn to the final objective, $R_1$. \eqref{eq:tauq_twoarm4} shows that $Q_{1,1} \geq 1- \delta_M^{-\epsilon}$. We have that for all large $M$
 \begin{equation}
 	R = \delta_M(1-Q_{1,1}) \leq \delta_M(1-q^*(\delta_M)) = \delta_M^{- (\epsilon-1)}. 
 \end{equation}
 Finally, notice that since $\alpha$ may take any value in $(1,2)$, so can $\epsilon-1$ take on any value in $(0,1)$. Considering that $\delta_M$ is a continuous increasing function of $M$ that tends to infinity as $M\to \infty$, we conclude that $R \prec 1/\delta$ as $\delta\to \infty$, as claimed.

 \textit{Case 2.} Next, we will look at the case where $\calE_1$ occurred. The structure of this part of the proof loosely mirrors Case 1 in the proof for the one-armed case (Section \ref{sec:proof:undersmooth_one}), where arm 1's reward is more certain compared to that of arm 2.  With this in mind, we will instead look at the drift of $Q_{2,\cdot}$: 
 \begin{equation}
 	d {Q_{2,t}} = \Phi\p{     \frac{Q_{1,t}}{\sqrt{t Q_{1,t}}} \p{\frac{-\delta Q_{2,t} - W_{2, Q_{2,t}}}{\sqrt{Q_{2,t}}}} + \sqrt{\frac{Q_{2,t}}{t}} \frac{W_{1,Q_{1,t}}}{\sqrt{Q_{1,t}}} } dt. 
 \end{equation}
 
 We now redefine    $g_1$ and $g_2$ in a fashion that mirrors symmetrically the first case: 
 \begin{equation}
 	g_1 (t) =  \frac{Q_{1,t}}{\sqrt{t Q_{1,t}}}, \quad g_2(t) =  \sqrt{\frac{Q_{2,t}}{t}} \frac{W_{1,Q_{1,t}}}{\sqrt{Q_{1,t}}}.
 \end{equation}
 We obtain
 \begin{equation}
 	d {Q_{2,t}} = \Phi\p{     g_1(t)\p{\frac{-\delta Q_{2,t} - W_{2, Q_{2,t}}}{\sqrt{Q_{2,t}}}} + g_2(t)} dt. 
 	\label{eq:g1g2_case2}
 \end{equation}
 Note that due to the symmetry, the bounds on $g_1$ and $g_2$ as laid out in Lemma \ref{lemm:g1g2prop} continue to hold. 
 
 We now define $\tau$ in terms of $Q_{2, \cdot}$:
 \begin{equation}
 	\tau(q) = \inf\{t: Q_{2,t} = q\}, 
 \end{equation}
 and define
 \begin{equation}
 	\tau^*_\delta =  \tau(\delta^{-\alpha}). 
 \end{equation}
 Because $\calE_1$ is assumed to have occurred, we have that 
 \begin{equation}
 	\tau^*_\delta \geq v. 
 \end{equation}
 The following property on the drift of $Q_{2,\cdot}$ will be used in the remainder of the proof; the proof is given in Section \ref{proof:lemm:dQ2bound}. 
 
 \begin{lemm}  Along almost all sample paths of $ W $, there exists a constant $C$ such that for all large $\delta$,
 	\label{lemm:dQ2bound}
 	\begin{equation}
 		\frac{d}{dt} Q_{2,t}\leq \Phi\p{  -g_1(t)\sqrt{Q_{2,t}} \delta   + C\sqrt{\log \log \delta }}, \quad \forall t \in [\tau^*_\delta,1).  
 		\label{eq:g1g2_boundcase2}
 	\end{equation}
 \end{lemm}

 Decompose the regret as follows: 
 \begin{equation}
 	R = R_{\tau(\delta^{-\alpha}) } + (R-R_{\tau(\delta^{-\alpha}) }). 
 	\label{eq:reg_decomp_twoarm}
 \end{equation}
 For the first term, we have that 
 \begin{equation}
 	R_{\tau(\delta^{-\alpha}) } \leq \delta \cdot \delta^{-\alpha} = \delta^{-(\alpha-1)}.
 \end{equation}
 To bound the second term, we will again aim to show that once $Q_{2,\cdot}$ reaches $\delta^{-\alpha}$, its drift would become overwhelmingly small as $\delta$ gets large. Compared to the one-armed setting, a major obstacle in this case is that the $\delta^{-\alpha/2}$ uniform lower bound on $g_1(t)$ in Lemma \ref{lemm:g1g2prop}  turns out to be too weak for our purpose. We will rely on the following stronger lower bound on $g_1(t)$; the proof is given in Section \ref{proof:lemm:g1bound_twoarm}. 
 
 \begin{lemm} 
 	\label{lemm:g1bound_twoarm}
Along almost all sample paths of $ W $,  we have that for all sufficiently large $\delta$,  under the conditioning of event $\calE_1$, 
 	\begin{equation}
 		g_1(t) \geq \sqrt{1/3}, \quad \forall t \in [\tau^*_\delta,1). 
 	\end{equation}
 \end{lemm}

 With the strengthened lower bound on $g_1(t)$ at hand, we are now ready to bound the second term in \eqref{eq:reg_decomp_twoarm}. Combining Lemmas \ref{lemm:dQ2bound} and \ref{lemm:g1bound_twoarm}, we have that  for all large $\delta$
 \begin{align}
 	\sup_{t \in [\tau^*_\delta, 1)} \frac{d}{dt}Q_{2,t} \leq & \sup_{t \in [\tau^*_\delta, 1)}\Phi\p{    -\sqrt{1/3}\sqrt{Q_{2,t}} \delta  + C\sqrt{\log \log \delta}} \nln
 	\leq & \Phi\p{    -\sqrt{1/3} \delta^{1-\alpha/2}+ C\sqrt{\log \log \delta}} \nln
 	\prec & \exp(- \delta^{1-\alpha/2}), 
 \end{align}
 where the second inequality follows from the definition of $\tau^*_\delta$.  We can write the regret term as
 \begin{align}
 	R - R_{\tau^*_\delta} = \delta \int_{\tau^*_\delta}^1 \frac{d}{dt}Q_{2,t}  \prec \delta \exp(- \delta^{1-\alpha/2}) \ll \delta^{-1}.
 \end{align}
 This shows that the regret $R$ is dominated by $R_{\tau^*_\delta}$, so that 
 \begin{equation}
 	R \prec \delta^{-(\alpha-1)}, 
 \end{equation}
 and the claim follows from the fact that $\alpha$ can be arbitrarily close to $2$.

\subsection{Proof of Theorem \ref{theo:sample_path}}
\label{proof:theo:sample_path}
\bpf Using the random-time change characterization of the diffusion limit from Theorem \ref{theo:timechange}, we have that the diffusion limit under $ \psi $ satisfy
\begin{align}\label{key}
	dQ_{t} =& \psi( S_t, Q_t) \, dt, \quad k = 1, \ldots, K,  \nln
S_{t} =& \mu Q_t + \sigma W_{Q_t}\, dt, \quad k = 1, \ldots, K, 
\end{align}
where $ W $ is a standard Brownian motion. 
In particular, we have that 
\begin{equation}\label{key}
	S_t / \sqrt{Q_t} = \mu \sqrt{Q_t} + \sigma W_{Q_t}/\sqrt{Q_t}.
\end{equation}
By the law of iterated logarithm (Lemma \ref{lemm:LIL}, Appendix \ref{app:tech}) of Brownian motion, as well as the fact that $ \lim_{t\downarrow0} Q_t = 0 $, we have that almost surely 
\begin{align}\label{key}
	\limsup_{t \downarrow 0} \frac{W_{Q_t}/\sqrt{Q_t}}{\sqrt{2 \log \log (1/Q_t)}  } =  \liminf_{t \downarrow 0} \frac{W_{Q_t}/\sqrt{Q_t}}{-\sqrt{2\log \log (1/Q_t)}  } = 1.
\end{align}
 Since the term $\frac{\mu\sqrt{Q_t}}{\sigma}  $ is always bounded, we conclude that, almost surely, the value of $ 	S_t / \sqrt{Q_t}$ will oscillate between arbitrarily large positive and negative values as $ t\to 0 $. Applying the definition of a sensitive sampling function immediately leads to the claim. 
\qed

\section{Technical Lemmas}
\label{app:tech}

We will use the following technical lemmas  repeatedly. 
\begin{lemm}[Gaussian Tail Bounds]
	\label{lemm:gaussian_cdf}
	For all $x<-\sqrt{2\pi/(9-2\pi)}$: 
	\begin{align}
		\Phi(x) \leq & \frac{1}{|x|}\exp(-x^2/2), \quad \Phi(x) \geq  \frac{1}{3|x|}\exp(-x^2/2). 
	\end{align}
	This immediately implies that for all $x > \sqrt{2\pi/(9-2\pi)}$: 
	\begin{equation}
		\Phi(x) \leq 1- \frac{1}{3x}\exp(-x^2/2),  \quad \Phi(x) \geq 1- \frac{1}{x}\exp(-x^2/2). 
	\end{equation}
\end{lemm}

\bpf For the lower bound, we have that for all $x<0$: 
\begin{align}
	\Phi(x) = \frac{1}{\sqrt{2\pi}}\int_{-x}^\infty  \exp(-s^2/2) ds \sk{a}{\leq} \frac{1}{\sqrt{2\pi}}\int_{-x}^\infty \frac{-x}{s} \exp(-s^2/2) ds < \frac{1}{|x|}\exp(x^2/2), \nnb
\end{align}
where $(a)$ follows from the fact that $-x/s\leq 1$ for all $s \geq -x$. For the upper bound, define $f(x) = \frac{x}{\sqrt{2\pi} (x^2+1)}\exp(-x^2/2) - \Phi(-x) $. We have that $f(0) =- \Phi(0)<0$, $\lim_{x \to \infty} = 0$, and 
\begin{equation}
	f'(x) =  \frac{1}{\sqrt{2\pi}(1+x^2)^2}\exp(-x^2/2)>0, \quad \forall x >0. 
\end{equation}
This implies that $f(x) <0$ for all $x >0$, which further implies that 
\begin{equation}
	\Phi(-x) \geq \frac{|x|}{\sqrt{2\pi} (x^2+1)} \exp(-x^2/2), \quad \forall x <0. 
\end{equation}
The claim follows by noting that $ \frac{|x|}{\sqrt{2\pi} (x^2+1)} \geq \frac{1}{3|x|}$ whenever $|x|\geq \sqrt{2\pi/(9-2\pi)}$.
\qed

The next result is a well-known, fundamental property of the Brownian motion. Proofs can be found in standard texts (e.g., \cite[Theorem 9.23]{karatzas2005brownian}).
\begin{lemm}[Law of  Iterated Logarithm]
	\label{lemm:LIL}
	Let $ W_t $ be a standard Brownian motion. Then, almost surely,
	\begin{align}\label{key}
	\limsup_{t \to \infty} \frac{W_t}{\sqrt{2t \log \log t}  } = 	\limsup_{t \downarrow 0} \frac{W_t}{\sqrt{2t \log \log (1/t)}  } = 1. \\
	\liminf_{t \to \infty} \frac{W_t}{\sqrt{2t \log \log t}  }  = \liminf_{t \downarrow 0} \frac{W_t}{\sqrt{2t \log \log (1/t)}  } = -1.
	\end{align}
Furthermore, 
\begin{align}\label{key}
	\limsup_{t \downarrow 0} \frac{W_t}{\sqrt{2t \log \log (1/t)}  } =  & 	\limsup_{t \downarrow 0} \sup_{0 \leq u \leq t} \frac{W_u}{\sqrt{2u \log \log (1/u)}  } = 1 \\
	\liminf_{t \downarrow 0} \frac{W_t}{\sqrt{2t \log \log (1/t)}  }  =  & 	\liminf_{t \downarrow 0} \inf_{0 \leq u \leq t} \frac{W_u}{\sqrt{2u \log \log (1/u)}  } = -1
\end{align}
	
	\end{lemm}

\section{Additional Proofs}

\subsection{Proof of Lemma \ref{lem:no_trunc_int}}
\label{app:lem:no_trunc_int}
\bpf The proof is based on \citet{durrett1996stochastic}, and we include it here for completeness. For $\Delta^n_\epsilon$, note that for all $\epsilon>0$, 
\begin{equation}
\Delta^n_\epsilon(z)  \leq \frac{1}{\epsilon^p} m^n_p(z). 
\end{equation}
The convergence of \eqref{eq:mnp} thus implies that of \eqref{eq:deltan}. For $b^n_k$, note that 
\begin{equation}
|\tilde  b^n_k(z) - b^n_k(z)| = \int_{x: |x-z|>1} |x-z|K^n(z,dx) \leq m^n_p(z), 
\end{equation}
where the last step follows from the assumption $p\geq 2$. We have thus proven that \eqref{eq:mnp} and \eqref{eq:tild_bn} together imply \eqref{eq:bn}.   Finally, for $a^n_{k,l}$, we have 
\begin{align}
|\tilde a^n_{k,l}(z) - a^n_{k,l}(z)| \leq & \int_{x: |x-z|>1} |(x_k - z_k)(x_l - z_l)| K^n(z,dx) \nln
\leq & \int_{x: |x-z|>1}  |x-z|^2 K^n(z,dx) \leq m^n_p(z), 
\end{align}
whenever $p\geq 2$, where the first step follows from the Cauchy-Schwartz inequality, and the second step from the observation that $(x_k - z_k)^2 \leq |x-z|^2$ for all $k$. This shows \eqref{eq:mnp} and \eqref{eq:tild_an} together imply \eqref{eq:an}, completing our proof. 
\qed

\subsection{Proof of Lemma \ref{lemm:trunc_int}}
\label{app:lem:trunc_int}
By the law of iterated logarithm (Lemma \ref{lemm:LIL} in Section \ref{app:tech}), we have that there exists constant $C>0$ such that for all large $u$
\begin{equation}
\tilde W_u /\sqrt{u} \leq -C\sqrt{\log \log u}, \quad a.s.
\end{equation}
Using the lower bound on the normal cdf:
\begin{equation}
\Phi(x) \geq \frac{1}{\sqrt{2\pi}}\frac{-x}{1+x^2}\exp(-x^2/2), \quad x<0, \nln
\end{equation}
We have that for all large $u$, almost surely:
\begin{align}
\xi(\mu, u) =   1/{\Phi\p{ \frac{\mu}{\sigma\sqrt{u}} +  \frac{\tilde W_u}{\sigma \sqrt{u}}}} \leq   1/{\Phi\p{ \frac{\tilde W_u}{\sigma \sqrt{u}}}}
\leq b_1 (\log u)^{b_2}, 
\end{align}
where $b_1$ and $b_2$ are positive constants that do not depend on $u$, which further implies that for any $\alpha \in (0,1)$
\begin{equation}
\int_{K}^\infty u^{-2}\xi(\mu, u)   \leq 1/K^{1-\alpha}, 
\end{equation}
for all large $K$.
\qed

\subsection{Proof of Lemma \ref{lemm:g1g2prop}}
\label{proof:lemm:g1g2prop}
\bpf The first claim follows directly from the definitions of $g_1$ and $\calE_2$. In particular, for all $ t \geq \nu $, $ t\leq 1 $, and under $ \calE_2 $, 
\begin{equation}\label{key}
	g_1(t) \geq \frac{\sqrt{Q_{2,t}}}{  \sqrt{t} } \geq \sqrt{Q_{2,t}} \geq \delta^{-\alpha/2}. 
\end{equation}

The second claim follows from  the law of iterated logarithm (Lemma \ref{lemm:LIL}) applied to $W_{2, \cdot}$  and the fact that $Q_{1,t} \leq t$.   
\qed

\subsection{Proof of Lemma \ref{lemm:trunc_int_twoarm}}
\label{proof:lemm:trunc_int_twoarm}
\bpf  The proof of the result will use the law of iterated logarithm of Brownian motion, along with the bounds we have developed for the functions $g_1$ and $g_2$ in Lemma \ref{lemm:g1g2prop}. We have that there exist constants $b_1, \ldots, b_4$ and $c_1, \ldots,c_6$, such that for all sufficiently large $M$: 
\begin{align}
	& \int_{M}^{1/Q_{1,v}} u^{-2}  \xi(\delta, u) du  \nln
	=& \int_{M}^{1/Q_{1,v}} u^{-2}  /\Phi \p{ g_1(\tau(u^{-1})) \p{ \frac{\delta}{\sqrt{u}} + \frac{\tilde W_u}{\sqrt{u}}} + g_2(\tau(u^{-1}))}du \nln
	\sk{a}{\leq} & \int_{M}^{\infty} u^{-2}  /\Phi \p{ -  \abs{ \frac{\tilde W_u}{\sqrt{u}}}- c_1 \sqrt{\log \log \delta}}du \nln
	\sk{b}{\leq}& \int_{M}^{\infty} u^{-2}  /\Phi \p{-  c_3(\sqrt{\log \log u}+\sqrt{\log \log \delta} ) }du \nln
	\sk{c}{\leq}& c_4 \int_{M}^{\infty} u^{-2} \p{\sqrt{\log \log u}+\sqrt{\log \log \delta} }  \exp \p{ c_3^2 \p{\sqrt{\log \log u}+\sqrt{\log \log \delta} }^2 }du \nln
	\sk{d}{\leq}& c_4 \int_{M}^{\infty} u^{-2} \p{\sqrt{\log \log u}+\sqrt{\log \log \delta} } {(\log u)^{b_1}  (\log \delta)^{b_2}}  du \nln
	{\leq}& c_5( \log  \delta)^{b_3} \int_{M}^{\infty} u^{-2}  (\log u)^{b_4} du \nln
	{\leq}& c_6 ( \log  \delta)^{b_3} M^{1-\beta}.
\end{align}
where the various steps are based on the following facts: 
\begin{enumerate}[(a):]
	\item $g_1(t) \leq 1$, $|g_2(t)| \leq c_1\sqrt{\log \log \delta}$ for all $t \geq v$ and that $\tau(u^{-1}) \geq v$ when $u \leq 1/Q_{1,u}$.
	\item Law of iterated logarithm applied to $\tilde{W}$.
	\item Lower bound on the normal cdf (Lemma \ref{lemm:gaussian_cdf} in Section \ref{app:tech}): $\Phi(x) \geq \frac{1}{\sqrt{2\pi}}\frac{-x}{1+x^2}\exp(-x^2/2)$, for $ x<0$. 
	\item Cauchy-Schwartz. 
\end{enumerate}
This proves our claim. 
\qed

\subsection{Proof of Lemma \ref{lemm:dQ2bound}}
\label{proof:lemm:dQ2bound}
\bpf
From \eqref{eq:g1g2_case2}, we have that
\begin{align}
	\frac{d}{dt}  Q_{2,t} = & \Phi\p{     g_1(t)\p{\frac{-\delta Q_{2,t} - W_{2, Q_{2,t}}}{\sqrt{Q_{2,t}}}} + g_2(t)} \nln
	\leq & \Phi\p{    -g_1(t)\sqrt{Q_{2,t}} \delta  + \abs{\frac{W_{2, Q_{2,t}}}{\sqrt{Q_{2,t}}} } +  g_2(t)} \nln
	\leq & \Phi\p{    -g_1(t)\sqrt{Q_{2,t}} \delta  + C\sqrt{\log \log \delta}}, 
\end{align}
where the first inequality follows from the fact that $g_1(t) \leq 1$. The second follows from applying the law of iterated logarithm to $ \abs{\frac{W_{2, Q_{2,t}}}{\sqrt{Q_{2,t}}} } $, the fact that $Q_{2,t} \geq \delta^{-\alpha}$ by the definition of $\tau^*_\delta$, and the upper bound on $|g_2(t)|$ from Lemma \ref{lemm:g1g2prop}. 
\qed

\subsection{Proof of Lemma \ref{lemm:g1bound_twoarm}}
\label{proof:lemm:g1bound_twoarm}
\bpf 
Define the function 
\begin{equation}
	h(t) = \frac{Q_{2,t}}{t}. 
\end{equation}
We have
\begin{equation}
	g_1(t) = \sqrt{\frac{Q_{1,t}}{t}} = \sqrt{1-h(t)}.
\end{equation}
Our goal is to show that $h(t)$ stays below $2/3$ for all $t \in [\tau^*_\delta,1)$. First, we verify that 
\begin{equation}
	h(\tau^*_\delta) \leq 1/2. 
	\label{eq:htaudelta}
\end{equation}
This is true under the condition of $\calE_1$, by noting that $Q_{1,\tau^*_\delta}  \geq Q_{2, \tau^*_\delta}$ and $ Q_{1,t}+Q_{2,t} = t $, and hence 
\begin{equation}
	h(\tau^*_\delta) = \frac{Q_{2,\tau^*_\delta}}{Q_{1,\tau^*_\delta} + Q_{2,\tau^*_\delta}}  \leq 1/2. 
\end{equation}

Suppose, for the sake of contradiction that $ h(t_0) = 1/2 $ for some $ t_0\geq \tau^*_\delta $. 
In light of \eqref{eq:htaudelta}, it suffices to show that $h'(t_0 )<0$, which would imply that $h$ will not increase beyond $1/2$ to ever reach $2/3$. We have that
\begin{equation}
	h'(t) = \frac{\frac{d}{dt} Q_{2,t}}{t} - \frac{Q_{2,t}}{t} \frac{1}{t} = \frac{\frac{d}{dt} Q_{2,t} - h(t)}{t}.
\end{equation}
Evaluating this derivative at $t=t_0$ yields
\begin{equation}
	h'(t_0) = \frac{\frac{d}{dt}Q_{2,t} - h(t_0)}{t_0} = \frac{\frac{d}{dt}Q_{2,t_0} - 1/2}{t_0}.
\end{equation}
It thus suffices to show that $\frac{d}{dt}Q_{2,t_0} < 1/2$. To this end, note that by the definition of $t_0$, 
\begin{equation}
	h(t) \leq 1/2, \quad \forall t \in [\tau^*_\delta,t_0]. 
	\label{eq:hupper}
\end{equation}
By Lemma \ref{lemm:dQ2bound}, we have that for all large $\delta$ and $t \in [\tau^*_\delta,t_0)$, 
\begin{align}
	\frac{d}{dt}Q_{2,t} \leq & \Phi\p{  -g_1(t)\sqrt{Q_{2,t}} \delta   + C\sqrt{\log \log \delta }} \nln
	= & \Phi\p{  -\sqrt{1-h(t)}\sqrt{Q_{2,t}} \delta   + C\sqrt{\log \log \delta }} \nln
	\sk{a}{\leq}  & \Phi\p{  -\sqrt{\frac{Q_{2,t}}{2}} \delta   + C\sqrt{\log \log \delta }} \nln
	\sk{b}{\leq}  & \Phi\p{  - \delta^{1-\alpha/2}/\sqrt{2}  + C\sqrt{\log \log \delta }} \nln
	< & 1/2, 
\end{align}
where step $(a)$ follows from \eqref{eq:hupper} and $(b)$ from the fact that $Q_{2,t} \geq \delta^{-\alpha}$ for all $t \geq \tau^*_\delta$. The final inequality follows from the fact that $\alpha < 2$ and hence $\delta^{1-\alpha/2} \gg \sqrt{\log \log \delta } $ for all large $\delta$. This proves the claim.  \qed

\section{Uniqueness of Diffusion Limit for One-Armed Undersmoothed Thompson Sampling}
\label{sec:proof:onearmTS}

Theorem \ref{theo:conv_triangular_psi} provides a generic recipe to obtain approximation guarantees for a continuous limiting sampling function that is non-Lipschitz. The resulting diffusion limits, however, may not be unique. In this section, we prove a  stronger characterization for the special case of undersmoothed one-armed Thompson sampling, with $ c=0 $, which leads to  a unique solution to the associated SDE. 

We will construct the SDE solution at $ c=0 $ by considering the sequence of solutions under a diminishing, but strictly positive, sequence of $ c $. The following theorem shows that the sequence of diffusion processes almost surely admits a unique limit as $c \to 0$.  The key step of the proof hinges on establishing that, as $ c\to0 $, the drift term of the diffusion does not exhibit too wild of an oscillation near $t=0$. This would further allow us to use the dominated convergence theorem and show the soundness of the limit solution. 

\begin{theo}
\label{theo:onearmTS} 
The diffusion limit $(Q_t)_{t \in [0,1]}$ under one-armed Thompson sampling converges uniformly	 to a limit $\tilde{Q}$ as $c \to 0$ almost surely. Furthermore, $\tilde{Q}$ is a strong solution to the stochastic differential equation: 
\begin{equation}
	d \tilde Q_t =   \Phi\p{\frac{\mu \sqrt{\tilde Q_t}}{\sigma}+  \frac{W_{\tilde Q_t}}{\sigma \sqrt{\tilde Q_t}}} dt, \quad \tilde Q_0=0. 
	\label{eq:sde_1d_ts}
\end{equation}
\end{theo}

\bpf We will use the characterization of the diffusion process in \eqref{eq:TS1_timechange}. Although our main focus is on the process $Z_t$ restricted to the $[0,1]$ interval, the diffusion process itself is in fact well defined on $ t\in [0,\infty)$. A useful observation we will make here is that, for any $c > 0$ and $\mu$, we have that almost surely
\begin{equation}
Q_t \to \infty, \quad \mbox{as $t \to \infty$}. 
\end{equation}
This fact can be verfied by noting that if $Q_t$ were bounded over $[0,\infty)$, then its drift $\Pi(c, Q_t)$ would have been bounded from below by a strictly positive constant, leading to a contradiction. Recall the stopping time from \eqref{eq:tau_stopping}, and here we emphasize the dependence on $ c $: 
\begin{equation}
\tau^c(q) = \inf\{t: Q_t \geq q \}.
\end{equation}
Note that for all $c>0$, $Q_t$ and $\tau^c$ are increasing and continuous and $Q_t = \tau^{-1}(t)$.

As $ c\to 0 $, we face the same challenge of the volatile behavior of the drift function near $ t=0 $, as was discussed in the paragraph preceding \eqref{eq:chage_var_int}, and we will resort to the same change of variable technique as that in \eqref{eq:chage_var_int}.  From \eqref{eq:TS1_timechange},  we have that for $q > 0$: 
\begin{align}
\tau^c(q) =&  \int_{0}^q 1/ \Phi\p{\frac{s\mu + W_{s}}{\sigma \sqrt{s + \sigma^2c}}} ds  \nln
= &  \int_{1/q}^\infty  u^{-2} / \Phi\p{\frac{ \mu+ u  W_{1/u}}{\sigma \sqrt{u + u^2 \sigma^2 c} } } du \nln
= &  \int_{1/q}^\infty  u^{-2} / \Phi\p{\frac{ \mu+ \tilde W_u}{\sigma \sqrt{u + u^2 \sigma^2 c} } } du  \nln
=&  \int_{1/q}^\infty  h^c(u) du 
\label{eq:tauc_eq}
\end{align}
where $\tilde W_t = t W_{1/t}$ is a standard Brownian motion, and 
\begin{equation}
h^c(u) = u^{-2} / \Phi\p{\frac{ \mu+ \tilde W_u}{\sigma \sqrt{u + u^2 \sigma^2 c} } }. 
\end{equation}
By the law of iterated logarithm of Brownian motion (Lemma \ref{lemm:LIL} in Section \ref{app:tech}), we have that 
\begin{equation}
\limsup_{t \to \infty}\frac{|\tilde W_t|}{\sqrt{2 t \log\log t}} = 1, \quad a.s.
\end{equation}
Fix a sample path of $W$ such that the above is satisfied. Then, there exist  $M \in (1/q, \infty)$ and $D >0$, such that 
\begin{equation}
\tilde W_t  \geq - D \sqrt{t \log\log t} , \quad \forall t \geq M. 
\label{eq:tildeWt_M}
\end{equation}

We now consider two cases depending on the sign of $\mu$. First, suppose that $\mu \geq 0$.
Define 
\begin{equation}
g(u) =   \begin{cases}
	u^{-2} / \Phi\p{\frac{  \mu- |\tilde{W_u}|}{\sigma \sqrt{u}}}, & 1/q \leq u < M,\\
	u^{-2} / \Phi\p{ - \frac{D}{\sigma}\sqrt{ \log\log u} }, & u \geq M. 
\end{cases} 
\label{eq:gudef}
\end{equation}
where $M \in (1/q, \infty)$ is chosen such that \eqref{eq:tildeWt_M} holds. It follows from the definitions that 
\begin{equation}
h^c(u) \leq g(u), \quad \forall c\geq 0, u \geq 1/q. 
\label{eq:hc_dom}
\end{equation}

We now show that $g$ is integrable over $u \in [1/q, \infty)$. It suffices to show that 
\begin{equation}
\int_{M}^\infty g(u) du < \infty. 
\label{eq:gu_intgrable}
\end{equation} 
Recall the following lower bound on the cdf of standard normal from Lemma \ref{lemm:gaussian_cdf} in Section \ref{app:tech}: for all sufficiently small $ x $
\begin{equation}
\Phi(x) \geq \frac{1}{\sqrt{2\pi}}\frac{-x}{1+x^2}\exp(-x^2/2), \quad x<0. 
\label{eq:Phx_LB}
\end{equation}
We thus have that, for all sufficiently large $u$, 
\begin{align}
g(u) = & u^{-2}/\Phi\p{ - \frac{D}{\sigma}\sqrt{\log\log u} } \nln
\leq & u^{-2}\sqrt{2\pi} \frac{1+ \frac{D^2}{\sigma^2}\log \log u}{\frac{D}{\sigma} \sqrt{\log \log u}}  \exp\p{ \frac{D^2}{2 \sigma^2} \log \log u} \nln
\leq & 2 u^{-2} \sqrt{2\pi} \frac{D}{\sigma }\sqrt{\log\log u} \p{\log u}^{D^2/2\sigma^2} \nln
\leq & b_1 u^{-2} ( \log u) ^{b_2}, 
\label{eq:gu_def}
\end{align}
where $b_1$ and $b_2$ are positive constants.  Noting that 
\begin{equation}
(\log u)^\alpha \ll \sqrt{u}
\end{equation} 
as $u \to \infty$ for any constant $\alpha>0$, we have that 
$b_1  u^{-2} (\log u)^{b_2}$ is integrable over $(M, \infty)$ for all sufficiently large $ M $.  This proves the integrability of $g$ in \eqref{eq:gu_intgrable}. 

Using \eqref{eq:hc_dom}, \eqref{eq:gu_intgrable} and the dominated convergence theorem, we thus conclude that, for all $q >0$,  
\begin{align}
\lim_{c \downarrow 0} \tau^c(q) =  &  \int_{1/q}^\infty \lim_{c\to \infty}  \p{u^{-2} / \Phi\p{\frac{ \mu+ \tilde W_u}{\sigma \sqrt{u + u^2 \sigma^2 c} } } } du  \nln 
= & \int_{0}^q  1/\Phi\p{\frac{s\mu + W_{s}}{\sigma \sqrt{s}}} ds : =  \tau^0(q), \quad a.s. 
\end{align}
Recall that $Q_t = (\tau^c)^{-1}(t)$, the above thus implies that, for all $t \in [0,1]$, 
\begin{equation}
Q_t \stackrel{c \downarrow 0}{\longrightarrow} \tilde Q_t := (\tau^0)^{-1}(t), \quad a.s.
\end{equation}
Finally, the point-wise convergence implies uniform convergence over the compact interval $[0, 1]$ since $Q_t$ is $1$-Lipschitz. This proves our claim in the case where $\mu \geq 0$. The case for $\mu<0$ follows an essentially identical set of arguments after adjusting the constants $M$ and $D$ in \eqref{eq:gudef}, recognizing that the behavior of $\mu +\tilde W^u$ is largely dominated by that of $\tilde W^u$ when $u$ is large, which can be in turn bounded by the law of iterated logarithm.  This completes the proof of Theorem \ref{theo:onearmTS}.
\qed

\ifx \useplain\undefined
\end{APPENDICES}
\else
\end{appendix}
\fi

\end{document}

%% file: title_info/informs_title.tex
\TITLE{\rev{Weak Signal} Asymptotics for \\ \rev{Sequentially Randomized} Experiments} 
\RUNTITLE{Weak Signal Asymptotics for Sequentially Randomized Experiments}

\RUNAUTHOR{Kuang and Wager}

\ARTICLEAUTHORS{%

\AUTHOR{Xu Kuang}
\AFF{Graduate School of Business\\ Stanford University\\
 \EMAIL{kuangxu@stanford.edu}}  
 
 \AUTHOR{Stefan Wager}
 \AFF{Graduate School of Business\\ Stanford University\\
 	\EMAIL{swager@stanford.edu}}  
 
} 

%% file: title_info/plain_title.tex
\usepackage{datetime}
\newdateformat{monthyeardate}{%
 \monthname[\THEMONTH], \THEYEAR}

\author{
Xu Kuang\\
\texttt{kuangxu@stanford.edu}
\and
Stefan Wager \\
\texttt{swager@stanford.edu}}
\date{Stanford Graduate School of Business}
\title{Weak Signal Asymptotics for \\ Sequentially Randomized Experiments}